\documentclass[draft, reqno,A4paper]{amsart}
\usepackage[english, francais]{babel}
\usepackage{amsfonts}
\usepackage[utf8]{inputenc}
\usepackage[T1]{fontenc}
\usepackage{graphicx}
\usepackage{color}
\usepackage{amsmath}
\usepackage{amssymb}
\usepackage{amsthm}
\usepackage{enumerate}
\usepackage{mathrsfs} 
\usepackage{eqlist}
\usepackage{array}

\theoremstyle{plain}

\newtheorem{thm}{Th\'eor\`eme}[section]
\newtheorem{pro}[thm]{Proposition}
\newtheorem{lem}{Lemme}[section]
\newtheorem{cor}{Corollaire}[thm]

\theoremstyle{definition}
\newtheorem{df}{D\'{e}finition}
\newtheorem{rem}{\textup{Remarque}} 
\newtheorem{exe}{\textit{Exemple}} 
\newtheorem{con}{\textit{Conjecture}}

\numberwithin{equation}{section}

\newcommand{\pre}{\noindent\noindent{\textit{Preuve. }}}
\newcommand{\di}{\displaystyle}
\def\sk{\vskip 0.2cm}
\def\Z{\mathbb{Z}}
\def\N{{\mathbb{N}}}
\def\squareforqed{\hbox{\rlap{$\sqcap$}$\sqcup$}}
\def\qed{\ifmmode\else\unskip\quad\fi\squareforqed}
\def\smartqed{\def\qed{\ifmmode\squareforqed\else{\unskip\nobreak\hfil
\penalty50\hskip1em\null\nobreak\hfil\squareforqed
\parfillskip=0pt\finalhyphendemerits=0\endgraf}\fi}}
\def\smartqed{\def\qed{\ifmmode\squareforqed\else{\unskip\nobreak\hfil
\penalty50\hskip1em\null\nobreak\hfil\squareforqed
\parfillskip=0pt\finalhyphendemerits=0\endgraf}\fi}}
\newcommand{\cqfd}{\textcolor{blue}{\smartqed \qed}}

\newcommand{\al}{\alpha}

\newcommand{\La}{\Lambda}
\newcommand{\be}{\beta}
\newcommand{\ep}{\varepsilon}

\newcommand{\Om}{\Omega}

\newcommand{\te}{\theta}

\begin{document}
\title[Extension de la $i$-modularit\'{e}]{Extension de la $i$-modularit\'{e}}
\author[El Hassane Fliouet]%
{El Hassane Fliouet*}

\newcommand{\acr}{\newline\indent}

\address{\llap{*\,}Regional Center for the Professions of Education and Training,\acr
                   Agadir, Morocco}
\email{fliouet@yahoo.fr}



\selectlanguage{english}
\begin{abstract}
Let $ K / k $ be a purely inseparable extension of characteristic $ p> 0 $ and of finite size. We recall that $K/k$ is modular if for every $n \in \N $,
$K^{p^n}$ and $k$ are $k\cap K^{p^ n}$-linearly disjoint.  A natural generalization of this notion is to say that $K/k$ is $lq$-modular if $K$ is modular over a finite extension of $k$.  Our main objective is to extend in definite form the results and definitions of the $lq$-modularity that have already been obtained in the case limited by the finiteness condition imposed on $[k :k^p]$ in a rather general framework (framework of extensions of finite size  called also $q$-finite extensions).

First, by means of invariants, we characterize the $lq$-modularity of a $q$-finite extension.
Next, we show that any intersection of a $q$-finite extensions covering $k$ or $K$ preserves the $lq$-modularity.
We also prove that any $q$-finite extension $ K/k $ contains a greater $lq$-modular and relatively perfect sub-extension.
In particular, this result is very useful for defining the modularity of order $i$ linked to a $q$-finite extension $ K/k $.
Moreover, we give a necessary and sufficient condition for $K/k$ to be $i$-modular. Certainly, the modularity level of $ K / k $ never exceeds the size
of $K/k$.
Notably, we explicitly describe the extension $K/k $ whose degree of modularity is the size of $ K / k $. In the end, we examine a particular decomposition of $ K/k $ defined by inverse chaining.
\end{abstract}
\maketitle
\selectlanguage{francais}
\section{Introduction}

Soit $K/k$ une extension purement ins\'{e}parable de caract\'{e}ristique $p>0$. Une partie $B$ de $K$ est dite $r$-base de $K/k$ si $K=k(K^p)(B)$ et pour tout $x\in B$, $x\not\in k(K^p)(B\setminus \{x\})$. En vertu du (\cite{N.B}, \textcolor{blue}{III},  p. \textcolor{blue}{49}, corollaire \textcolor{blue}{{3}}) et de la propri\'{e}t\'{e} d'\'{e}change des $r$-bases,  nous d\'{e}duisons que toute extension admet une $r$-base et que le cardinal d'une $r$-base  est invariant. Si de plus $K/k$ est d'exposant fini, on v\'{e}rifie aussit\^{o}t que $B$ est une $r$-base de $K/k$ si et seulement si $B$ est un g\'{e}n\'{e}rateur minimal de $K/k$. Compte tenu de (\cite{N.B}, \textcolor{blue}{III}, p. \textcolor{blue}{25}, proposition \textcolor{blue}{{2}}),  nous pouvons  contr\^{o}ler  la taille de toute extension purement ins\'{e}parable $K/k$ au moyen du degr\'{e} d'irrationalit\'{e} de $K/k$   d\'{e}fini  par $di(K/k)=\di\sup_{n\in\N}(|G_n|)$ o\`{u} $G_n$ est un g\'{e}n\'{e}rateur minimal de $k^{p^{-n}}\cap K/k$.  Notamment,  la mesure de la taille d'une extension est une fonction croissante par rapport \`{a} l'inclusion. Plus pr\'{e}cis\'{e}ment, pour toute chaine d'extensions purement ins\'{e}parables $k\subseteq L\subseteq L' \subseteq K$, on a $di(L'/L)\leq di(K/k)$.  D\'{e}sormais, toute extension de taille finie sera appel\'{e}e extension $q$-finie. Il est clair que les extensions $q$-finies contiennent strictement les extensions de $k$ dont le degr\'{e} $[k :k^p]$ est fini. De plus, on montre que toutes cha\^{i}ne d\'{e}croissante d'extensions $q$-finies est stationnaire.
On rappelle \'{e}galement que $K/k$ est dite modulaire si et seulement si pour tout $n\in\N$, $K^{p^{n}}$ et $k$ sont $K^{p^{n}}\cap k$-lin\'{e}airement disjointes. Cette notion a \'{e}t\'{e} d\'{e}finie pour la premi\`{e}re fois par Swedleer dans {\cite{Swe}}, elle caract\'{e}rise les extensions purement ins\'{e}parables qui sont produit tensoriel sur $k$ d'extensions simples sur $k$. Dans \cite{Wat},  Waterhouse  a montr\'{e} que la modularit\'{e} est stable par une intersection quelconque portant sur $k$ ou $K$, et qu'une r\'{e}union croissante d'extensions modulaires est aussi modulaire. En particulier, il existe une plus petite sous-extension de $K/k$ not\'{e}e   $lm(K/k)$  telle que  $K/lm(K/k)$  est modulaire.
Dans cette note nous continuons \`{a} \'{e}tudier le positionnement de $lm(K/k)$  par rapport \`{a} $k$ et $K$. A cet \'{e}gard, $K/k$ est dite $lq$-modulaire  si $lm(K/k)/k$  est finie. Il s'agit de la modularit\'{e} \`{a} une extension finie pr\`{e}s.
Sachant que \cite{Che-Fli5}, \cite{Che-Fli1}, et \cite{Che-Fli3}  sont enti\`{e}rement consacr\'{e}s \`{a} l'\'{e}tude de cette notion dans  le cas local d\'{e}limit\'{e} par l'hypoth\`{e}se $[k :k^p]$ est fini, dans le pr\'{e}sent  papier, nous d\'{e}sirons \'{e}tendre sous une forme d\'{e}finitive les  r\'{e}sultats et les d\'{e}finitions de la $q$-modularit\'{e} qui ont \'{e}t\'{e} d\'{e}j\`{a} obtenues localement dans un cadre  assez r\'{e}nov\'{e} (c'est le cadre des extensions $q$-finies).
Dans un premier temps, nous traitons les questions de stabilit\'{e}  de la $q$-modularit\'{e} li\'{e}e \`{a} une extension $q$-finie. D'abord, comme dans {\cite{Che-Fli1}} nous commen\c{c}ons par caract\'{e}riser la $lq$-modularit\'{e} d'une extension $q$-finie au moyen d'invariants.  Ensuite,   on montre qu'une intersection quelconque des extensions $q$-finies portant sur $k$ ou $K$ pr\'{e}serve la $lq$-modularit\'{e}. Plus particuli\`{e}rement, il existe des plus petites extensions $k \longrightarrow lqm(K/k) \longrightarrow K \longrightarrow ulqm(K/k)$ telles que $K/lqm(K/k)$ et $ulqm(K/k)/k)$ sont $lq$-modulaires. De plus, nous montrons  que $lqm(K/k)$  est exactement la cl\^{o}ture relativement parfaite de $lm(K/k)/k$, cependant la position de $ulqm(K/k)$ varie selon le choix des exemples, donc on ne peut pas situer avec pr\'{e}cision $ulqm(K/k)$. Dans la seconde \'{e}tape, nous prolongeons le th\'{e}or\`{e}me de la cl\^{o}ture $lq$-modulaire (cf. \cite{Che-Fli3})  \`{a} une extension $q$-finie quelconque. Autrement dit, toute extension $q$-finie $K/k$ admet une plus grande sous-extension relativement parfaite et $lq$-modulaire
que l'on note $H(K/k)$ et que l'on appelle cl\^{o}ture $lq$-modulaire de $K/k$. D'ailleurs, ce r\'{e}sultat  s'av\`{e}re fort utile pour la  construction de la $i$-\`{e}me cl\^{o}ture
$lq$-modulaire d'une extension $q$-finie. Dans cette vue, on note $H_1(K/k)=K$ si $K/k$ est
finie, et $H_1(K/k)=H(K/k)$ si $K/k$ est d'exposant non born\'{e}. Par
r\'{e}currence, on pose $H_i(K/k)=K$ si $K/H_{i-1}(K/k)$ est finie, et
$H_i(K/k)=H(K/H_{i-1}(K/k))$ si $K/H_{i-1}(K/k)$ est d'exposant non
born\'{e}. Comme $K/k$ est $q$-finie, la suite $(H_i(K/k))_{i\geq 1}$
est stationnaire \`{a} partir d'un certain rang $n_0$. En particulier, $k=H_0(K/k)\longrightarrow   H_1(K/k)\longrightarrow\dots\longrightarrow H_{n_0}=K$, et donc toute extension $q$-finie se d\'{e}compose en extensions $lq$-modulaires. Ce qui permet de d\'{e}finir la $lq$-modularit\'{e} d'ordre $i$. A cette occasion, une
extension $q$-finie $K/k$ est dite $i$-modulaire si elle se
d\'{e}compose en $i$ extensions $lq$-modulaires~  ; en d'autres termes, s'il
existe une suite croissante d'extensions $(F_j/k)_{0\leq j\leq i}$ telle que
$F_0=k\subseteq F_1\subseteq \cdots\subseteq F_i=K$, avec $F_{j+1}/F_j$ est $lq$-modulaire
pour tout entier $j \in [0, i[$. Le plus petit entier
$i$ tel que $K/k$ est $i$-modulaire s'appelle le degr\'{e} de
modularit\'{e} de $K/k$, et se note $dm(K/k)$. Il permet de mesurer
le niveau de modularit\'{e} de $K/k$. Il est clair que toute
extension finie est $1$-modulaire, donc ce cas est trivial. Par
contre, si $K/k$ est d'exposant non born\'{e},  les
r\'{e}sultats suivants sont v\'{e}rifi\'{e}es :
\sk

\begin{itemize}{\it
\item[{\rm (1)}] Le niveau de modularit\'{e} de $K/k$ ne d\'{e}passe
jamais l'entier $di(rp(K/k)/k$ o\`{u} $rp(K/k)/k$ est la cl\^{o}ture relativement parfaite de $K/k$.
\item[{\rm (2)}] $K/k$ est $i$-modulaire si et seulement si la
cl\^{o}ture relativement parfaite de $K/k$ que l'on note $rp(K/k)$
co\"{\i}ncide avec celle de $H_i(K/k)/k$. En outre,
$dm(K/k)=\inf\{i\in\N$ tel que $rp(H_i(K/k)/k)=rp(K/k)\}$.
\item[{\rm (3)}] La $i$-modularit\'{e} est stable par le produit.
C'est-\`{a}-dire si $K_1/k$ et $K_2/k$ sont deux sous-extensions
$i$-modulaires de $K/k$, il en est de m\^{e}me de $K_1(K_2)/k$.
\item[{\rm (4)}] La $i$-modularit\'{e} est respect\'{e}e si on change
le corps de base dans le sens ascendant. Plus pr\'{e}cis\'{e}ment, pour toute extension
$L/k$,
on a $K(L)/L$ est $i$-modulaire si $K/k$ l'est.
}
\end{itemize}
\sk

Au plus des r\'{e}sultats cit\'{e}s en haut, on d\'{e}crit explicitement
les extensions $q$-finies $K/k$ dont le niveau de modularit\'{e} atteint la
taille de $K/k$.
\sk

Enfin, nous examinons une d\'{e}composition particuli\`{e}re d'une extension $q$-finie  d\'{e}finie par cha\^{i}nage inverse. D'abord, on montre que toute extension $q$-finie $K/k$ se d\'{e}compose sous forme $K=m_0(K/k)\supseteq m_1(K/k)\supseteq \cdots \supseteq m_{i-1}(K/k)\supseteq m_i(K/k)=k$ o\`u
$m_{j+1}(K/k)=lqm(m_{j}(K/k)/k)$ pour tout entier $j \in [0,i[$. Il est \'{e}videment clair que la d\'{e}composition ainsi obtenue $k=m_i(K/k)\subset m_{i-1}(K/k)\subseteq\cdots\subseteq m_0(K/k)=K$ est une suite de d\'{e}composition associ\'{e}e \`{a} $K/k$. Elle sera appel\'{e}e d\'{e}composition $l$-$i$-modulaire de $K/k$.  Comme dans le cas de la $i$-modularit\'{e}, on s'int\'{e}resse aux questions de stabilit\'{e} et de  pr\'{e}servation de cette notion.
\section{D\'{e}finitions et r\'{e}sultats pr\'{e}liminaires}

Nous commencerons par donner quelques notations que nous utiliserons souvent tout au long de cette note.
\sk

\begin{itemize}{\it
\item $k$ d\'{e}signe toujours un corps commutatif de caract\'{e}ristique $p>0$, et $\Om$ une cl\^{o}ture alg\'{e}brique de $k$.
\item $k^{p^{-\infty}}$ indique la cl\^{o}ture purement ins\'{e}parable de $\Om/k$.
\item Pour tout $a\in \Om$, pour tout $n\in \N^*$, on symbolise la racine du polyn\^{o}me $X^{p^n}-a$ dans $\Om$  par $a^{p^{-n}}$. En outre, on pose $k(a^{p^{-\infty}})=k(a^{p^{-1}},\dots, a^{p^{-n}},$ $\dots )=\di\bigcup_{n\in \N^*} k(a^{p^{-n}})$ et
$k^{p^{-n}}=\{a\in \Om\,|\ , a^{p^{n}}\in k\}$.
\item Pour toute famille $B=(a_i)_{i\in I}$  d'\'{e}l\'{e}ments de $\Om$, on note $k(B^{p^{-\infty}})= k({({a_i}^{p^{-\infty}}}$ $)_{i\in I})$.
\item Enfin, |.| sera employ\'{e} au lieu du terme cardinal.}
\end{itemize}
\sk

 Il est \`{a} signaler aussi que toutes les extensions qui interviennent dans ce papier sont des sous-extensions purement ins\'{e}parables de $\Om$, et il est commode de noter $[k,K]$ l'ensemble des corps interm\'{e}diaires d'une extension $K/k$.
\begin{df}Soit $K/k$ une extension. Une partie $G$ de $K$
est dite $r$-g\'{e}n\'{e}ra\-t\-eur de $K/k$, si $K=k(G)$  ; et si de plus pour
tout $x\in G$, $x\not \in k(G\backslash x)$, $G$ sera appel\'{e}e
$r$-g\'{e}n\'{e}rateur minimal de $K/k$.
\end{df}
\begin{df} Etant donn\'{e}es  une extension $K/k$  de caract\'{e}ristique $p>0$  et
une partie $B$ de $K$. On dit que $B$ est une  $r$-base de $K/k$, si $B$
est un $r$-g\'{e}n\'{e}rateur minimal de $K/k(K^p)$. Dans le m\^{e}me ordre d'id\'{e}es, on dit que $B$ est
$r$-libre sur $k$, si $B$ est une $r$-base de $k(B)/k$  ; dans le cas
contraire $B$ est dite $r$-li\'{e}e sur $k$.
\end{df}
Pour des raisons de coh\'{e}rence th\'{e}matique, toute $r$-base de $k/k^p$ s'appelle $p$-base de $k$, et toute partie d'\'el\'ements de $k$, $r$-libre sur $k^p$ sera appel\'ee $p$-ind\'epend\-a\-n\-te (ou $p$-libre) sur $k^p$. Consid\'{e}rons maintenant
une partie $B$  de $k$. On v\'{e}rifie imm\'{e}diatement que :
\sk

\begin{itemize}{\it
\item[\rm{(1)}] $B$  est $p$-base de $k$ si et seulement si pour tout $n\in \Z$, $B^{p^n}$ l'est \'{e}galement  de $k^{p^n}$.
\item[\rm{(2)}] $B$ est $r$-libre sur $k^p$ si et seulement si pour tout $n\in \Z$, $B^{p^n}$ l'est aussi sur  $k^{p^{n+1}}$.
\item[\rm{(3)}] $B$  est $p$-base de $k$ si et seulement si $B$ est un $r$-g\'{e}n\'{e}rateur minimal de $k/k^p$.
\item[\rm{(4)}] $B$  est $p$-base de $k$ si et seulement si pour tout $n\in \N^*$, $k^{p^{-n}}=\otimes_k (\otimes_k k(a^{p^{-n}}$ $))_{a\in B}$
et  pour tout $a\in B$, $a\not\in k^p$. En particulier, $B$  est $p$-base de $k$ si et seulement si $k^{p^{-\infty}}=\otimes_k (\otimes_k k(a^{p^{-\infty}}))_{a\in B}$ et pour tout $a\in B$, $a\not\in k^p$.}
\end{itemize}
\sk

Il est \`{a} noter que le produit tensoriel  est  utilis\'{e} conform\'{e}ment  \`{a}  (\cite{N.B2}, \textcolor{blue}{III}, p. \textcolor{blue}{42}, d\'{e}finition \textcolor{blue}{50}). Il est vu comme une limite inductive du produit tensoriel d'une famille finie de $k$-alg\`{e}bre.
\sk

Comme dans l'alg\`{e}bre lin\'{e}aire, le th\'{e}or\`{e}me suivant dit th\'{e}or\`{e}me de la $r$-base incompl\`{e}te permet 
de compl\'{e}ter une partie $r$-libre en une $r$-base.
\begin{thm} [\cite{Che1}, Th\'{e}or\`{e}me \textcolor{blue}{2.7}, Th\'{e}or\`{e}me de la $r$-base incompl\`{e}te] {\label{thm1}}  Etant donn\'{e}es  une extension $K/k$ de caract\'{e}ristique $p>0$, et une partie $B$ de $K$,  $r$-libre sur $k(K^p)$. Pour tout $r$-g\'{e}n\'{e}rateur $G$ de $K/k(K^p)$, il existe un sous-ensemble $G_1$ de $G$ tel que $B\cup G_1$ est une $r$-base de $K/k$.
\end{thm}
Voici quelques cons\'{e}quences imm\'{e}diates :
\sk

\begin{itemize}{\it
\item[\rm{(1)}]  De tout $r$-g\'{e}n\'{e}rateur de $K/k(K^p)$ on
peut en extraire  une $r$-base de $K/k$.
\item[\rm{(2)}]  Toute partie $r$-libre sur $k(K^p)$ peut
\^{e}tre compl\'{e}t\'{e}e en une $r$-base de $K/k$. En particulier, toute partie $p$-ind\'{e}pendente sur $k^p$ peut \^{e}tre \'{e}tendue en une $p$-base de $k$.
\item[\rm{(3)}]  Toute extension $K/k$ admet une $r$-base. En outre, tout corps commutatif de caract\'{e}ristique $p>0$ admet une $p$-base.}
\end{itemize}
\sk

Par ailleurs, toutes les $r$-bases d'une m\^{e}me extension ont m\^{e}me cardinal comme le pr\'{e}cise le  r\'{e}sultat suivant :
\begin{thm} [\cite{Che1}, Th\'{e}or\`{e}me \textcolor{blue}{2.8}] {\label{thm2}} Soit $K/k$ une extension de caract\'{e}ristique $p>0$. Si $B_1$
et $B_2$  sont deux $r$-bases de $K/k$, on a $|B_1|=|B_2|$.
\end{thm}

On rappelle que $K$ est  d'exposant fini sur $k$,
s'il existe $e\in \N$ tel que $K^{p^e}\subseteq k$, et le plus petit entier qui satisfait cette relation sera appel\'{e} exposant (ou hauteur) de $K/k$. Certes,
la proposition  suivante  permet de ramener l'\'{e}tude des propri\'{e}t\'{e}s  des $r$-g\'{e}n\'{e}rateurs minimaux des extensions 
d'exposant fini au cas des extensions  de hauteur $1$,  lesquelles sont  plus riches.
\begin{pro} [\cite{Che1}, Proposition \textcolor{blue}{2.12}] {\label{pr5}} Soit $K/k$ une extension purement ins\'{e}parable d'exposant fini. Pour qu'une  partie de $K$ soit une $r$-base de $K/k$ il faut et il suffit que elle soit $r$-g\'{e}n\'{e}rateur minimal de $K/k$.
\end{pro}

Cette proposition et le th\'{e}or\`{e}me ci-dessus s'av\`{e}rent fort utile pour d\'{e}finir la taille d'une extension et la longueur d'un corps commutatif.
Pour cela, consid\'{e}rons une extension purement ins\'{e}parable $K/k$ de caract\'{e}ristique $p>0$, donc pour tout $n\in \N$, $k^{p^{-n}}\cap K/k$ admet un exposant fini et $K=\di\bigcup_{n\in\N}
k^{p^{-n}}\cap K$.
\begin{df}
l'invariant $ di(K/k)=\sup_{n \in \N}(|B_n|)$ sera appel\'{e} degr\'{e} d'irrationalit\'{e} de $K/k$.
\end{df}

Ici le $\sup$ est pris dans le sens de (\cite{N.B}, \textcolor{blue}{III}, p. \textcolor{blue}{25}, Proposition \textcolor{blue}{2}).
Toutefois,  on reprend la notation $di(k)$ du degr\'{e} d'imperfection de $k$  (cf. respectivement \cite{Che-Fli} et \cite{Che-Fli2}) qui sera d\'{e}finie cette fois-ci par $di(k)=di(k/k^p)$. Syst\'{e}matiquement ces deux invariants permet de contr\^{o}ler la taille de $K/k$ et la longueur de $k$. Par ailleurs, la mesure de la taille d'une extension croit en fonction de l'inclusion. Plus pr\'{e}cis\'{e}ment, on a :
\begin{thm} [\cite{Che1}, Th\'{e}or\`{e}me \textcolor{blue}{3.8}] {\label{thm6}} Pour toutes extensions purement ins\'{e}parables $k\subseteq L\subseteq L'\subseteq K$, $di(L'/L)\leq di(K/k)$. En outre, $di(K/k)=\di\sup(di(L/k))_{L\in [k, K]}$.
\end{thm}

Il est clair que la mesure de la taille d'une extension est vue comme une limite inductive de ces corps interm\'{e}diaires. De plus,
nous v\'{e}rifions aussit\^{o}t que toute famille croissant d'extensions purement ins\'{e}parables  $\di(K_n)_{n\in\N}$ satisfait $di(\di\bigcup_{n\in\N}(K_n)/k)= \di\sup_{n\in\N}(di(K_n/k))$, et comme cons\'{e}quence type on a :
\begin{cor}  {\label{thm5}}
Pour toute extension purement ins\'{e}parable $K/k$, on a $di(K/k)$ $\leq di(k)$. En outre, $di(K)\leq di(k)$.
\end{cor}

Consid\'{e}rons maintenant  deux sous-extensions $K_1/k$ et $K_2/k$ d'exposant fini d'une m\^{e}me extension purement ins\'{e}parable  $K/k$. On v\'{e}rifie aussit\^{o}t que si $B_1$ et $B_2$ sont deux $r$-bases respectivement de $K_1/k$ et $K_2/k$, alors $B_1$ et $B_1\cup B_2$ sont  deux $r$-g\'{e}n\'{e}rateurs respectivement de $K_1(K_2)/K_2$ et $K_1(K_2)/k$. En outre, $di(K_1(K_2)/K_2)\leq di(K_1/k)$ et $di(K_1(K_2)/k)\leq di(K_1/k)+di(K_2/k)$. D'une fa\c{c}on  plus pr\'{e}cise, on a :
\begin{pro} [\cite{Che1}, Proposition \textcolor{blue}{3.5}] {\label{pr8}} Sous les conditions ci-dessus, et si de plus $K_1/k$ et $K_2/k$ sont
$k$-lin\'{e}airement disjointes, on a :
\sk

\begin{itemize}{\it
\item[\rm{(i)}] $B_1\cup B_2$ est une $r$-base de $K_1(K_2)/k$.
\item[\rm{(ii)}] $B_1$ est une $r$-base de $K_1(K_2)/K_2$.}
\end{itemize}
\end{pro}

Comme cons\'{e}quence imm\'{e}diate, on a :
\begin{cor}  {\label{cor4}}  Soient $K_1$ et $K_2$ deux corps interm\'{e}diaires d'une m\^{e}me extension purement ins\'{e}parable $\Om/k$. Alors :
\sk

\begin{itemize}{\it
\item[\rm{(i)}] $di(K_1(K_2)/k)\leq di(K_1/k)+di(K_2 /k)$, et il y'a \'{e}galit\'{e} si $K_1$ et $K_2$ sont $k$-lin\'{e}airement disjoints.
\item[\rm{(ii)}] $di(K_1(K_2)/K_2)\leq di(K_1/k)$, et il y'a \'{e}galit\'{e} si $K_1$ et $K_2$ sont $k$-lin\'{e}air\-e\-ment disjoints.}
\end{itemize}
\end{cor}
 \subsection{Extensions relativement parfaites}

Au cours de cette section, on reprend quelques notions et r\'{e}sultats de {\cite{Che-Fli4}}, puisqu'ils sont utilis\'{e}s fr\'{e}quemment ici.
\sk

Un corps $k$ de caract\'{e}ristique $p$ est dit parfait si $k^{p}=k$ ; dans le m\^{e}me ordre d'id\'{e}es,
on dit que $K/k$ est relativement parfaite  si $k(K^{p})=K$. On v\'{e}rifie ais\'{e}ment que :
\sk

\begin{itemize}{\it
\item La relation "\^{e}tre relativement parfaite" est transitive, c'est-\`{a}-dire si $K/L$ et $L/k$ sont relativement parfaites, alors $K/k$ l'est aussi.
\item Si $K/k$ est relativement parfaite, il en est de m\^{e}me de $L(K)/k(L)$.
\item La propri\'{e}t\'{e} "\^{e}tre relativement parfaite" est stable par un produit quelconque portant sur $k$. Autrement dit, pour toute famille $(K_i/k)_{i\in I}$ d'ext\-e\-nsions relativement parfaites, on a alors $\displaystyle \di \prod_{i}^{}K_{i}/k$ est aussi relativement parfaite.}
\end{itemize}
\sk

Par suite, il existe une plus grande sous-extension relativement parfaite de $K/k$ appel\'{e}e cl\^{o}ture relativement parfaite de $K/k$, et se note $rp(K/k)$.
Certes, la propri\'{e}t\'{e} "\^{e}tre relativement parfaite" v\'{e}rifie les relations d'associativit\'{e}-transitivit\'{e} suivantes :
\begin{pro} [\cite{Che-Fli4}, Proposition \textcolor{blue}{5.2}] Soit $L$ un corps interm\'{e}diaire de
$K/k$. Alors
$$
rp(rp(K/L)/k)=rp(K/k) \quad \mbox{ et } \quad
rp(K/rp(L/k))=rp(K/k).
$$
\end{pro}
\begin{cor} Pour tout $L\in
[k,K]$, on a $K/L \hbox{ finie} \Longrightarrow rp(K/k) \subset  L.$
\end{cor}

En particulier, si $K/k$ est relativement parfaite, on a $K/L \hbox{
$finie$ } \Longrightarrow   L=K.$
Sch\'{e}matiquement on a un $trou$
\sk

\[\begin{array}{rcl}
k\longrightarrow &&K ;\\
&\uparrow&\\
&\hbox{$trou$ }&
\end{array}\]

\noindent et ce $trou$ caract\'{e}rise le fait que $K/k$ est
relativement parfaite. En effet,  supposons que
$K/k$ v\'{e}rifie le $trou$ et soit $B$ une $r$-base de $K/k$.
Supposons $B\neq \emptyset$ ;
soit $x\in B$ et $L=k(K^{p})(B\setminus\left \{x\right \})$ ; on a $K/L$
finie, donc $K=L$ ce qui est absurde.
\begin{pro} [\cite{Che-Fli4}, Lemme \textcolor{blue}{1.2}] {\label{arpaa1}} Soit $K/k$ une extension purement ins\'{e}parable telle que $[K :k(K^p)]$ est fini. Alors on a :
\sk

\begin{itemize}{\it
\item[{\rm (i)}] $K$ est relativement parfaite sur une extension finie de $k$.
\item[{\rm (ii)}] La suite d\'{e}croissante  $(k(K^{p^{n}}))_{n \in
\N}$ est stationnaire sur $k(K^{p^{n_{0}}})=rp(K/k)$.}
\end{itemize}
\end{pro}

Comme  cons\'{e}quence de la proposition pr\'{e}c\'{e}dente, on a :
\begin{pro}[{\cite{Che-Fli4}},  Proposition \textcolor{blue}{6.2}] \label{6.2}  Soit $K/k$ une extension purement ins\'{e}parable telle que $[K :k(K^p)]$ est fini. Pour tout $L\in [k,K]$, on a $rp(K/L)=L(rp(K/k)).$
\end{pro}

Le r\'{e}sultat suivant exprime une condition n\'{e}cessaire et suffisant pour que $K/rp(K$ $/k)$ soit finie. Plus pr\'{e}cis\'{e}ment, on a :
\begin{pro}  {\label{arp1}} Soit $K/k$ une extension purement ins\'{e}parable, alors $K/rp(K/k)$ est finie si est seulement il en est de m\^{e}me de $K/k(K^p)$.
\end{pro}
\pre R\'{e}sulte de la proposition \textcolor{blue}{\ref{arpaa1}}. \cqfd
\subsection{Extensions quasi-finies}
\begin{df}  Toute extension de degr\'{e} d'irrationalit\'{e}  fini s'appelle extension $q$-finie.
\end{df}

En d'autres sens, la $q$-finitude est synonyme de la finitude horizontale. Toutefois, la finitude se traduit par la finitude horizontale et verticale, il s'agit de la finitude au point de vue taille et hauteur. Autrement dit, $K/k$ est finie si et seulement si $K/k$ est $q$-finie d'exposant born\'{e}.
\sk

Dans la suite, pour tout $n\in \N$, $k_n$ d\'{e}signe toujours  $k^{p^{-n}}\cap K$. Consid\'{e}rons maintenant   une sous-extension $L/k$ d'une extension $q$-finie $K/k$. On v\'{e}rifie aussit\^{o}t que :
\sk

\begin{itemize}{\it
\item[\rm{(i)}] La $q$-finitude est transitive, en particulier, pour tout $n\in \N$, $K/k(K^{p^n})$ et $k_n/k$ sont finies.
\item[\rm{(ii)}] Il existe $n_0\in \N$, pour tout entier $n\geq n_0$, $di(k_n/k)=di(K/k)$.}
\end{itemize}
\sk

Par ailleurs, nous allons voir quelques applications imm\'{e}diates des propositions \textcolor{blue}{\ref{arpaa1}} et \textcolor{blue}{\ref{arp1}}.
\begin{pro} {\label{pr9}} Soit $K/k$ une extension $q$-finie. La suite $(k(K^{p^n}))_{n\in\N}$ s'arr\-\^{e}\-te sur $rp(K/k)$ \`{a} partir d'un $n_0$. En particulier,  $K/rp(K/k)$ est finie.
\end{pro}

On obtient en particulier le r\'{e}sultat suivant :
\begin{cor}  {\label{cor7}} La cl\^{o}ture relativement parfaite d'une extension $q$-finie $K/k$ n'est pas triviale. Plus pr\'{e}cis\'{e}ment,  $rp(K/k)/k$ est d'exposant non born\'{e} si  $K/k$ l'est.
\end{cor}
\begin{pro}  {\label{pr10}} Pour toute  extension $q$-finie $K/k$, il existe $n\in \N$ tel que $K/k_n$ est relativement parfaite. En outre,  $k_n(rp(K/k))=K$.
\end{pro}

Les sous-extensions relativement parfaites servent comme  nœuds  de liaisons pour la taille d'une extension $q$-finie, c'est-\`{a}-dire :
\begin{pro} [\cite{Che1}, Proposition \textcolor{blue}{4.8}] {\label{pr12}} Pour toute suite de sous-extensions relativement parfaites  $k=K_0 \subseteq K_1\subseteq \dots\subseteq K_n$ d'une extension $q$-finie $K/k$,  on a $di(K/k)=\di \sum_{i=0}^{n-1} di(K_{n+1}/K_n) +di(K/K_n)$.
\end{pro}
Il s'ensuit le r'{e}sultat suivant :
\begin{cor} Soient $k\subseteq L\subseteq K$ des extensions $q$-finies et relativement parfaites sur $k$. Alors $L=K$ si et seulement si $di(L/k)=di(K/k)$.
\end{cor}
\pre Imm\'{e}diat du fait que $di(K/k)=di(K/L)+di(L/k)$, et donc  $di(L/k)=di(K/k)$ \'{e}quivaut \`{a} ce que $di(K/L)=0$, ou encore $K=L$.\cqfd
\sk

Voici une application extr\^{e}mement importante  du th\'{e}or\`{e}me \textcolor{blue}{\ref{thm6}} et la proposition ci-dessus.
\begin{pro}  {\label{pr41}}Toute suite d\'{e}croissante d'une extension $q$-finie est stationnaire.
\end{pro}
\pre  Soient $(K_n/k)_{n\in \N}$ une suite d\'{e}croissante de sous-extensions de $K/k$ et $(F_i/k)_{i\in \N}$ la suite associ\'{e}e \`{a} leurs cl\^{o}tures relativement parfaites. Compte tenu du th\'{e}or\`{e}me \textcolor{blue}{\ref{thm6}}, la suite des entiers $(di(F_i/k))_{i\in \N}$ est  d\'{e}croissante, donc stationnaire \`{a} partir d'un entier $n_0$, ou encore pour tout entier $n\geq n_0$, $F_i=F_{n_0}$. En vertu de la monotonie, pour tout entier $n\geq n_0$,  $[K_{n+1} :F_{n_0}]\leq [K_{n_0} :F_{n_0}]$. Autrement dit, la suite des entiers $([K_{n} :F_{n_0}])_{n\geq n_0}$ est d\'{e}croissante, donc stationnaire \`{a} partir d'un entier $e$, ou encore
pour tout entier $n\geq e$, $[K_{n} :F_{n_0}]=[K_{e} :F_{n_0}]$. Comme pour tout entier $n\geq e$, $K_n\subseteq K_e$, on en d\'{e}duit que $K_n=K_e$ pour tout entier $n\geq e$. 
\subsection{Exposants d'une extension $q$-finie}

Dans ce paragraphe, nous reprenons quelques d\'{e}finitions et notations de base telles qu'elles sont mentionn\'{e}es dans \cite{Che-Fli}.
Tout d'abord, on consid\`{e}re dans un premier temps   une extension purement ins\'{e}parable finie ${K/k}$.
Pour $x\in K$, posons $o( x/k ) = \inf\{m \in \N|\; x^{p^m}\in k \}$
et $o_1(K/k) = \inf\{m\in \N|\; K^{p^{m}}\subset k\}$. Une
$r$-base $B=\{a_{1},a_{2},\dots, a_{n}\}$ de $K/k$ est dite
canoniquement ordonn\'{e}e si pour $j=1,2,\dots,n$, on a
$o(a_{j}/k(a_{1},a_{2},$ $\dots$ $ ,a_{j-1}))=
o_{1}(K/k(a_{1},a_{2},\dots ,a_{j-1})).$ Ainsi, l'entier
$o(a_j/k(a_1,\dots, $ $a_{j-1}))$ d\'{e}fini ci-dessus v\'{e}rifie
$o(a_j/k(a_1,\dots, a_{j-1}))=\inf\{m\in \N|
\; di(k(K^{p^m})/k)\leq j-1\}$ (cf. {\cite{Che-Fli2}}, p.
\textcolor{blue}{138}, lemme \textcolor{blue}{1.3}). On en
d\'{e}duit aussit\^{o}t  le r\'{e}sultat de ({\cite{Pic}}, p.
\textcolor{blue}{90}, satz \textcolor{blue}{14}) qui confirme
l'ind\'{e}p\-e\-n\-d\-a\-n\-ce des entiers $o(a_i/k(a_1,$ $\dots ,
a_{i-1}) )$, $(1\leq i\leq n)$, vis-\`{a}-vis au choix des $r$-bases
canoniquement ordonn\'{e}es $\{a_1,\dots , $ $a_n\}$  de $K/k$. Par
suite, on pose  $o_i(K/k)=o(a_i/k(a_1,$ $\dots , a_{i-1}) )$ si
$1\leq i\leq n$, et $o_i(K/k)=0$ si $i>n$, o\`u $\{a_1,\dots , a_n\}$
est une $r$-base canoniquement ordonn\'{e}e de $K/k$. L'invariant $o_i(K/k)$
ci-dessus s'appelle le $i$-\`{e}me exposant de $K/k$.
\begin{pro} [\cite{Che-Fli2}, Proposition \textcolor{blue}{5.3}] {\label{pr15}}
 Soient $\{\al_1 ,\dots,\al_n\}$
une $r$-base canoniquement ordonn\'{e}e de $K/k$, et $m_j$ le $j$-\`{e}me exposant de $K/k$,  $1\leq j\leq n$. On a :
\sk

\begin{itemize}{\it
\item[{\rm (1)}]  $k(K^{p^{m_j}})=k(\al_{1}^{p^{m_j}},\dots,
\al_{j-1}^{p^{m_j}})$.
\item[{\rm (2)}] Soit $\Lambda_j=\{(i_1,\dots, i_{j-1})$ tel que
$0\leq i_1<p^{m_1-m_j},\dots,0\leq i_{j-1}<p^{m_{j-1}-m_j}\}$,
alors $\{{(\al_1,\dots, \al_{j-1})}^{{p^{m_j}}\xi}$ tel que $\xi\in
\Lambda_j\}$ est une base de $k(K^{p^{m_j}})$ sur $k$.
\item[{\rm (3)}] Soient $n\in\N$ et $j$ le plus grand entier tel que
$m_j>n$. Alors $\{\al_{1}^{p^{n}},\dots, $ $\al_{j}^{p^{n}}\}$ est une
$r$-base canoniquement ordonn\'{e}e de $k(K^{p^n})/k$, et sa liste
des exposants est $(m_1-n,\dots, m_j-n)$}.
\end{itemize}
\end{pro}

Par ailleurs, voici un algorithme qui permet de construire une $r$-base.
\begin{pro} [\cite{Che-Fli2}, Proposition \textcolor{blue}{4.3}, Algorithme de la compl\'{e}tion des r-bases] {\label{pr17}}  Soient K/k une
extension purement ins\'{e}parable finie, $G$ un $r$-g\'{e}n\'{e}rateur de $K/k$, et $\{\al_1,\dots, $ $\al_s\}$
un syst\`{e}me de $K$ tel que pour tout $j\in \{1,\dots, s\}$, $o(\al_j,k(\al_1,\dots,\al_{j-1}))$ $=o_j(K/k)$.
Pour toute suite $\al_{s+1}, \al_{s+2}, \dots,$  d'\'{e}l\'{e}ments de $G$  v\'{e}rifiant
$o(\al_m,$ $k(\al_1,$ $\dots,\al_{m-1}))=\di\sup_{a\in G}(o(a,k(\al_1,\dots,\al_{m-1})))$,
la suite $(\al_i)_{i\in \N^*}$ s'arr\^{e}te sur un plus grand entier $n$  tel que $o(\al_n, k(\al_1,\dots,\al_{n-1}))>0$. En particulier, $\{\al_1,\dots,$ $\al_n\}$ est une r-base canoniquement ordonn\'{e}e de $K/k$.
\end{pro}

De m\^{e}me, la lin\'{e}arit\'{e} disjointe pr\'{e}serve les exposants d'une extension. Il s'agit d'une forme de  stabilit\'{e} des hauteurs.
\begin{pro} [\cite{Che-Fli}, Proposition \textcolor{blue}{7}] {\label{pr16}}
Soient $K_1/k$ et $K_2/k$ deux
sous-extensions purement ins\'{e}p\-a\-rables de $K/k$. $K_1$ et $K_2$
sont $k$-lin\'{e}airement disjointes si et seulement si
$o_j(K_1(K_2)/K_2)=o_j(K_1/k)$ pour tout $ j\in\N$.
\end{pro}

Dans la deuxi\`{e}me \'{e}tape, on consid\`{e}re que $K/k$ est $q$-finie. Pour tout $n\in \N$, notons $k_n=k^{p^{-n}}\cap K$.
En vertu de ({\cite{Che-Fli}},  Proposition \textcolor{blue}{6}),  pour tout $j\in \N^*$, la suite des entiers naturels $(o_j(k_n/k))_{n\geq 1}$ est croissante, et donc
$(o_j(k_n/k))_{n\geq 1}$ converge vers $+\infty$, ou $(o_j(k_n/k))_{n\geq 1}$ est stationnaire \`{a} partir d'un certain rang. Lorsque $(o_j(k_n/k))_{n\geq 1}$ est born\'{e}e,  par construction,  pour tout entier $t\geq j$, $(o_t(k_n/k))_{n\geq 1}$ est aussi born\'{e}e (et donc stationnaire).
\begin{df} Soient $K/k$ une extension $q$-finie  et $j$ un entier naturel non nul. On appelle le $j$-\`{e}me exposant de $K/k$ l'invariant $o_j(K/k)=\di\lim_{n\rightarrow +\infty} (o_j(k_n/k))$.
\end{df}

En particulier, si $K$ est r\'{e}union croissante d'une famille d'extensions $(K_n)_{n\in \N}$, alors pour tout $j\in \N^*$, on a $o_j(K/k)=\di\lim_{n\rightarrow +\infty} (o_j(K_n/k))$.
\begin{lem} [\cite{Che1}, Lemme \textcolor{blue}{4.14}] {\label{lem2}} Soit $K/k$ une extension $q$-finie, alors  $o_s(K/k)$ est fini si et seul\-ement s'il existe un entier naturel $n$ tel que $di(k(K^{p^n})/k)<s$, et on a $o_s(K/k)=\inf\{m\in \N\,|\, di(k(K^{p^m})/k)<s\}$. En particulier, $o_s(K/k)$ est infini si et seulement si pour tout $m\in \N$, $di(k(K^{p^m})/k)\geq s$.
\end{lem}

Le r\'{e}sultat ci-dessous permet de ramener l'\'{e}tude des propri\'{e}t\'{e}s des exposants des extensions $q$-finies aux extensions finies par le biais  des cl\^{o}tures relativement parfaites.
\begin{thm} [\cite{Che1}, Th\'{e}or\`{e}me \textcolor{blue}{4.15}] {\label{thm8}} Soit $K_r/k$ la cl\^{o}ture relativement parfaite de degr\'{e} d'irratio\-n\-alit\'{e} $s$ d'une extension $q$-finie $K/k$, alors on a :
\sk

\begin{itemize}{\it
\item[\rm{(i)}] Pour tout entier $t\in [1, s]$, $o_t(K/k)=+\infty$.
\item[\rm{(ii)}] Pour tout entier $t> s$, $o_t(K/k)=o_{t-s}(K/K_r)$.}
\end{itemize}
En outre, $o_t(K/k)$ est fini si et seulement si $t> s$.
\end{thm}

Voici une liste de cons\'{e}quences imm\'{e}diates :
\begin{pro} [\cite{Che1}, Proposition \textcolor{blue}{4.16}]  {\label{pr18}} Soient $K$ et $L$ deux corps interm\'{e}diaires  d'une extension $q$-finie $M/k$. Pour tout  $j\in \N^* $, on a $o_j(L(K)/L)\leq o_j(K/k)$.
\end{pro}
\begin{pro} [\cite{Che1}, Proposition \textcolor{blue}{4.17}] {\label{pr19}} Etant donn\'{e}es des extensions $q$-finies $k\subseteq L\subseteq L' \subseteq K$. Pour tout  $j\in \N^*$, on a $o_j(L'/L)\leq o_j(K/k).$
\end{pro}

Par ailleurs la taille d'une extension relativement parfaite reste invariante, \`{a} une extension finie pr\`{e}s comme l'indique le r\'{e}sultat suivant :
\begin{pro} [\cite{Che1}, proposition \textcolor{blue}{4.18}] {\label{pr20}}  Etant donn\'{e}e  une sous-extension $K/k$ relativement parfaite   d'une extension $q$-finie $M/k$. Pour toute sous-extension finie $L/k$ de $M/k$, on a $di(L(K)/L)=di(K/k)$.
\end{pro}
\subsection{Extensions modulaires}

On rappelle qu'une extension
$K/k$ est dite modulaire si et seulement si pour tout
$n\in\N$, $K^{p^{n}}$ et $k$ sont $K^{p^{n}}\cap
k$-lin\'{e}airement disjointes. Cette notion a \'{e}t\'{e} d\'{e}finie
pour la premi\`{e}re fois par Swedleer dans {\cite{Swe}}, elle
caract\'{e}rise les extensions purement ins\'{e}parables, qui sont
produit tensoriel sur $k$ d'extensions simples sur $k$. Par ailleurs, toute $r$-base $B$ (si elle existe) d'une extension purement ins\'{e}parable $K/k$ telle que $K\simeq \otimes_k (\otimes_k k(a))_{a\in B}$ sera appel\'{e}e
$r$-base modulaire. En particulier, d'apr\`{e}s le th\'{e}or\`{e}me de Swedleer,  si $K/k$ est d'exposant born\'{e}, il est \'{e}quivalent de dire que :
\sk

\begin{itemize}{\it
\item[\rm{(i)}] $K/k$ admet une $r$-base modulaire.
\item[\rm{(ii)}] $K/k$ est modulaire.}
\end{itemize}
\sk

Soient $m_j$ le $j$-\`{e}me exposant d'une extension purement ins\'{e}parable finie $K/k$ et $\{\al_1 ,\dots,\al_n\}$
une $r$-base canoniquement ordonn\'{e}e de $K/k$,
donc d'apr\`{e}s la proposition $\textcolor{blue}{\ref{pr15}}$, pour tout $j\in \{2,\dots,  n\}$, il existe des constantes uniques
$C_{\ep}\in k$ telles que ${\al_{j}}^{p^{m_j}}=\di\sum_{\ep \in \Lambda_j}{C_{\ep}}{(\al_1,
\dots, \al_{j-1})}^{p^{m_j}\ep}$, o\`u $\Lambda_j=\{(i_1,\dots, i_{j-1})$ tel que
$0\leq i_1<p^{m_1-m_j},\dots,0\leq i_{j-1}<p^{m_{j-1}-m_j}\}$.
Ces relations s'appellent les \'{e}quations de d\'{e}finition de $K/k$.
\sk

Le crit\`{e}re ci-dessous permet de tester la modularit\'{e} d'une extension.
\begin{thm} [\cite{Che-Fli1}, Proposition \textcolor{blue}{10}] {\label{thm9}}
{\bf [Crit\`{e}re de modularit\'{e}]} Sous les notations ci-dessus,
les propri\'{e}t\'{e}s suivantes sont \'{e}quivalentes :
\sk

\begin{itemize}{\it
\item[\rm{(1)}] $K/k$ est modulaire.
\item[\rm{(2)}] Pour toute $r$-base  canoniquement ordonn\'{e}e $\{\al_1,\dots , \al_n\}$
de $K/k$, les $C_{\ep}\in k \cap K^{p^{m_j}}$ pour tout $j\in \{2,\dots,  n\}$.
\item[\rm{(3)}] Il existe une $r$-base canoniquement ordonn\'{e}e  $\{\al_1,\dots , \al_n\}$
de $K/k$ telle que les $C_{\ep}\in k \cap K^{p^{m_j}}$ pour tout $j\in \{2,\dots,  n\}$.}
\end{itemize}
\end{thm}

Le r\'{e}sultat suivant est cons\'{e}quence imm\'{e}diate de la  modularit\'{e}.
\begin{pro}  {\label{apr2}} Soient $m,n\in \Z$ avec $n\geq m$. Si
$K/k$ est modulaire, alors
 $K^{p^{m}}/k^{p^{n}}$  est modulaire.
\end{pro}
\begin{pro} [{\cite{Che-Fli2}},  Proposition \textcolor{blue}{8.4}] {\label{proa1}} Soit $K/k$ une extension purement ins\'{e}parable finie (respectivement, et modulaire),
et soit $L/k$ une sous-extension de $K/k$ (respectivement, et modulaire) avec
$di (L/k) = s$. Si $K^p\subseteq L$,  il existe une r-base canoniquement ordonn\'{e}e (respectivement, et modulaire)
$ (\al_1, \al_2, \dots, \al_n)$ de $K/k$, et $e_1, e_2,\dots, e_s\in \{1,p\}$ tels que $({\al_1}^{e_1}, {\al_2}^{e_2},\dots, {\al_s}^{e_s})$ soit une r-base canoniquement
ordonn\'{e}e (respectivement, et modulaire) de $L/k$. De plus, pour tout $j\in \{1,\dots,  s\}$, on a $o_j (K/k) = o_j (L/k)$, auquel cas $e_j = 1$, ou $o_j (K/k) = o_j (L/k) + 1$, auquel cas
$e_j = p$.
\end{pro}

Le th\'{e}or\`{e}me suivant de  Waterhouse joue un r\^{o}le important
dans l'\'{e}tude des extensions  modulaires (cf. \cite{Wat} Th\'{e}or\`{e}me \textcolor{blue}{1.1}).
\begin{thm}  Soient $(K_{j})_{j\in I}$ une famille
de sous-corps d'un corps commutatif $\Omega$, et $K$ un autre sous-corps de $\Omega $. Si pour tout $j\in I$, $K$ et $K_j$ sont  $K\cap K_j$-lin\'{e}airement
disjoints,
alors $K$ et  $\di\bigcap_{j} K_{j}$ sont  $K\cap (\di\bigcap_{j} K_{j})$-lin\'{e}airement disjoint.
\end{thm}

Comme cons\'{e}quence, la modularit\'{e} est stable par une intersection quelconque portant soit au dessus ou en dessous d'un corps commutatif. Plus pr\'{e}cis\'{e}ment, on a :
\begin{cor} {\label{apr4}} Sous les m\^{e}mes hypoth\`{e}ses du th\'{e}or\`{e}me ci-dessus, on a :
\sk

\begin{itemize}{\it
\item[\rm{(i)}] Si pour tout $j\in I$, $K_{j}/k$ est modulaire, il en est de m\^{e}me de $\di\bigcap_{j} K_{j}/k$.
\item[\rm{(ii)}] Si pour tout $j\in I$, $K/K_j$ est modulaire, il en est de m\^{e}me de $K/\di\bigcap_{j} K_{j}$.}
\end{itemize}
\end{cor}

D'apr\`{e}s le th\'{e}or\`{e}me de Waterhouse, il existe une plus petite sous-extension $m/k$ de $K/k$ (respectivement une plus petite extension $M/K$) telle que $K/m$ (respectivement $M/k$) est modulaire. Dor\'enavant, on note $m=lm(K/k)$ et $M=um(K/k)$. Toutefois, l'extension $um(K/k)$ sera appel\'{e}e cl\^{o}ture modulaire de $K/k$.

Le r\'{e}sultat suivant est bien connu (cf. \cite{Kim}).
\begin{pro}  {\label{pr23}} Soit $K/k$ une extension purement
ins\'{e}parable et modulaire ; soit pour tout $n\in\N$,
$K_n=k(K^{p^n})$. Alors $k_n/k$, $K/k_n$,  $K_n/k$ et $K/K_n$ sont modulaires.
\end{pro}
\section{Extensions $lq$-modulaires}
\subsection{Invariant de la $lq$-modularit\'{e} d'une extension}

D\'{e}sormais, et sauf mention du contraire,  $K/k$ d\'{e}signe une extension $q$-finie d'exposant non born\'{e}, et pour tout $j\in \N$, $k_j= k^{p^{-j}}\cap K$ et $U_{s}^j(K/k)=j-o_s(k_j/k)$ pour tout $s\in \N^*$.
\begin{df}{\label{d1}}
Le premier entier naturel $i_0$ pour lequel la suite
$(U_{s}^{j}(K/k))_{j\in\N}$ est non born\'{e}e s'appelle l'invariant de la $lq$-modularit\'{e} de $K/k$ et se note $Ilqm(K/k)$.
\end{df}

On v\'{e}rifie aussit\^{o}t que $2\leq Ilqm(K/k)$. Par ailleurs, le r\'{e}sultat suivant est une cons\'{e}quence imm\'{e}diate des propositions \textcolor{blue}{\ref{pr18}} et \textcolor{blue}{\ref{pr19}}.
\begin{pro} {\label{prr1}} Pour tout entier naturel non nul $s$,  la suite $(U_{s}^j(K/k))_{j\in\N}$ est croissante.
\end{pro}
\pre
\smartqed
Il est imm\'{e}diat que $o_s(k_n/k)\leq o_s(k_{n+1}/k)
\leq o_s(k_n/k)+1$, car ${k_{n+1}}^p\subseteq k_n$. Donc $n+1-
o_s(k_{n+1}/k)\geq n-o_s(K_n/k)$ ; et par suite
$(U_{s}^j(K/k))_{j\in\N}$ est croissante. 
\sk

Il en r\'{e}sulte  aussit\^{o}t que :
\sk

\begin{itemize}
\item[$\bullet$] Pour tout entier  $s\geq Ilqm(K/k)$, $\di\lim_{n \rightarrow
+\infty}(U_{s}^j(K/k))=+\infty$.
\item[$\bullet$]  Pour tout entier $s\in [1,Ilqm(K/k)[$,  la suite
$(U_{s}^j(K/k))_{j\in\N}$ est born\'{e}e ; et par suite, pour tout entier
$n\geq \di\sup_{j\in\N} (\di\sup(U_{s}^{j}(K/k)))$ ({s<Ilqm(K/k)}), on
obtient $U_{s}^n(K/k)=U_{s}^{n+1}(K/k)$. Ou encore,
$o_s(k_{n+1}/k)=o_s(k_n/k)+1$.
\end{itemize}
\sk

Cela conduit \`{a} :
\begin{pro} {\label{pr2}} Soient $k\subseteq K_1\subseteq K_2$ des extensions $q$-finies d'exposant non born\'{e}. Pour tout $s\in \N^*$, pour tout $n\in \N$,
 on a $U_{s}^n(K_1/k)\geq U_{s}^n(K_2/k)$ . En outre, $Ilqm(K_1/k)\leq
Ilqm(K_2/k)$, et il y'a \'{e}galit\'{e} si $K_2/K_1$ est finie.
\end{pro}
\pre
\smartqed
En vertu de la proposition \textcolor{blue}{\ref{pr19}}, pour tout $n\in \N$, on a
$k^{p^{-n}}\cap K_1\subseteq k^{p^{-n}}\cap K_2$. Par passage aux
exposants, on obtient $o_s(k^{p^{-n}}\cap K_1/k)\leq
o_s(k^{p^{p^{-n}}}\cap K_2/k)$, et donc $U_{s}^{n}(K_2/k)\leq
U_{s}^{n}(K_1/k)$. Il en r\'{e}sulte que $Ilqm(K_1/k)\leq Ilqm
(K_2/k)$.
\sk

Dans le cas o\`{u} $K_2/K_1$ est finie, posons $e=o_1(K_2/K_1)$, et
donc pour tout entier $n>e$,  $k(k^{p^{-n}}\cap K_{2}^{p^{e}})\subseteq
k^{p^{-n}}\cap K_1$. Compte tenu de des propositions \textcolor{blue}{\ref{pr19}} et \textcolor{blue}{\ref{pr15}}, on obtient
$o_s(k^{p^{-n-e}}\cap K_2/k)-e\leq o_s(k^{p^{-n}}\cap K_1/k)$. Soit
$e+n-o_s(k^{p^{-n-e}}\cap K_2/k)\geq n-o_s(k^{p^{-n}}\cap K_1/k)$, ou encore $U_{s}^{n}(K_1/k)\leq
U_{s}^{n+e}(K_2/k)$.
D'o\`{u} $Ilqm(K_2/k)\leq Ilqm(K_1/k)$ ; et par suite,
$Ilqm(K_1/k)=Ilqm(K_2/k)$.\cqfd
\subsection{Caract\'{e}risation de la $lq$-modularit\'{e}  d'une extension $q$-finie}
\begin{df}{\label{d2}}  Soit $K/k$ une extension $q$-finie. $K/k$ est dite $lq$-modulaire si
$K$ est modulaire sur une extension finie de $k$.
\end{df}

D'apr\`{e}s le th\'{e}or\`{e}me de Waterhouse, il existe une plus
petite sous-extension $m/k$ de $K/k$ telle que $K/m$ est modulaire.
Par suite, $K/k$ est $lq$-modulaire si et seulement si $m/k$  est
finie.
\sk

Le r\'{e}sultat qui suit caract\'{e}rise les extensions
$lq$-modulaires au moyen de la variation des exposants de certains
corps interm\'{e}diaires. Plus pr\'{e}cis\'{e}ment, on a :
\begin{thm} {\label{thm10}}  Soit $K/k$ une extension $q$-finie et $t=di(rp(K/k)/k)$. Les affirmations suivantes sont \'{e}quivalentes :
\sk

\begin{itemize}{\it
\item[\rm{(1)}] $K/k$ est $lq$-modulaire.
\item[\rm{(2)}] Il existe un entier naturel $j$ tel que $K/k^{p^{-j}}\cap
K$ est modulaire.
\item[\rm{(3)}]  Pour tout entier $s\in[1,t]$,  la suite  $(U_{s}^{j}(K/k))_{j\in\N}$  est born\'{e}e.
\item[\rm{(4)}] $Ilqm(K/k)=t+1$.}
\end{itemize}
\end{thm}
\pre
Il est clair que $(2) \Leftrightarrow (3)$. Par ailleurs, compte tenu de la proposition \textcolor{blue}{\ref{pr10}}, il existe un entier naturel non nul  $j_0$ tel que $K/k_{j_0}$ est relativement parfaite et $k_{j_0}(rp(K/k))=K$. En particulier, d'apr\`{e}s la proposition \textcolor{blue}{\ref{pr20}}, $di(K/k_{j_0})= di(rp(K/k)/k)=t$. Supposons ensuite que la condition $(1)$ est v\'{e}rifi\'{e}e. On distingue alors deux cas :
\sk

 Si $K/k$ est modulaire, en vertu de (\cite{Che1}, proposition \textcolor{blue}{6.3}), pour tout entier $j\geq j_0$, on a $k_j/k_{j_0}$ est  \'{e}quiexponentielle d'exposant $j-j_0$ et $di(k_j/k_{j_0})=t$. D'o\`{u} pour tout $s\in \{1,\dots,t\}$, on a $U_{s}^{j}(K/k)=U_{s}^{j+1}(K/k)$.
 \sk

Si $K$ est modulaire sur une extension finie $L$ de $k$,  compte tenu de la finitude de  $L/k$, il existe un entier naturel $e_1$ tel que $L\subseteq k_{e_1}$. Par suite, $L^{p^{-j}}\cap K\subseteq k_{e_1+j}$, et par passage aux exposants pour tout $s\in \N^*$, $o_s(L^{p^{-j}}\cap K)\leq o_s(k_{e_1+j})$ ; soit  donc $U_{s}^{e_1+j}(K/k)\leq e_1+U_{s}^{j}(K/L)$. D'o\`{u}, la suite $(U_{s}^{j}(K/k))_{j\in \N}$ est stationnaire pour tout $s\in\{1,\dots, t\}$.
\sk

 Inversement, si la  condition $(2)$ est v\'{e}rifi\'{e}e, il existe $m_0\geq \sup(e(K/k),j_0)$,  pour tout entier $j\geq m_0$, pour tout $s\in \{1,\dots,t\}$,  on a $o_s(k_{j+1}/k)=o_s(k_j/k)+1$ (et $di(k_j/k_{m_0})=t$). Par suite, $k_j/k_{j_0}$ est \'{e}quiexponentielle, donc modulaire. D'o\`{u} $K=\di\bigcup_{j>m_0} k_j$ est modulaire sur $k_{j_0}$. \cqfd
\sk

Du fait que la suite $(U_{s}^{j}(K/k))_{s\in \N^*}$ est croissante ($j$ \'{e}tant un entier naturel fixe), la condition $(3)$ du th\'{e}or\`{e}me ci-dessus  se r\'{e}duit  \`{a}  $(U_{t}^{j}(K/k))_{j\in \N}$ est born\'{e}e, et par suite $K/k$ est $lq$-modulaire si et seulement si  la suite $(U_{t}^{j}(K/k))_{j\in \N}$ est born\'{e}e.
\sk

Sans perdre de g\'{e}n\'{e}ralit\'{e}, le r\'{e}sultat qui suit
am\'{e}liore naturellement les hypoth\`{e}ses de la condition
suffisante du th\'{e}or\`{e}me ci-dessus.
\begin{pro} {\label{pr3}} Soit $(K_n/k)_{n\in\N}$ une suite croissante de sous-extensions finies d'une extension $q$-finie $L/k$ et $K=\displaystyle\bigcup_{n\in\N}K_n$.
Si les conditions suivantes sont v\'{e}rifi\'{e}es :
\sk

\begin{itemize}
\item [{\rm {(i)}}] La suite $(o_1(K_n/k))_{n\in\N}$ est non born\'{e}e,
\item[{\rm {(ii)}}] Pour tout entier $s\in [1,t]$ o\`{u} $t=di(rp(K/k)/k)$,
la suite $(o_1(K_n/k)-o_s(K_n/k))_{n\in\N}$ est born\'{e}e ;
\end{itemize}
alors $K/k$ est $lq$-modulaire.
\end{pro}
\pre Tout d'abord, on pose $e_n=o_1(K_n/k)$, notamment $K_n\subseteq
k^{p^{-e_n}}\cap K$. Par passage aux exposants, on obtient
$o_t(K_n/k)\leq o_t(k^{p^{-e_n}}\cap K/k)$, et donc
$e_n-o_t(K_n/k)\geq e_n-o_t(k^{p^{-e_n}}\cap K/k)$, ou encore $e_n-o_t(K_n/k)\geq U_{t}^{e_n}(K/k)$. Compte tenu de la monotonie des suites $(e_n)_{n\in\mathbb{N}}$ et $(U_{t}^{j}(K/k))_{j\in\mathbb{N}}$, et vue que $(e_n)_{n\in\mathbb{N}}$ est non born\'{e}e, on en d\'{e}duit que $\di\sup_{j\in \N}(U_{t}^{j}(K/k))=\di\sup_{j\in \N}(U_{t}^{e_j}(K/k))=\di\sup_{j\in \N}(e_j-o_t(K_j/k))$, et par cons\'{e}quent la
suite $(U_{t}^{j}(K/k))_{j\in\mathbb{N}}$ est born\'{e}e.  D'apr\`{e}s le th\'{e}or\`{e}me pr\'{e}c\'{e}dent, $K/k$ est $lq$-modulaire.\cqfd
\sk

La condition (2) est n\'{e}cessaire, mais non suffisante comme le montre l'exemple suivant :
\begin{exe} Soient $P$ un corps parfait de caract\'{e}ristique
$p>0$, et $k=P(X,Y)$ le corps des fractions rationnelles aux
ind\'{e}termin\'{e}es $X,Y$. Posons
$K=\displaystyle\bigcup_{n\in
\N}K_n$ o\`{u} $K_n=k(X^{p^{-2n}}, Y^{p^{-n}})$. On a $K/k$ est modulaire (donc
$lq$-modulaire), mais la suite $(o_1(K_n/k)-o_2(K_n/k))_{n\in\N}$ est non born\'{e}e.
\end{exe}

Soient $k$ un corps commutatif de caract\'eristique $p>0$ et $\Om$ une cl\^{o}ture alg\'{e}brique de  $k$. Dans  $[k :\Om]$ on d\'{e}finit la relation $\sim$ de la fa\c{c}on suivante : $k_1\sim k_2$ si et seulement si $k_1\subseteq k_2$ et $k_2/k_1$ est finie ou $k_2\subseteq k_1$ et $k_1/k_2$ est finie.
On v\'{e}rifie aussit\^{o}t que $\sim$ est r\'{e}flexive, sym\'{e}trique, cependant $\sim$ est g\'{e}n\'{e}ralement non transitive.
De plus, pour toute extension $q$-finie $K_1/k$, l'application de modularit\'{e} inf\'{e}rieure :
\sk

\begin{eqnarray*}
lm : [k :K_1] & \longmapsto& [k :K_1]\\
L & \longrightarrow & lm(K_1/L),
\end{eqnarray*}
\noindent est  compatibles avec la relation $\sim$. Plus pr\'{e}cis\'{e}ment, on a :
\begin{pro} {\label{pr43}} Soient $k_1\subseteq k_2\subseteq K_1$  des sous-extensions  $q$-finies.
Si $k_1\sim k_2$, alors $lm(K_1/k_1)\sim lm(K_1/k_2)$.
\end{pro}
\pre Il suffit de remarquer que  $lm(K_1/k_1)\subseteq lm(K_1/k_2)$, et si  de plus  $o_1(k_2/k_1)$ $=e_1$, alors $k_2 \subseteq (lm(K_1/{k_1))}^{p^{-e_1}}\cap K_1$, avec  $K_1/{(lm(K_1/k_1))}^{p^{-e_1}}\cap K_1$ est modulaire, (cf.  la proposition \textcolor{blue}{\ref{pr23}}). Soit donc $lm(K_1/k_2)\subseteq (lm(K_1/{k_1))}^{p^{-e_1}}\cap K_1$.\cqfd
\sk

Comme cons\'{e}quence,  la $lq$-modularit\'e est stable \`a une extension  finie pr\`es du choix du corps de base comme le pr\'{e}cise le r\'{e}sultat qui suit :
\begin{pro} {\label{p34}} Soit $K/k$ une extension $q$-finie. On a :
\sk

\begin{itemize}
\item[\rm{(1)}] Si $k'\sim k$ et
$k'\subset K$, $K/k$ est $lq$-modulaire si et seulement il en est de
m\^{e}me de $K/k'$.
\item[\rm{(2)}] Si $K\sim K'$ et $k\subset
K'$, $K/k$ est $lq$-modulaire si et seulement si $K'/k$ l'est aussi.
\item[\rm{(3)}] Si $k'\sim k$ et $K\sim K'$, avec $k'\subset
K'$, alors $K/k$ est $lq$-modulaire si et seulement si il en est de
m\^{e}me de $K'/k'$.
\end{itemize}
\end{pro}
\pre \smartqed Cela r\'{e}sulte imm\'{e}diatement des propositions
\textcolor{blue}{\ref{pr2}} et \textcolor{blue}{\ref{pr43}}  et du th\'{e}or\`{e}me \textcolor{blue}{\ref{thm10}}. 
\sk

Comme cons\'{e}quence, le r\'{e}sultat ci-dessous permet de ramener l'\'{e}tude de la
$lq$-modularit\'{e} au cas des extensions relativement parfaites.
\begin{cor} {\label{c35}}
Soient $K/k$ une extension $q$-finie
et $H/k$ la cl\^{o}ture relativement parfaite de $K/k$. Alors :
\sk

\begin{itemize}
\item[{\rm (i)}] $K/k$ est $lq$-modulaire si et seulement si il en est de m\^{e}me de
$H/k$.
\item[{\rm (ii)}] Soit $F/k$ une sous-extension de $K/k$. $K/F$ est
$lq$-modulaire si et seulement si il en est de m\^{e}me de
$rp(K/k)/rp(F/k)$ et de $K/rp(F/k)$.
\end{itemize}
\end{cor}
\pre \smartqed Il suffit  de remarquer que $K\sim rp(K/k)$
et $F\sim rp(F/k)$. \cqfd \sk

Soient $K/k$ une extension purement ins\'{e}parable et $P$ le
sous-corps parfait maximal de $k$. On a $K/P$ est modulaire ; donc
on peut esp\'{e}rer que $K/k$ est $lq$-modulaire! Autrement dit,
toute extension purement ins\'{e}parable serait $lq$-modulaire.
Cependant $K/P$ est transcendante. Voici un exemple non \'{e}vident
d'extensions purement ins\'{e}parables et non $lq$-modulaires.
\begin{exe} {\label{e1}} Soient $k_0$ un corps parfait de caract\'{e}ristique
$p>0$ et $(X,Z_1,Z_2)$ une famille alg\'{e}briquement
ind\'{e}pendante sur $k_{0}$. Notons $k=k_0(X,Z_1,Z_2)$, et pour tout $n\in \N^*$,
$K_{n}=k(X^{p^{-2n}},$ $\theta _{n})$, avec $\theta _{1}=Z_1^{p^{-1}}X^{p^{-2}}+Z_2^{p^{-1}}$, et pour tout entier $n\geq 2$,
\[\theta _{n}=Z_1^{p^{-1}}X^{p^{-2n}}+{(\theta _{n-1})}^{p^{-1}}=Z_1^{p^{-1}}X^{p^{-2n}}+Z_{1}^{p^{-2}}X^{p^{-2n+1}}+
\cdots+Z_1^{p^{-n}}X^{p^{-n-1}}+Z_2^{p^{-n}}.\]
\end{exe}

On a $\theta _{n+1}^{p}=Z_{1}X^{p^{-2n-1}}+\theta _{n}$, donc
$K_{n}\subset K_{n+1}$. Soit  $K=\displaystyle\bigcup_n K_{n}$.
\begin{thm} Sous les hypoth\`{e}ses ci-dessus, $k({X}^{p^{-\infty}})/k$
est la plus petite sous-exten\-s\-ion $m/k$ de $K/k$ telle que $K/m$
est modulaire ($lm(K/k)=k({X}^{p^{-\infty}})$).
\end{thm}

Pour la preuve de ce th\'{e}or\`{e}me, on  se servira du r\'{e}sultat
suivant qui est une cons\'{e}quence imm\'{e}diate du crit\`{e}re de la
modularit\'{e}. Il intervient fr\'{e}quemment dans le reste de ce
travail.
\begin{lem} {\label{l37}} Soient $K/k$ une extension modulaire de caract\'{e}ristique
$p\not = 0$ et $((a,b),(e_1,e_2))\in k^2\times K^2$ tels que
${e_2}^{p^j}=a{e_{1}}^{p^{j }}+b$ ($j\in \N$). Si  ${e_{1}}^{p^j}\not\in k$, alors ${a}^{p^{-j}}\in K$ et
${b}^{p^{-j}}\in K$.
\end{lem}
\pre \smartqed ${e_{1}}^{p^j}\not\in k$ est synonyme de
$(1,{e_{1}}^{p^j})$ est libre sur $k$. Notamment, $(1,{e_{1}}^{p^j})$ est libre sur $k\cap K^{p^j}$. Prolongeons ensuite ce
syst\`{e}me en une base $B$ de $K^{p^{j}}$ sur $k\cap K^{p^j}$.
Comme $K^{p^j}$ et $k$ sont $k\cap K^{p^j}$-lin\'{e}airement
disjointes (car $K/k$ est modulaire), $B$ est aussi une base de
$k(K^{p^j})$ sur $k$ ; et comme ${e_2}^{p^j}=a{e_{1}}^{p^{j}}+b$ ;
par identification, on a $a\in k\cap K^{p^{j}}$ et $b\in k\cap
K^{p^j}$. Il en r\'{e}sulte que ${a}^{p^{-j}}\in K$ et
${b}^{p^{-j}}\in K$. \cqfd
\sk

\noindent\noindent{\bf Preuve du th\'{e}or\`{e}me} \smartqed Soit $m/k$ la plus petite
sous-extension de $K/k$ telle que $K/m$ est modulaire. Il est clair que
$K/k(X^{p^{-\infty}})$ est modulaire, donc $m\subseteq k(X^{p^{-\infty}})$. Supposons
l'existence d'un entier naturel non nul $n$ tel que $X^{p^{-n+1}}\in m$ et
$X^{p^{-n}}\not\in m$. Par construction, on a  ${\theta_{1}}^p=Z_1{(X^{p^{-2}})}^p+Z_2=Z_1X^{p^{-1}}+Z_2$ et,
$$
{[{Z_1}^{p^{-1}}X^
{p^{-2n}}+\cdots+{Z_1}^{p^{-n}}X^{p^{-n-1}}+{Z_2}^{p^{-n}}]}^{p^n}={Z_1}^{p^
{n-1}}X^{p^{-n}}+\cdots + {Z_1}X^{p^{-1}}+Z_2
$$
pour tout entier $n\geq 2$, d'apr\`{e}s le lemme ci-dessus, on a  ${Z_1}^{p^{-1}} \in K$. Il en
r\'{e}sulte que ${Z_2}^{p^{-1}}\in K$,  puisque
${Z_1}^{p^{-1}}X^{p^{-2}}+{Z_2}^{p^{-1}}\in K$ ; et par suite, $
di(k(X^{p^{-1}},{Z_1}^{p^{-1}},{Z_2}^{p^{-1}})/k)$ $=3 \leq di(K/k)=2$,
ce qui est absurde. \cqfd
\sk

Dans ce qui suit, nous illustrons ces propos par quelques exemples d'extensions $lq$-modulaires extraits de
\cite{Che-Fli1}. On rappelle  qu'une extension
$K/k$ est dite $w_0$-g\'{e}n\'{e}r\'{e}e, si pour toute sous-extension propre $L$ de $K/k$, on a $L/k$
est finie (cf. \cite{Che-Fli2} et \cite{Dev2} pour des exemples non
triviaux d'extensions $w_0$-g\'{e}n\'{e}r\'{e}es).
\begin{pro} Toute extension
$w_0$-g\'{e}n\'{e}r\'{e}e et $q$-finie est
$lq$-modulaire.
\end{pro}
\pre \smartqed Imm\'{e}diat. \cqfd
\sk

Une extension alg\'{e}brique $K/k$ est dite simple si $K=k(\theta)$
$(\theta \in K)$. Dans le cas purement ins\'{e}parable, cette d\'{e}finition se traduit en terme de taille et de hauteur par $K/k$ est simple si et seulement si $K/k$ est d'exposant fini et de degr\'{e} d'irrationalit\'{e} $di(K/k)=1$.  Une g\'{e}n\'{e}ralisation naturelle de cette notion est de dire que $K/k$ est $q$-simple si l'ensemble des corps interm\'{e}diaires de $K/k$ est totalement ordonn\'{e} (cf. \cite{Che-Fli2}) ; dans le cas purement ins\'{e}parable, cela \'{e}quivaut \`{a} ce que toute sous-extension propre de $K/k$ est simple ou encore $di(K/k)=1$.
\sk

Soient $F/k$ et $L/k$ deux sous-extensions d'une extension $\Omega/k$. On
v\'{e}rifie imm\'{e}d\-i\-a\-t\-e\-m\-ent en se ramenant \`{a} des
extensions simples, que nous avons les propri\'{e}t\'{e}s suivantes :
\sk

\begin{itemize}{\it
\item[$\bullet$] Si $F/k$ est $q$-simple, alors $L(F)/L$ est $q$-simple.
\item[$\bullet$]  Si $F/k$ est $q$-simple, alors $L(F)=L\otimes  _{L\cap F}F$.}
\end{itemize}
\sk

Une extension $K/k$ est dite semi-simple de genre fini si elle est produit fini d'extensions $q$-simples (Pour plus d'informations \`{a} propos des extensions semi-simples il est sugg\'{e}r\'{e} de se r\'{e}f\'{e}rer au \cite{Che-Fli7}).
\begin{pro} Toute extension   semi-simple de genre fini $K/k$ est produit tensoriel
d'e\-x\-t\-ensions $q$-simples sur une extension finie de  $k$.
\end{pro}
\pre \smartqed Application imm\'{e}diate de la transitivit\'{e} et
l'associativit\'{e} du produit tensoriel. \cqfd
\sk

Comme cons\'{e}quence, on a :
\begin{cor} {\label{c310}}
Toute extension semi-simple de genre fini  est
$lq$-modulaire.
\end{cor}
\pre \smartqed Imm\'{e}diat. \cqfd
\sk

Dans la suite, on montre comme dans le cas de la modularit\'{e} que
la $lq$-modularit\'{e} est pr\'{e}serv\'{e}e par intersection.
\subsection{Stabilit\'{e} de la $lq$-modularit\'{e} d'une extension $q$-finie}

D'une fa\c{c}on assez g\'{e}n\'{e}rale, la $lq$-modularit\'{e} est respect\'{e}e non seulement \`{a} une extension finie pr\`{e}s, mais \`{a} toute d\'{e}placement du corps de base dans le sens ascendant. Il s'agit d'une propri\'{e}t\'{e} absolue de la $lq$-modularit\'{e}.
\begin{pro} {\label{le2}}
Soient $k\subseteq L \subseteq K$ des extensions $q$-finies. Si
$K/k$ est $lq$-modulaire,  il en est de m\^{e}me de $K/L$.
\end{pro}
\pre D'apr\`{e}s la proposition \textcolor{blue}{\ref{p34}}, il suffit de faire la d\'{e}monstration  lorsque  $K/k$ et $L/k$ sont relativement parfaites.
D'abord, compte tenu de la propri\'{e}t\'{e} "\^{e}tre relativement parfaite",  pour tout $n\in \N$, $k(L^{p^n})=L$. Or, $k(L^{p^{n}})\subset k(L\cap K^{p^{n}})\subset L$, donc $L=k(L\cap K^{p^n})$, et par suite il est trivialement  \'{e}vident  que $k(K^{p^{n}})=K$ et $L$ sont $k(L\cap K^{p^{n}})$-lin\'{e}airement
disjoints. En vertu de (\cite{Mor-Vin}, lemme \textcolor{blue}{1.60}),  si $K/k$ est modulaire, il en est de m\^{e}me de  $K/L$. Sinon, par la $lq$-modularit\'{e} il existe une sous-extension
finie $k_1/k$ de $K/k$ telle que $K/k_1$ est modulaire, donc il en est de m\^{e}me de $K/k_1(L)$ puisque $rp(k_1(L)/k_1)=k_1(rp(L/k))=k_1(L)$, et comme $k_1(L)\sim L$, alors $K/L$ sera $lq$-modulaire.
\cqfd
\begin{pro} {\label{p316}}
Soient $K/k$ une extension $q$-finie et $(k_{i})_{i\in I}$ une famille de corps interm\'{e}diaires de
$K/k$. Si pour tout $i\in I$, $K/k_{i}$ est $lq$-modulaire, il en est
de m\^{e}me de $K/\displaystyle\bigcap_{i\in I}^{}k_{i}$.
\end{pro}
\pre \smartqed Gr\^{a}ce \`{a} la proposition \textcolor{blue}{\ref{pr41}},  on se ram\`{e}ne \`{a} $I=\{1,2\}$.  Pour tout $j=1,2$, on note  $m_j=lm(K/k_j)$. Compte tenu de la $lq$-modularit\'{e}, il existe $e\in \N$ tel que $m_j\subseteq {k_j}^{p^{-e}}\cap K$ pour $j=1,2$ ; et donc $m_1\cap m_2\subseteq {k_1}^{p^{-e}}\cap {k_2}^{p^{-e}}\cap K={(k_1\cap k_2)}^{p^{-e}}\cap K$. Il en r\'{e}sulte que $m_1\cap m_2/k$ est finie (\`{a} savoir $K/k$ est $q$-finie). Par ailleurs, en vertu du corollaire \textcolor{blue}{\ref{apr4}},  $K/m_1\cap m_2$ est modulaire, donc $K/k_1\cap k_2$ est $lq$-modulaire. \cqfd
\sk

Comme cons\'{e}quence, on a :
\begin{cor} Pour toute extension $q$-finie $K/k$, il existe une plus petite
sous-extension $m/k$ de $K/k$ telle que $K/m$ est $lq$-modulaire.
\end{cor}
\pre \smartqed Imm\'{e}diat. \cqfd
\sk

On peut situer avec pr\'{e}cision $m$. Pour cela, comme dans le cas de la modularit\'{e}, on commence  par d\'{e}signer par
$lqm(K/k)$ la plus petite sous-extension $m/k$ de $K/k$ telle que
$K/m$ est $lq$-modulaire.
\begin{pro} {\label{p318}} Sous les conditions ci-dessus, on a  $lqm(K/k)=rp(lm(K/k)/k)$.
En outre, $lqm(K/k)/k$ est relativement parfaite.
\end{pro}
\pre Soient $m_0=lqm(K/k)$, $m_1=lm(K/m_0)$ et $m_2=lm(K/k)$. Il est
clair que $m_0\subseteq m_2\subseteq m_1$ et $m_1/m_0$ est finie,
puisque $K/m_0$ est $lq$-modulaire. D'o\`u $rp(m_0/k)$ $=rp(m_2/k)$.
Or, d'apr\`{e}s le corollaire \textcolor{blue}{\ref{c35}},
$K/rp(m_0/k)$ est aussi $lq$-modulaire, donc $m_0\subseteq rp(m_0/k)$, ou
encore $m_0=rp(m_0/k)=rp(m_2/k)$. \cqfd \sk

Cela conduit \`{a} :
\begin{cor} {\label{c319}}
Soient $K$ et $K'$ deux corps interm\'{e}diaires d'une extension
$q$-finie $K/k$. Alors :
\sk

\begin{itemize}{\it
\item[{\rm (i)}] Si $K\sim K'$,  on a $lqm(K/k)=lqm(K'/k)$. En particulier, $lqm(K/k)=lqm(rp(K/k)/k)$.
\item[{\rm (ii)}] Pour toute sous-extension $L/k$ de $K/k$, on a $L(lqm(K/k))=lqm(K/L)$.}
\end{itemize}
\end{cor}
\pre \smartqed Compte tenu  de la proposition
\textcolor{blue}{\ref{p34}}, il suffit de prouver l'assertion $(ii)$.
Notons $m_1=lqm(K/k)$ et $m_2=lqm(K/L)$. Il est clair que
$m_1\subseteq m_2$, donc $m_1(L)\subseteq m_2$. D'autre part,
d'apr\`{e}s la proposition \textcolor{blue}{\ref{le2}}, $K/m_1(L)$
est $lq$-modulaire, donc $m_2\subseteq m_1(L)$. D'o\`{u}
$m_1(L)=m_2$. \cqfd
\sk

La pr\'{e}servation de la $lq$-modularit\'{e} de $K/k$ par
intersection portant sur $K$ semble beaucoup moins \'{e}vidente. D'abord, nous
aurons besoin des r\'{e}sultats  suivants.
Nous commen\c{c}ons par la proposition  suivante qui est une cons\'{e}quence bien connue de la transitivit\'{e} de la lin\'{e}arit\'{e} disjointe.
\begin{pro} {\label{pr7}} Soient $K_1/k$ et $K_2/k$ deux sous-extensions d'une m\^eme extension $K/k$, $k$-lin\'{e}airement disjointes.  Pour touts corps interm\'{e}diaires  $L_1$ et $L_2$  respectivement de $K_1$ et $K_2$, on a $L_2(K_1)$ et $L_1(K_1)$ sont $k(L_1, L_2)$-lin\'{e}airement-disjointes. En particulier, $L_2(K_1)\cap L_1(K_2)=k(L_1,L_2)$.
\end{pro}

Cela conduit \`{a} :
\begin{cor} {\label{gcor1}} Soient $L$ et $K$ deux corps interm\'{e}diaires  d'une m\^{e}me extension $M/k$, $k$-lin\'{e}airement disjoints. Tout  sous-ensemble $G$  de $K$ tel que $L(K)=L(G)$ v\'{e}rifie $K=k(G)$. En outre, si une partie $G$ de $K$ est une base de $L(K)/L$, alors $G$ l'est \'{e}galement pour  $K/k$.
\end{cor}
\pre
D'apr\`{e}s la transitivit\'{e} de la lin\'{e}arit\'{e} disjointe, on aura $K$ et $L(G)$ sont $k(G)$-lin\'{e}airement disjoints, et par suite $k(G)=K\cap L(G)=K\cap L(K)=K$. \cqfd\sk

Consid\'{e}rons  maintenant  deux sous-extensions finies $L/k$ et $m/k$  d'une m\^{e}me extension
purement ins\'{e}parable $K/k$. Soit $B=(\alpha
_{1},\alpha _{2},\dots,\alpha _{n})$ une $r$-base canoniquement
ordonn\'{e}e de $L/k$. Posons $e_{i}=o_{i}(L/k)$ et $e=o_{1}(m/k)$.
On suppose qu'il existe un entier $s\in [1,n]$ tel que
$e_{s-1}-e_{s}> e$ (en particulier, $e_{s}\neq e_{s-1}$). D'apr\`{e}s la proposition \textcolor{blue}{\ref{pr15}},
il existe des constantes uniques $a_{\xi}\in k(L^{p^{e_{s}+1}})$ telles que $$\alpha
_{s}^{p^{e_{s}}}=\di\sum_{\xi\in \La }^{}a_{\xi }(\alpha _{1},\alpha
_{2},\dots,\alpha _{s-1} )^{\xi p^{e_{s}}},$$
o\`{u} $\La=\left \{(i_{1},i_{2},\dots,i_{s-1})|\;0\leq
i_{j}\leq p-1\right \}$. Supposons de plus que
$K^{p^{e_{s}}}$ et $k(m\cap K^{p^{e_{s}}})$ sont $m\cap
K^{p^{e_{s}}}$-lin\'{e}airement disjointes. En outre, cette
condition est remplie si $K/m$ est modulaire.
\begin{lem} {\label{l320}} Sous les hypoth\`{e}ses ci-dessus, on a :
\[\forall \xi\in \left \{(i_{1},i_{2},\dots,i_{s-1})|\; 0\leq i_{j}\leq
p-1\right \}, \quad a_{\xi}^{p^{-e_{s}}} \in k^{p^{-e_{1}+1}}\cap
K.\]
\end{lem}
\pre \smartqed Montrons d'abord que ${((\alpha _{1},\alpha
_{2},\dots,\alpha _{s-1} )^{\xi p^{e_{s}}})}_{\xi\in \La}$ est une
base de $k(m\cap
K^{p^{e_s}})({\al_1}^{p^{e_s}},\cdots,{\al_{s-1}}^{p^{e_s}})$ sur
$k(m\cap
K^{p^{e_s}})({\al_1}^{p^{e_s+1}},\cdots,{\al_{s-1}}^{p^{e_s+1}})$.
Pour cela, il suffit de montrer que $di(k(m\cap
K^{p^{e_s}})({\al_1}^{p^{e_s}},\cdots,{\al_{s-1}}^{p^{e_s}})/k(m\cap K^{p^{e_s}}))=s-1$. S'il existe $i\in \{1,\cdots,s-1\}$
tel que ${\al_i}^{p^{e_s}}\in k(m\cap
K^{p^{e_s}})({\al_1}^{p^{e_s}},\cdots,$
${\al_{i-1}}^{p^{e_s}},{\al_{i+1}}^{p^{e_s}} ,$
$\cdots,{\al_{s-1}}^{p^{e_s}})$, alors
${\al_i}^{p^{e_s+e}}={({\al_i}^{p^{e_s}})}^{p^e}\in k(m^{p^e}\cap
K^{p^{e_s+e}})({\al_1}^{p^{e_s+e}},$
$\cdots,{\al_{i-1}}^{p^{e_s+e}},$ ${\al_{i+1}}^{p^{e_s+e}} ,$
$\cdots,{\al_{s-1}}^{p^{e_s+e}}$ $)\subseteq
k(m^{p^e})({\al_1}^{p^{e_s+e}},\cdots,$ ${\al_{i-1}}^{p^{e_s+e}},$
${\al_{i+1}}^{p^{e_s+e}},$ $\cdots,{\al_{s-1}}^{p^{e_s+e}})$
$=k({\al_1}^{p^{e_s+e}},\cdots,{\al_{i-1}}^{p^{e_s+e}},$ $
{\al_{i+1}}^{p^{e_s+e}} ,\cdots,{\al_{s-1}}^{p^{e_s+e}})$, et donc $di(k({\al_1}^{p^{e_s+e}},\cdots,$ ${\al_{i-1}}^{p^{e_s+e}},{\al_{i+1}}^{p^{e_s+e}} ,\cdots,$
${\al_{s-1}}^{p^{e_s+e}})/k)< s-1$. Il en
r\'{e}sulte  d'apr\`{e}s le lemme \textcolor{blue}{\ref{lem2}} que  $o_{s-1}(L/k)=o_{s-1}(k(\al_1,$ $\cdots,\al_{s-1})/k)=e_{s-1}\leq e_s+e$, ou
encore $ e_{s-1}-e_s\leq e $, ce qui contredit l'hypoth\`{e}se.
Par ailleurs, on a $K^{p^{e_s}}$ et $k(m\cap K^{p^{e_s}})$ sont
$m\cap K^{p^{e_s}}$-lin\'{e}airement disjointes, donc  en particulier $m\cap K^{p^{e_s}}(L^{p^{e_s}})$ et $k(m\cap K^{p^{e_s}})$ sont
$m\cap K^{p^{e_s}}$-lin\'{e}airement disjointes. Comme, $k(L^{p^{e_s}})=k({\al_1}^{p^{e_s}},\cdots,{\al_{s-1}}^{p^{e_s}})$, donc plus particuli\`{e}rement
$k(m\cap K^{p^{e_s}})(L^{p^{e_s}})=k(m\cap K^{p^{e_s}})({\al_1}^{p^{e_s}},\cdots,{\al_{s-1}}^{p^{e_s}})$. D'apr\`{e}s le corollaire pr\'{e}c\'{e}dent $m\cap K^{p^{e_s}}({L}^{p^{e_s}})=m\cap K^{p^{e_s}}({\al_1}^{p^{e_s}},\cdots,{\al_{s-1}}^{p^{e_s}})$, et par suite ${\al_s}^{p^{e_s}}\in   m\cap K^{p^{e_s}}({\al_1}^{p^{e_s}},\cdots,{\al_{s-1}}^{p^{e_s}})$ (\`{a} savoir ${\al_s}^{p^{e_s}}\in  {L}^{p^{e_s}}$).
D'autre part,  en vertu de la
transitivit\'{e} de la lin\'{e}arit\'{e} disjointe, $m\cap K^{p^{e_s}}({L}^{p^{e_s}})$  et
$k(m\cap K^{p^{e_s}} )(L^{p^{e_s+1}})$ sont $m\cap
K^{p^{e_s}}({L}^{p^{e_s+1}})$-lin\'{e}airement disjointes.
Puisque $B_1= {((\alpha_{1},\alpha _{2},\dots,\alpha _{s-1} )^{\xi
p^{e_{s}}})}_{\xi\in\La}$ est une base de $k(m\cap
K^{p^{e_s}})({L}^{p^{e_s}})$
sur $k(m\cap K^{p^{e_s}})({L}^{p^{e_s+1}})$,  alors $B_1$ est aussi une base de
$m\cap K^{p^{e_s}}({L}^{p^{e_s}})$ sur $m\cap
K^{p^{e_s}}({L}^{p^{e_s+1}})$.
Or, $$\alpha _{s}^{p^{e_{s}}}=\di\sum_{\xi }^{}a_{\xi }(\alpha
_{1},\alpha _{2},\dots,\alpha _{s-1} )^{\xi p^{e_{s}}}$$  avec les
$a_{\xi}\in k(L^{p^{e_{s}+1}})\subseteq k(m\cap K^{p^{e_s}}
)(L^{p^{e_s+1}})$ ; par identification, pour tout $\xi\in \La$,
$a_{\xi}\in m\cap K^{p^{e_s}}({L}^{p^{e_s+1}})\subseteq K^{p^{e_s}}$. Il en r\'{e}sulte que
${a_{\xi}}^{p^{-e_s}}\in K$. On a aussi $a_{\xi}\in
k(L^{p^{e_s+1}})\subseteq k^{p^{-e_1+e_s+1}}$ (car
$o_1(k(L^{p^{e_s+1}})/k)=e_1-e_s-1$). D'o\`u $a_{\xi}\in
k^{p^{-e_1+1}}$, et donc $a_{\xi}\in k^{p^{-e_1+1}}\cap
K$. \cqfd
\begin{thm} Soient $(K_{i}/k)_{i\in I}$ une famille de
sous-extensions $lq$-modulaires d'une extension $q$-finie $H /k$. Alors $\di\bigcap_{i\in I}K_{i}/k$ est $lq$-modulaire.
\end{thm}
\pre \smartqed Gr\^ace \`{a} la proposition \textcolor{blue}{\ref{pr41}}, on se ram\`{e}ne \`{a}
$I=\{1,2\}$. On a $K_1/k$ et $K_2/k$ sont $lq$-modulaires, donc
d'apr\`{e}s le th\'{e}or\`{e}me \textcolor{blue}{\ref{thm10}} et la proposition
\textcolor{blue}{\ref{pr23}}, il existe $e\in \N$ tel que $K_1/k^{p^{-e}}\cap
K_1$ et $K_2/k^{p^{-e}}\cap K_2$ sont modulaires.
Le cas $q$-simple et le cas fini sont trivialement \'{e}vidents, donc il suffit d'\'{e}tablir le r\'{e}sultat lorsque
$K_1\cap K_2/k$ est d'exposant non born\'{e} et de degr\'{e} d'irrationalit\'{e}  $1<di(K_1\cap K_2/k)$, et par cons\'{e}quent $1\leq di(rp(K_1\cap K_2)/k)$ et
$o_1(k^{p^{-j}}\cap (K_1\cap K_2)/k)=j$ pour tout $j\in \N$. Posons ensuite,  $t=di(rp(K_1\cap K_2)/k)$, et pour tout entier
$s \in [1, t]$, $U_{s}^j(K_1\cap K_2/k)=j-o_s(k^{p^{-j}}\cap
K_1\cap K_2/k)$, $e_{s}^j=o_s(k^{p^{-j}}\cap K_1\cap K_2/k)$,
$\ep_{s}^j=e_{s}^{j+1}-e_{s}^j$ ($\ep_{s}^j\in \{0,1\}$), et $Ilqm(K_1\cap K_2/k)=i_0$.
Si $i_0=t+1$, en vertu du th\'{e}or\`{e}me \textcolor{blue}{\ref{thm10}},  $K_1\cap K_2/k$ est $lq$-modulaire. Supposons ensuite que $1<i_0\leq t$. Par d\'{e}finition l'invariant de la $lq$-modularit\'{e} d'une extension, on a
$\di\lim_{j\longrightarrow +\infty}(U_{i_0}^j(K_1\cap
K_2/k))=+\infty$, et pour tout entier $s\in [1, i_0-1]$, la suite
croissante d'entiers $(U_{s}^j(K_1\cap K_2/k))_{j\in \N}$ est born\'{e}e, donc stationnaire. D'o\`u il existe
$e_1\in \N$, pour tout entier $j\geq e_1$, pour tout entier $s\in [1, i_0-1]$,
$U_{s}^{j}(K_1\cap K_2/k)=U_{s}^{j+1}(K_1\cap K_2/k)$. En particulier, $\di\lim_{j\longrightarrow +\infty}(U_{i_0}^j(K_1\cap
K_2/k)-U_{i_0-1}^j(K_1\cap K_2/k))=\di\lim_{j\longrightarrow
+\infty}(e_{i_0-1}^j-e_{i_0}^j)=+\infty$. Par suite, il existe
$n_0>e_1$, pour tout entier $n\geq n_0$, $e_{i_0-1}^n-e_{i_0}^n>>>e$. Si
$\ep_{i_0}^n=1$ pour tout entier $n\geq n_0$, alors la suite
$(U_{i_0}^n(K_1\cap K_2/k))_{n\in \N}$ est born\'{e}e, c'est une contradiction. Donc il existe $n_1\geq n_0$ tel que $\ep_{i_0}^{n_1}=0$,
ou encore $e_{i_0}^{n_1}=e_{i_0}^{n_1+1}$. Soit $n_2$ le plus grand
entier tel que $e_{i_0}^{n_2}=e_{i_0}^{n_1}$ ($n_2$ existe, car $1<i_0\leq t$, et donc par le th\'{e}or\`{e}me \textcolor{blue}{\ref{thm8}}
$\di\lim_{n \longrightarrow +\infty}(e_{i_0}^n)=+\infty$). D'apr\`{e}s
la proposition \textcolor{blue}{\ref{proa1}}, il existe une $r$-base
$\{\al_1,\cdots,\al_{m}\}$ canoniquement ordonn\'{e}e de
$k^{p^{-n_2-1}}\cap K_1\cap K_2/k$, il existe $\ep_{i_0+1}\in
\{1,p\},\cdots, \ep_{m'}\in \{1,p\}$, ($m'=di(k^{p^{-n_2}}\cap K_1\cap K_2/k)$),
tels que $\{{\al_1}^p,\cdots,
{\al_{i_0}}^p, {\al_{i_0+1}}^{\ep_{i_0+1}},\cdots,
{\al_{m'}}^{\ep_{m'}}\}$ est aussi une $r$-base canoniquement
ordonn\'{e}e de $k^{p^{-n_2}}\cap K_1\cap K_2/k$. Ecrivons :
$${\al_{i_0}}^{p^{e_{i_0}^{n_2+1}}}=\di\sum_{\xi}a_{\xi}{({\al_{1}},
\cdots,{\al_{i_0-1}})}^{\xi p^{e_{i_0}^{n_2+1}}},$$ o\`{u} $a_{\xi}\in
k({\al_{1}}^{p^{e_{i_0}^{n_2+1}+1}},\cdots,{\al_{i_0-1}}^{p^{e_{i_0}^{n_2+1}+1}})$
et $\La=\{(i_{1},i_{2},\dots,i_{i_0-1})|\; 0\leq i_j \leq p-1 \}$.
Comme $K_1/k^{p^{-e}}\cap K_1$ et $K_2/k^{p^{-e}}\cap K_2$ sont
modulaires,  en vertu du lemme \textcolor{blue}{\ref{l320}}
ci-dessus, ${a_{\xi}}^{p^{-e_{i_0}^{n_2+1}}}\in
k^{p^{-n_2}}\cap K_1$ et ${a_{\xi}}^{p^{-e_{i_0}^{n_2+1}}}\in
k^{p^{-n_2}}\cap K_2$. D'o\`u $\al_{i_0}\in
k(\al_1,\cdots,$ $\al_{i_0-1},((a_{\xi})^{p^{-e_{i_0}^{n_2+1}}})_{\xi})$, et donc $e_{i_0}^{n_2}=e_{i_0}^{n_2+1}-1=o({\al_{i_0}}^p,
k({\al_1}^p,$ $\cdots,{\al_{i_0-1}}^p))$ $\geq
o_1(k({\al_1}^p,\cdots,{\al_{i_0-1}}^p,((a_{\xi})^{p^{-e_{i_0}^{n_2+1}}})_{\xi})/k({\al_1}^p,\cdots,{\al_{i_0-1}}^p))\geq o_1(k(\al_1,$
$\cdots,$ $\al_{i_0-1},$ $((a_{\xi})^{p^{-e_{i_0}^{n_2+1}}})_{\xi})/
k(\al_1,\cdots,\al_{i_0-1}))\geq o(\al_{i_0}, k(\al_1,
\cdots,\al_{i_0-1}))= e_{i_0}^{n_2+1}$ (cf. Propositions \textcolor{blue}{\ref{pr18}} et
\textcolor{blue}{\ref{pr19}}). Par suite, $e_{i_0}^{n_2+1}-1\geq
e_{i_0}^{n_2+1}$, ce qui est absurde. \cqfd
\sk

Nous  en d\'{e}duisons le r\'{e}sultat suivant :
\begin{cor} Pour toute extension $q$-finie $K/k$,
il existe une plus petite extension $M/K$ telle que $M/k$ est
$lq$-modulaire.
\end{cor}
\pre  \smartqed Imm\'{e}diat. \cqfd
\sk

D\'{e}sormais, $ulqm(K/k)$ d\'{e}signe la plus petite
extension $M/K$ telle que $M/k$ est $lq$-modulaire. Malheureusement, on ne sait pas situer avec pr\'{e}cision $ulqm(K/k)$.
Tout ce que l'on peut affirmer est que $ulqm(K/k)\subset rp(um(K/k)/k)$. Mais l'inclusion peut \^{e}tre stricte comme
l'indique l'exemple suivant :
\begin{exe}  Soient $P$ un corps parfait de caract\'{e}ristique
$p\neq 0$ et $(X,Y,Z)$ une famille alg\'{e}briquement
 libre sur $P$. Notons $k=P(X,Y,Z)$, $F_{1}=k(X^{p^{-\infty }})$,
$F_{2}=k((X^{p^{-2}}Y^{p^{-1}}$ $+Z^{p^{-1}})^{p^{-\infty }})$ et
$K=F_{1}(F_{2})$. Posons $\theta_1 =Y^{p^{-1}}X^{p^{-2}}+Z^{p^{-1}}$
et $\theta_2 =X^{p^{-2}}$.
\end{exe}
\begin{pro}
Avec les notations pr\'{e}c\'{e}dentes, on a $um(K/k)/ulqm(K/k)$ est d'exposant non born\'{e} (et donc $ulqm(K/k)\not=um(K/k)$).
\end{pro}
\pre \smartqed Comme $F_{1}/k$ et  $F_{2}/k$ sont $q$-simples, selon
le corollaire \textcolor{blue}{\ref{c310}}, $K/k$ est $lq$-modulaire ; d'o\`u
$ulqm(K/k)=K$. Soit $N=um(K/k)$.  Si $N/K$ est finie, comme $K/k$ est relativement parfaite, donc $K=rp(N/k)=\di\bigcap_{n\in \N} k(N^{p^n})$, et par suite $K/k$ sera modulaire. Or,
${\theta_1} ^{p}=Y{\theta_2}^{p}+Z$ ; en vertu du lemme
\textcolor{blue}{\ref{l37}}, $Y^{p^{-1}},Z^{p^{-1}}\in K$ ; et par
suite, $2=di(K/k)\geq di(k(X^{p^{-2}},Y^{p^{-1}},Z^{p^{-1}})/k)=3$,
contradiction. D'o\`u $N/ulqm(K/k)$ est d'exposant non born\'{e}, et donc $N\not=K=ulqm(K/k)$. \cqfd
\section{Prolongement du th\'{e}or\`{e}me de la cl\^{o}ture $lq$-modulaire}
\subsection{Cl\^{o}ture $lq$-modulaire d'une extension $q$-finie}

Soit $K/k$ une extension $q$-finie d'exposant non born\'{e}. Dans ce qui suit, on utilise les notations suivantes : $k_n=k^{p^{-n}}\cap K$,  $o_s(k_n/k)=e_{s}^{n}$, $di(K/k)=t_1$, $Ilqm(K/k)=i_0$, et
$e(K/k)=\di\sup_{s\in I}(\di\sup_{j\in\N}(U_{s}^{j}(K/k)))$ o\`u
$I=\{1,\cdots, i_0-1\}$.
\begin{thm} {\label{th2}} Pour tout entier $n\geq e$, on a $k(k_{n}^{p^{e_{i_0}^{n}}})
\subseteq k(k_{n+1}^{p^{e_{i_0}^{n+1}}})$, et l'extension
$M=\displaystyle\bigcup_{n\geq e}k(k_{n}^{p^{e_{i_0}^n}})$
remplie les conditions ci-dessous :
\sk

\begin{itemize}{\it
\item[{\rm (1)}] $di(M/k)=Ilqm(K/k)-1=i_0-1$.
\item[{\rm (2)}] $M/k$ est relativement parfaite.
\item[{\rm (3)}] $M/k$ est $lq$-modulaire.}
\end{itemize}
\end{thm}
\pre
\smartqed
D'abord, pour tout
entier $j\geq e(K/k)$,  pour tout entier $s\in [1, i_0-1]$, on a
$e_{s}^{j+1}=e_{s}^j+1$, et comme $k(k_{j+1}^p)\subseteq k_j$,
d'apr\`{e}s la proposition \textcolor{blue}{\ref{proa1}}, pour tout entier $j\geq e(K/k)$, il existe une
$r$-base $\{\al_1,\cdots, \al_{t_1}\}$ canoniquement ordonn\'{e}e de
$k_{j+1}/k$, il existe $e_s\in \{1,p\}$ tels que
$\{\al_{1}^{p},\cdots, \al_{i_0-1}^p,$ $ \al_{i_0}^{e_{i_0}},\cdots,
\al_{t_1}^{e_{t_1}}\}$ est une $r$-base canoniquement ordonn\'{e}e
de $k_j/k$.
Ensuite, pour tout entier $j>e(K/k)$, notons $K_j = k(k_{j}^{p^{e^{j}_{i_0}}})$. D'une part,  $K_j=k(\al_{1}^{p^{e^{j}_{i_0}+1}},\dots,\al_{i-1}^{p^{e^{j}_{i_0}+1}})$ et $K_{j+1}= k(\al_{1}^{p^{e^{j+1}_{i_0}}},\dots,\al_{i-1}^{p^{e^{j+1}_{i_0}}})$. D'autre part, on a $e^{j+1}_{i_0}=e^{j}_{i_0}+\ep $, avec $\ep=0$ ou $1$, cela conduit \`a $K_j\subseteq K_{j+1}$. Toutefois,  par d\'{e}finition de $Ilqm(K/k)$, on a $1+ e^{j}_{i_0}> e^{j+1}_{i_0}$ (c'est-\`a-dire $e^{j}_{i_0} = e^{j+1}_{i_0}$) pour une infinit\'{e} de valeurs de $j$. Pour ces valeurs, on a $di(K_{j+1}/k) = i_0-1$, sinon d'apr\`es le lemme \textcolor{blue}{\ref{lem2}},  $e^{j+1}_{i_0} = e^{j+1}_{i_0-1} = 1 + e^{j}_{i_0-1} = e^{j}_{i_0}$ , et donc $e^{j}_{i_0} > e^{j}_{i_0-1}$, ce qui contredit la d\'efinition des exposants.
Comme $(di(K_j/k))_{j>e(K/k)}$ est une suite croissante d'entiers born\'{e}e par $di(K/k)$, donc elle stationne sur $i_0-1$. De plus, $K_j \not= K_{j+1}$, en effet si $K_{j+1}=K_j =k( K^{p}_{j+1})$,  comme $K_{j+1}/k$ est d'exposant born\'e, on aura $K_{j+1} = k$, contradiction.
Cela entra\^{i}ne aussit\^ot que $M/k$ est d'exposant non born\'{e} et $di(M/k) = i_0 -1$. De plus, $M/k$ est relativement parfaite car $k(K^{p}_{j+1})= K_j$ pour une infinit\'{e} de $j$.
\cqfd
\subsection{Th\'{e}or\`{e}me d'existence}
\begin{thm}[th\'{e}or\`{e}me de la cl\^{o}ture $lq$-modulaire] L'extension $M/k$
ci-de\-s\-s\-us est la plus grande sous-extension $lq$-modulaire et
relativement parfaite de $K/k$.
\end{thm}

Pour d\'{e}montrer ce th\'{e}or\`{e}me, nous aurons besoin des
terminologies suivantes :
\begin{df} {\label{d3}} Une suite $(e_n)_{n\in\N}$ d'\'{e}l\'{e}ments de $\N$ est
 dite $q$-enti\`{e}re
si pour tout $n\in\N$, $(e_{n+1}-e_n)\in \{0,1\}$.
\end{df}
\begin{lem} {\label{le1}} Soient $K/k$ une extension $lq$-modulaire et relativement parfaite,
et $(e_n)_{n\in\N}$ une suite $q$-enti\`{e}re telle que $\di
\lim_{n\rightarrow +\infty}(n-e_n)=+\infty$. Alors, pour tout entier $n\geq
e(K/k)$, on a :
\sk

\begin{itemize}{\it
\item[{\rm (i)}] $k(k_{n}^{p^{e_n}}) \subseteq k(k_{n+1}^{p^{e_{n+1}}})$.
\item[{\rm (ii)}] $\displaystyle\bigcup_{n\geq e(K/k)}k(k_{n}^{p^{e_n}})=K$.}
\end{itemize}
\end{lem}
\pre
\smartqed
Compte tenu de la $lq$-modularit\'{e} et de la propri\'{e}t\'{e} "\^{e}tre relativement parfaite", pour
tout $n\geq e(K/k)$, pour tout $s\in\{1,\cdots,t\}$ o\`{u}
$t=di(K/k)$,  on a $e_{s}^{n}+1=e_{s}^{n+1}$ ;
c'est-\`{a}-dire $k(k_{n+1}^{p})=k_n$. Comme $e_{n+1}-e_n\in \{0,1\}$, alors
$k(k_{n}^{p^{e_{n}}}) \subseteq k(k_{n+1}^{p^{e_{n+1}}})$.
Ensuite,  pour tout entier $j>e(K/k)$, posons $H=\displaystyle\bigcup_{n\geq e(K/k)} k(k_{n}^{p^{e_{n}}})$. Comme la suite $(U_{t}^{n}(K/k))_{n\in \N}$ est born\'{e}e (cf. th\'{e}or\`{e}me \textcolor{blue}{\ref{thm10}}) et $\di\lim_{n\rightarrow +\infty}(n-e_n)=+\infty$, alors $\di\lim_{n\rightarrow +\infty}((n-e_n)-(U_{t}^{n}(K/k)))= \di\lim_{n\rightarrow +\infty}((n-e_n)-(n-e_{t}^n))=\di\lim_{n\rightarrow +\infty}(e_{t}^n-e_n)=+\infty$. Or, en vertu de la proposition \textcolor{blue}{\ref{pr15}}, $(n-e_n,e_{2}^{n}-e_n,\cdots, e_{t}^{n}-e_n)$ est
exactement la liste des exposants de $k(k_{n}^{p^{e_n}})/k$, et donc d'apr\`{e}s le th\'{e}or\`{e}me \textcolor{blue}{\ref{thm8}}, $H/k$ est relativement parfaite et $di(H/k)=t=di(K/k)$. Par suite $H=K$ puisque $K/k$ est aussi relativement parfaite.\cqfd
\sk

\noindent \noindent {\bf Preuve du th\'{e}or\`{e}me de la cl\^{o}ture $lq$-modulaire.}
\smartqed  Soit $J/k$ une sous-extension
$lq$-modulaire et relativement parfaite de $K/k$. Posons
$J_n=k^{p^{-n}}\cap J$. Comme $\di\lim_{n\rightarrow
+\infty}(n-e_{i_0}^n)=+\infty$, et $(e_{i_0}^n)_{n\in\N}$ est
$q$-enti\`{e}re, en  vertu du lemme pr\'{e}c\'{e}dent, on aura $J=\displaystyle\bigcup_{n\in\N}k({J_n}^{p^{e_{i_0}^n}})$. Par ailleurs,
$k({J_{n}}^{p^{e_{i_0}^n}})\subseteq k(k_{n}^{p^{e_{i_0}^n}})$ ; donc
$J\subseteq M$.\cqfd
\sk

Une \'{e}tude d\'{e}taill\'{e}e autour de ce th\'{e}or\`{e}me dans le cas o\`{u} $di(k)$ est fini se trouve dans
\cite{Che-Fli3}.
\begin{df} Soit $K/k$ une extension $q$-finie.
La plus grande sous-extension $lq$-modulaire et relativement
parfaite de $K/k$ s'appelle la cl\^{o}ture
$lq$-modulaire de $K/k$, et se note $H(K/k)$.
\end{df}
\begin{rem} Soit $K/k$ une extension $q$-finie. Alors :
\sk

\begin{itemize}{\it
\item  Si $K/k$ est finie, $k$ est la cl\^{o}ture $lq$-modulaire  de
$K/k$ ; donc ce cas est trivial. Cependant, si $K/k$ est d'exposant
non born\'{e}, il en est de  m\^{e}me de $H(K/k)/k$. Pour cela, on
s'int\'{e}resse uniquement aux extensions d'exposant non born\'{e}.
\item Pour toute sous-extension $L/k$ de $K/k$, on a $rp(L/k)\subseteq H(K/k)$ si $L/k$
est $lq$-modulaire.
\item $H(K/k)=H(rp(K/k)/k)$ ; en particulier, si $K/k$ est $lq$-modulaire,
$H(K/k)$ est la cl\^{o}ture relativement parfaite de $K/k$.
\item La cl\^{o}ture $lq$-modulaire d'une extension purement ins\'{e}parable peut ne pas \^{e}tre triviale
comme le montre l'exemple ci-dessous.}
\end{itemize}
\end{rem}
\begin{exe} Soient $Q$ un corps  parfait de caract\'{e}ristique
$p>0$, et $k=Q(X,Z_1,Z_2)$ le corps des fractions rationnelles aux
ind\'{e}termin\'{e}es $X,Z_1,Z_2$. Posons $K=\displaystyle\bigcup_{n\in\N^*}K_n$
avec $K_n=k(X^{p^{-2n}},{Z_1}^{p^{-1}}{X}^{p^{-2n}}+\cdots+{Z_1}^
{p^{-n}}X^{p^{-n-1}}+{Z_2}^{p^{-n}})$. L'extension $L=\displaystyle\bigcup_{n\in {\N}^*}k(X^{p^{-n}})$ est la
plus petite sous-extension de $K/k$ telle que $K/L$ est modulaire
(cf. exemple \textcolor{blue}{\ref{e1}}), donc $K/k$ n'est pas $lq$-modulaire.
Il en r\'{e}sulte que $L=H(K/k)$, (cf. th\'{e}or\`{e}me \textcolor{blue}{\ref{th4}} ci-apr\`{e}s).
\end{exe}
\begin{pro} {\label{p313}}
Soit $K/k$ une extension $q$-finie. Alors :
\sk

\begin{itemize}{\it
\item[{\rm (i)}] Pour toute sous-extension finie $k'/k$ de $K/k$, on a $k'(H(K/k))=H(K/k')$.
\item[{\rm (ii)}] Pour toute extension finie $K'/K$, on a $H(K/k)=H(K'/k)$.}
\end{itemize}
\end{pro}
\pre \smartqed L'assertion $(ii)$ est imm\'{e}diate, il suffit de
remarquer que $rp(K/k)=rp(K'/k)$. Pla\c{c}ons-nous dans les
conditions de $(1)$. En vertu de la proposition
\textcolor{blue}{\ref{p34}}, $k'(H(K/k))/k'$ est $lq$-modulaire et
relativement parfaite, puisque $H(K/k)/k$ l'est ; donc
$k'(H(K/k))$ $\subseteq H(K/k')$. D'autre part, comme $k'/k$ est fini, alors
$H(K/k')/k$ est $lq$-modulaire, et par suite
$rp(H(K/k')/k)=H(K/k)$. D'o\`{u} $k'(rp(H(K/k')/k))=k'(H(K/k))$, ou
encore $k'(H(K/k))=k'rp(H(K/k')/k)=H(K/k')$ (cf.
proposition \textcolor{blue}{\ref{6.2}}). \cqfd
\subsection{Premi\`{e}res applications du th\'{e}or\`{e}me de la cl\^{o}ture $lq$-modulaire}

La premi\`{e}re cons\'{e}quence imm\'{e}diate du th\'{e}or\`{e}me de la cl\^{o}ture
$lq$-modulaire est :
\begin{thm} {\label{th4}} Soit $K/k$ une extension $q$-finie.
Pour toute sous-extension $lq$-modulaire $J/k$ de $K/k$, on a
$di(H(J/k))\leq Ilqm(K/k)-1$, il y'a \'{e}galit\'{e} si $J/k$
contient la cl\^{o}ture $lq$-modulaire de $K/k$.
\end{thm}
\pre
\smartqed
On a $H(J/k)\subseteq H(K/k)$, donc
$di(H(J/k)/k)\leq di(H(K/k)/k)=Ilqm(K/k)$ $-1$. Avec
l'\'{e}galit\'{e} si et seulement si
$di(H(J/k)/k)=di(H(K/k)/k)$, ou encore
$H(J/k)=H(K/k)\subseteq J$. \cqfd \sk

Egalement, comme application du th\'{e}or\`{e}me de la cl\^{o}ture $lq$-modulaire, le produit pr\'{e}serve la $lq$-modularit\'{e}.
Plus pr\'{e}cis\'{e}ment, on a :
\begin{pro} {\label{pr5}} Le produit de deux extensions $lq$-modulaires est $lq$-modulaire.
\end{pro}
\pre
\smartqed
Soient $K_1/k$ et $K_2/k$ deux sous-extensions purement ins\'{e}parables d'une m\^{e}me extension $K/k$.
 Il est imm\'{e}diat que $H(K_1/k)(H(K_2/k))\subseteq H(K/k)$. Si $K_1/k$ et $K_2/k$
sont $lq$-modulaires,  $H(K_1/k)/k$ et $H(K_2/k)/k$ sont
les cl\^{o}tures relativement parfaites respectives de $K_1/k$ et
$K_2/k$, et par suite, $H(K_1/k)(H(K_2/k))$ est la cl\^{o}ture
relativement parfaite de $K_1(K_2)/k$. Il en r\'{e}sulte que
$H(K_1/k)(H(K_2/k$ $))=H(K_1(K_2)/k)=rp(K_1(K_2)/k)$, d'o\`{u} $K_1(K_2)/k$ est
$lq$-modulaire. \cqfd \sk

Plus particuli\`{e}rement, on obtient  le résultat suivant :
\begin{pro} {\label{pr6}} Soient $K_1$ et $K_2$ des corps interm\'{e}diaires d'une extension
$q$-finie $K/k$, $k$-lin\'{e}airement disjoints. Alors, $Ilqm(K_1/k)+Ilqm(K_2/k)-1\leq Ilqm(K_1(K_2)$ $/k)$.
\end{pro}
\pre
\smartqed
Soient $M$, $M_1$, et $M_2$ les cl\^{o}tures
$lq$-modulaires respectives de $K_1(K_2)/k$, $K_1/k$, et $K_2/k$.
On a $M_1(M_2)\subseteq M$ avec $M_1/k$ et $M_2/k$ sont
$k$-lin\'{e}airement disjointes (transitivit\'{e} de la
lin\'{e}arit\'{e} disjointe) ; donc $di(M_1\otimes_k
M_2/k)=di(M_1/k)+di(M_2/k)\leq di(M/k)$ (cf. Corollaire \textcolor{blue}{\ref{cor4}}). Ou encore
$Ilqm(K_1/k)-1+Ilqm(K_2/k)$ $-1\leq Ilqm(K_1(K_2)/k)-1$. Soit
$Ilqm(K_1/k)+Ilqm(K_2/k)-1\leq Ilqm(K_1(K_2)/k)$.
\sk

Soit $K/k$ une extension $q$-finie. Par application r\'{e}p\'{e}t\'{e}e du th\'{e}or\`{e}me  de la cl\^{o}ture $lq$-modulaire,
on d\'{e}finit la $i$-\`{e}me cl\^{o}ture $lq$-modulaire de la fa\c{c}on
suivante : $H_1(K/k)=K$ si $K/k$ est finie, et $H_1(K/k)=H(K/k)$ si
$K/k$ est d'exposant non born\'{e}. Par r\'{e}currence, on pose
$H_i(K/k)=K$ si $K/H_{i-1}(K/k)$ est finie, et
$H_i(K/k)=H(K/H_{i-1}(K/k))$ si $K/H_{i-1}(K/k)$ est d'exposant non
born\'{e}. Par convention, on note  $H_0(K/k)=k$.

On v\'{e}rifie imm\'{e}diatement que :
\sk

\begin{itemize}{\it
\item[$\bullet$] Pour tout $i\in\N$, $H_{i+1}(K/k)/H_i(K/k)$ est $lq$-modulaire.
\item[$\bullet$] Soit $(i,j)\in {\N}^2$. Si $i\leq j$, on a $H_i(H_j(K/k)/k)=H_i(K/k)$ ; et si $i\geq j$,
$H_i(H_j(K$ $/k)/k)=H_j(K/k)$.
}
\end{itemize}
\begin{pro} {\label{ppr7}}
La suite $(H_i(K/k))_{i\in \N}$ est stationnaire.
\end{pro}
\pre
R\'{e}sulte imm\'{e}diatement de la proposition \textcolor{blue}{\ref{pr41}}. \cqfd
\begin{rem} Soit $K/k$ une extension $q$-finie.
Il est clair que $H_1(K/k)=K$ si et seulement si $K/k$ est finie ou
$K/k$ est $lq$-modulaire et relativement parfaite.
\end{rem}

Le r\'{e}sultat ci-dessous peut \^{e}tre consid\'{e}r\'{e} comme une g\'{e}n\'{e}ralisation naturelle de la proposition \textcolor{blue}{\ref{p313}}.
\begin{thm} {\label{t45}} Soient $K/k$ une extension $q$-finie d'exposant non bor\-n\'{e} et $n_0$ le plus petit entier tel que
$K=H_{n_0}(K/k)$. Alors :
\sk

\begin{itemize}{\it
\item[\rm(i)] Pour toute extension finie $K'/K$, pour tout entier  naturel $j$, on a $rp(H_j(K/k))$ $=rp(H_j(K'/k))$.
\item[\rm(ii)] Pour toute  sous-extension finie $k'/k$ de $K/k$, pour tout $i\in\N$, on a $H_i(K/k')$ $=k'(H_i(K/k))$.}
\end{itemize}
\end{thm}
\pre \smartqed
L'item $(i)$ s'en d\'{e}duit par application r\'{e}p\'{e}t\'{e}e de la proposition \textcolor{blue}{\ref{p313}}. Par ailleurs, d'apr\`{e}s la m\^{e}me proposition on a encore $H_1(K/k')=k'(H_1(K/k))$. Supposons que l'on
$k'(H_i(K/k))=H_i(K/k')$. Puisque $k'/k$ est finie, il en
est de m\^{e}me de $k'(H_i(K/k))/H_i(K/k)$. Si $K/H_i(K/k)$ est finie, alors $H_{i+1}(K/k)=H_{i+1}(K/k')=K$. Sinon, en vertu de la
proposition \textcolor{blue}{\ref{p313}}, $k'(H_{i+1}(K/k))$
$=k'(H(K/H_i(K/k)))=k'(H_i(K/k))(H(K/H_i(K/k)))=
H(K/k'(H_i(K/k)))=H(K/H_i(K/k'))=H_{i+1}$ $(K/k')$. Le r\'{e}sultat d\'{e}coule ainsi par r\'{e}currence .\cqfd
\sk

Contrairement aux applications cit\'{e}es en haut, le r\'{e}sultat suivant n'est pas une cons\'{e}quence directe du th\'{e}or\`{e}me de la cl\^{o}ture $lq$-modulaire. Il semble beaucoup moins \'{e}vident.
\begin{pro} {\label{pr8}} Soient $K_1$ et $K_2$ deux corps interm\'{e}diaires d'une
extension $q$-finie. Alors
$K_1(K_2)/K_2$ est $lq$-modulaire si $K_1/k$ l'est.
\end{pro}
\pre
\smartqed
Par induction, d'apr\`{e}s la proposition \textcolor{blue}{\ref{ppr7}} ci-dessus, on se ram\`{e}ne au cas o\`{u} $K_2/k$ est $lq$-modulaire. En vertu de la proposition \textcolor{blue}{\ref{pr5}}, $K_1(K_2)/k$ est $lq$-modulaire, et par la proposition \textcolor{blue}{\ref{le2}}, il en est de m\^{e}me de $K_1(K_2)/K_2$. \cqfd
\sk

L'in\'{e}galit\'{e} suivante r\'{e}sulte imm\'{e}diatement de la proposition ci-dessus.
\begin{pro} {\label{pr9}} Soient $K_1$ et $K_2$ des corps interm\'{e}diaires d'une extension
$q$-finie $K/k$, $k$-lin\'{e}airement disjointes. On a :
$$Ilqm(K_1/k)\leq Ilqm(K_2(K_1)/K_2).$$
\end{pro}
\pre
\smartqed
Imm\'{e}diat \cqfd \sk

Les deux in\'{e}galit\'{e}s qui figurent dans les propositions \textcolor{blue}{\ref{pr6}}  et
\textcolor{blue}{\ref{pr9}}
peuvent \^{e}tre strictes comme le montre l'exemple suivant :
\begin{exe} Soient $Q$ un corps parfait de caract\'{e}ristique
$p>0$, et $k=Q(X,Z_1,Z_2)$ le corps des fractions rationnelles  aux
ind\'{e}termin\'{e}es $X,Z_1,Z_2$. Posons $K_{1,1}=k(X^{p^{-2n}}, \te_1)$ o\`{u} $\te_1={Z_1}^{p^{-1}}X^{p^{-2}}+{Z_2}^{p^{-1}}$, et pour tout  $n\geq 2$, $K_{1,n}=k(X^{p^{-2n}}, \te_n)$ avec $\te_n={Z_1}^{p^{-1}}X^{p^{-2n}}+{(\te_{n-1})}^{p^{-1}}={Z_1}^{p^{-1}}X^{p^{-2n}}+\cdots +{Z_1}^{p^{-n}}{X}^{p^{-n-1}}+ {Z_2}^{p^{-n}}$. Notons \'{e}galement $K_{1}=\displaystyle\bigcup_{n\in\N^*}K_{1,n}$ et
$K_2=\displaystyle\bigcup_{n\in\N}k({Z_1}^{p^{-n}})$, et soit $K=K_1(K_2)$. On v\'{e}rifie aussit\^{o}t que $K = \displaystyle
\bigcup_{n\in\N}k^{p^{-n}}$, et
$K_1/k$ n'est pas $lq$-modulaire, car l'extension $L=k(X^{p^{-\infty}})=lm(K_1/k)$
(cf. exemple \textcolor{blue}{\ref{e1}}), et donc $k\subseteq L\subseteq
H(K_1/k)\subseteq K_1$. Comme $di(K_1/k)=2$, d'apr\`{e}s le th\'{e}or\`{e}me \textcolor{blue}{\ref{thm6}}, $di(H(K_1/k)/k)=1$ ou $2$. Or,  $K_1\not= H(K_1/k)$ (\`{a} savoir $K_1/k$ n'est pas $lq$-modulaire), compte tenu de la propri\'{e}t\'{e} "\^{e}tre relativement parfaite", on aura $di(H(K_1/k)/k)=di(L/k)=1$, et par suite
$L=H(K_1/k)$. On en d\'{e}duit que
$H(K_1/k)(H(K_2/k))=H(K_1/k)(K_2)=L(K_2)\not= H(K/k)=K.$
\end{exe}
\section{Prolongement de la $i$-modularit\'{e}}

Exclusivement, et sauf mention expresse du contraire,
toutes les extensions consid\'{e}r\'{e}es dans le reste de ce
travail sont $q$-finies d'exposant non born\'{e}.
\subsection{Suite de d\'{e}composition d'une extension $q$-finie}

Par analogie aux groupes \`{a} op\'{e}rateur, et plus particuli\`{e}rement aux suites de Jordan-Hölder, nous adoptons la d\'{e}finition suivante :
\begin{df} Toute suite d'extensions $q$-finies $k=F_0\subseteq F_1\subseteq
\cdots\subseteq F_n$ telle que $F_{i+1}/F_i$ est $lq$-modulaire pour
tout entier $i\in [0,n[$ s'appelle suite de d\'{e}composition de rang $n$.
\end{df}

En outre, soit $K/k$ une extension $q$-finie, toute suite de d\'{e}composition $k=F_0\subseteq F_1\subseteq
\cdots\subseteq F_n$ telle que $F_n=K$ sera appel\'{e}e  suite de d\'{e}composition de rang $n$ associ\'{e}e \`{a} $K$.
\sk

Comme r\'{e}sultat \'{e}l\'{e}mentaire, on a :
\begin{pro} Toute extension $q$-finie admet une suite de d\'{e}composition. En particulier, toute suite de
d\'{e}composition associ\'{e}e \`{a} une sous-extension $L/k$ de $K/k$ se prolonge en une suite de
d\'{e}composition associ\'{e}e \`{a} $K/k$.
\end{pro}
\pre R\'{e}sulte imm\'{e}diatement des propositions \textcolor{blue}{\ref{ppr7}} et \textcolor{blue}{\ref{pr8}}.\cqfd
\sk

Soient $(F_i/k)_{1\leq i\leq n}$ et $(F'_i/k')_{1\leq i\leq n}$ deux
suites d'extensions purement ins\'{e}parables respectivement de $k$
et $k'$. On dit que $(F_i/k)_{1\leq i\leq n}$ est \'{e}quivalente
\`{a} $(F'_i/$ $k')_{1\leq i\leq n}$, et on note $(F_i/k)_{1\leq
i\leq n}\sim (F'_i/k')_{1\leq i\leq n}$, si $k\sim k'$, et si pour
tout $i\in \{1\cdots, n\}$, $F_i\sim F'_i$. On rappelle que la suite
$(F_i/k)_{1\leq i\leq n}$ est dite croissante si pour tout $i\in
\{1,\cdots,n-1\}$, on a $F_i\subseteq F_{i+1}$.  On \'{e}tablit
imm\'{e}diatement le r\'{e}sultat suivant :
\begin{pro} {\label{p41}}
Soient  $(F_i/k)_{1\leq i\leq n}$ et $(F'_i/k')_{1\leq i\leq n}$
deux suites croissantes d'extensions $q$-finies
respectivement de $k$ et $k'$. Si $(F_i/k)_{1\leq i\leq n}\sim
(F'_i/k')_{1\leq i\leq n}$, alors $k\subseteq F_1\subseteq \cdots
\subseteq F_n$ est une suite de d\'{e}composition de rang $n$ si et
seulement si il en est de m\^{e}me de $k'\subseteq F'_1\subseteq
\cdots \subseteq F'_n$.
\end{pro}
\pre \smartqed Cela r\'{e}sulte aussit\^{o}t de la proposition
\textcolor{blue}{\ref{p34}}. \cqfd
\sk

Comme cons\'{e}quence imm\'{e}diate, le r\'{e}sultat suivant permet de ramener l'\'{e}tude portant sur les suites de d\'{e}composition \`{a} celle des suites de d\'{e}composition dont les termes sont relativement parfaits.
\begin{cor} {\label{corr1}} $k\subseteq F_1\subseteq \cdots\subseteq F_n$
est une suite de d\'{e}composition si et seulement si il en est de m\^{e}me de
$k\subseteq rp(F_1/k)\subseteq\cdots\subseteq rp(F_n/k)$ ;
\end{cor}

Le r\'{e}sultat suivant est une g\'{e}n\'{e}ralisation naturelle de la proposition \textcolor{blue}{\ref{pr8}}.
\begin{pro} {\label{p42}} Soit $k\subseteq F_1\subseteq\cdots\subseteq F_n$ une suite de
d\'{e}composition. Pour toute extension $q$-finie
$L/k$, $L\subseteq L(F_1)\subseteq\cdots\subseteq L(F_n)$ est aussi
une suite de d\'{e}composition.
\end{pro}
\pre \smartqed Imm\'{e}diat. \cqfd
\sk

Egalement, voici une extension non triviale de la proposition
\textcolor{blue}{\ref{pr5}}.
\begin{pro} {\label{p43}} Le produit de deux suites de d\'{e}composition est une suite de
d\'{e}composition, c'est-\`{a}-dire si $k\subseteq F_1\subseteq
\cdots\subseteq F_n$ et $k\subseteq F'_1\subseteq\cdots\subseteq
F'_n$ sont deux suites de d\'{e}compositions, il en est de m\^{e}me de
$k\subseteq F_1(F'_1)\subseteq \cdots\subseteq F_n(F'_n)$.
\end{pro}
\pre \smartqed En vertu de la proposition
\textcolor{blue}{\ref{pr8}}, on a aussit\^{o}t
$F_{i+1}(F'_i)/F_i(F'_i)$ et $F_i(F'_{i+1})$ $/F_i(F'_i)$ sont
$lq$-modulaires ; et d'apr\`{e}s la proposition
\textcolor{blue}{\ref{pr5}},  $F_{i+1}(F'_{i+1})/F_i(F'_i)$ est
aussi $lq$-modulaire. \cqfd
\begin{pro} {\label{p44}} Soit $K/k$ une extension $q$-finie.
Pour toute suite de d\'{e}composition $k\subseteq F_1\subseteq
\cdots\subseteq F_n=K$,  pour tout entier $i\in [1,n]$, on a $rp(F_i/k)\subseteq H_i(K/k)$.
\end{pro}
\pre \smartqed Pour tout $i\in \N$,
on pose $H_i(K/k)=H_i$. On a $F_1/k$ est $lq$-modulaire, donc
$rp(F_1/k)\subseteq H_1$ par d\'{e}finition de la cl\^{o}ture $lq$-modulaire. De proche en proche si $rp(F_i/k)\subseteq
H_i$, comme $F_{i+1}/F_i$ est $lq$-modulaire,  en vertu du corollaire \textcolor{blue}{\ref{c35}}, il en est de m\^{e}me de
$F_{i+1}/rp(F_i/k)$.  D'apr\`{e}s la proposition
\textcolor{blue}{\ref{pr8}}, $H_i(F_{i+1})/H_i$ est aussi
$lq$-modulaire. D'o\`{u} $rp(F_{i+1}/k)\subseteq
rp(H_i(F_{i+1})/k)\subseteq  H_{i+1}$. \cqfd
\sk

 En vertu des propositions
\textcolor{blue}{\ref{p44}},  \textcolor{blue}{\ref{p42}},
\textcolor{blue}{\ref{p43}}, et le corollaire \textcolor{blue}{\ref{corr1}}, on a :
\begin{cor} Soient $K_1$ et $K_2$ deux corps interm\'{e}diaires d'une extension $q$-finie $K/k$.  Les assertions suivantes sont v\'{e}rifi\'{e}es :
\sk

\begin{itemize}{\it
\item[{\rm (i)}] Pour toute sous-extension $L/k$ de $K/k$, pour tout $i\in \N$, on a $rp(H_i(L/k)/k)$ $\subseteq
H_i(K/k)$.
\item[{\rm (ii)}] $rp(H_i(K_1/k)(H_i(K_2/k))/k)
\subseteq H_i(K_1(K_2)/k)$ et  $K_1(H_i(K_2/k))\subseteq
H_i(K_1(K_2)$ $/K_1)$ pour tout $i\in \N$.}
\end{itemize}
\end{cor}

Le th\'{e}or\`{e}me suivant permet de caract\'{e}riser le rang minimal d'une suite de d\'{e}composition.
\begin{thm} {\label{t45}} Soient $K/k$ une extension $q$-finie d'exposant non bor\-n\'{e} et $n_0$ le plus petit entier tel que
$K=H_{n_0}(K/k)$. Alors :
\sk

\begin{itemize}{\it
\item[{\rm (i)}] Si $K/k$ est relativement parfaite,
$k\subseteq H_1(K/k)\subseteq\cdots\subseteq H_{n_0}(K/k)=K$ est une
suite de d\'{e}composition de rang minimal associ\'{e}e \`{a} $K/k$.
\item[{\rm (ii)}] Si $K/k$ n'est pas relativement parfaite, alors $2\leq n_0$ et $k\subseteq H_1(K/k)\subseteq
\cdots\subseteq H_{n_0-2}(K/k)\subseteq H_{n_0}(K/k)=K$ est une suite de
d\'{e}composition de rang minimal associ\'{e}e \`{a} $K/k$.}
\end{itemize}
\end{thm}
\pre \smartqed Compte tenu de la proposition
\textcolor{blue}{\ref{p44}}, les assertions propos\'{e}es sont valides. 
\sk

On en d\'{e}duit imm\'{e}diatement le r\'{e}sultat suivant :
\begin{cor} {\label{corrr1}}
Pour toute  extension $q$-finie  $K/k$,  $lqm(K/k)$  est non triviale
($K\not= lqm(K/k)$). Plus pr\'{e}cis\'{e}ment, si $K/k$ est
d'exposant non born\'{e}, il en est de m\^{e}me de $K/lqm(K/k)$. En
outre, $di(lqm(K/k)/k)< di(K/k)$.
\end{cor}
\pre \smartqed Si $K/k$ est finie, on a $lqm(K/k)=k$, donc ce cas
est trivial. Si $K/k$ est d'exposant non born\'{e}, d'apr\`{e}s la
proposition \textcolor{blue}{\ref{p34}}, on se ram\`{e}ne au cas
o\`u $K/k$ est relativement parfaite. Soit $n_0$ le plus petit
entier tel que $H_{n_0}(K/k)=K$. D'apr\`{e}s le th\'{e}or\`{e}me
\textcolor{blue}{\ref{t45}} ci-dessus, $k=H_0(K/k)\subseteq
H_1(K/k)\subseteq\cdots\subseteq H_{n_0}(K/k)=K$ est une suite de
d\'{e}composition de rang minimal associ\'{e}e \`{a} $K/k$. Il en
r\'{e}sulte que $lqm(K/k)\subseteq H_{n_0-1}(K/k)$, puisque
$K/H_{n_0-1}(K/k)$ est $lq$-modulaire. De plus, $K/lqm(K/k)$ est
d'exposant non born\'{e}, car $K/H_{n_0-1}(K/k)$ l'est par
construction. 
\sk

Plus particuli\`{e}rement, cela conduit  \`{a} :
\begin{cor} La plus petite sous-extension $m/k$ d'une extension $q$-finie $K/k$ telle que
$K/m$ est modulaire n'est pas triviale ($m\not= K$). Plus
pr\'{e}cis\'{e}ment, si $K/k$ est d'exposant non born\'e, il en est de m\^eme de $K/m$.
\end{cor}
\pre Puisque $rp(lm(K/k)/k)=lqm(K/k)$ (cf. proposition
\textcolor{blue}{\ref{p318}}) et $K/K(K^p)$ est modulaire, alors
$lm(K/k)$ est aussi non triviale. Si de plus $K/k$ est d'exposant
non born\'{e}, il en est de m\^{e}me de $K/lm(K/k)$. \cqfd
\subsection{Degr\'{e} de modularit\'{e} d'une extension}
\begin{df} Soit $K/k$ une extension $q$-finie.
$K/k$ est dite $i$-modulaire si elle admet une suite de
d\'{e}composition de rang $i$ ; autrement dit, s'il existe une suite d'extensions
$(F_j/k)_{0\leq j\leq i}$ telle que $k=F_0\subseteq
F_1\subseteq \cdots\subseteq F_i=K$ et $F_{j+1}/F_j$ est
$lq$-modulaire pour tout entier $j\in [0,i[$. Le plus petit entier $i$ tel
que $K/k$ est $i$-modulaire s'appelle le degr\'{e} de modularit\'{e}
de $K/k$, et se note $dm(K/k)$.
\end{df}
\begin{rem} Le degr\'{e} de modularit\'{e}
permet de mesurer le niveau  de modularit\'{e} d'une extension.
\end{rem}

Il est imm\'{e}diat que :
\sk

\begin{itemize}{\it
\item[$\bullet$] Toute extension $lq$-modulaire est $1$-modulaire. En particulier,  $dm(K/k)=1$ si et seulement si $K/k$ est $lq$-modulaire.
\item[$\bullet$] $K/k$ est $i$-modulaire pour tout entier $i\geq dm(K/k)$.}
\end{itemize}
\sk

Le r\'{e}sultat qui suit caract\'{e}rise les extensions $i$-modulaires. D'une mani\`{e}re pr\'{e}cise, on a :
\begin{thm} {\label{t47}}Soit $K/k$ une extension $q$-finie d'exposant non born\'{e}. $K/k$ est $i$-modulaire
si et seulement si $H_i(K/k)$ co\"{\i}ncide approximativement avec
$K$, c'est-\`{a}-dire $H_i(K/k)\sim K$. En outre, $dm(K/k)$
$=\inf\{i\in \N$ tel que $H_i(K/k)\sim K\}$.
\end{thm}
\pre \smartqed Cela r\'{e}sulte  du th\'{e}or\`{e}me \textcolor{blue}{\ref{t45}} et de la
proposition \textcolor{blue}{\ref{p44}}. \cqfd
\sk

On v\'{e}rifie ais\'{e}ment que :
\sk

\begin{itemize}{\it
\item $dm(K/k)=n$ si
est seulement si $H_n(K/k)=rp(K/k)$.
\item Pour tout entier $j$, $H_j(K/k)/k$ est $j$-modulaire.
\item Pour tout entier  $j\in [1,dm(K/k)[$, $dm(K/H_j(K/k))=dm(K/k)-j$ et
$dm(H_j(K/k)/k)=j$.}
\end{itemize}
\sk

Le r\'{e}sultat qui suit g\'{e}n\'{e}ralise la proposition \textcolor{blue}{\ref{p34}}.
En d'autres termes, comme dans le cas de la $lq$-modularit\'{e}, la $i$-modularit\'{e}
est pr\'{e}serv\'{e}e \`{a} une extension finie pr\`{e}s. Plus pr\'{e}cis\'{e}ment, on a :
\begin{pro} {\label{p46}} Soit $K/k$ une extension $q$-finie d'exposant non born\'{e}. Alors :
\sk

\begin{itemize}{\it
\item[{\rm (1)}] Si $k'\sim k$ et
$k'\subset K$, $K/k$ est $i$-modulaire si et seulement si il en est de
m\^{e}me de $K/k'$. En outre, $dm(K/k)=dm(K/k')$.
\item[{\rm (2)}] Si $K\sim K'$ et $k\subset
K'$, $K/k$ est $i$-modulaire si et seulement si $K'/k$ l'est aussi. En outre, $dm(K/k)=dm(K'/k)$.
\item[{\rm (3)}] Si $k'\sim k$ et $K\sim K'$, avec $k'\subset
K'$, alors $K/k$ est $i$-modulaire si et seulement si il en est de
m\^{e}me de $K'/k'$. En outre, $dm(K/k)=dm(K'/k')$.}
\end{itemize}
\end{pro}
\pre \smartqed $(1)$ et $(2)$ d\'{e}coulent imm\'{e}diatement de la proposition
\textcolor{blue}{\ref{p41}} et du th\'{e}or\`{e}me
\textcolor{blue}{\ref{t47}}, et l'assertion $(3)$ r\'{e}sulte aussit\^{o}t de $(1)$ et $(2)$. \cqfd
\sk

Comme cons\'{e}quence imm\'{e}diate, le r\'{e}sultat suivant permet de ramener l'\'{e}tude de la $i$-modularit\'{e} au cas des
extensions relativement parfaites.
\begin{cor} {\label{c49}} Soit $K/k$ une extension $q$-finie d'exposant non born\'{e}. Alors :
\begin{itemize}{\it
\item[{\rm (i)}] $K/k$ est $i$-modulaire si et seulement si il en est de m\^{e}me de $rp(K/k)/k$.
En outre, $dm(K/k)=dm(rp(K/k)/k)$.
\item[{\rm (ii)}] Pour toute extension finie $L/k$, on a $K/k$ est $i$-modulaire si et seulement si
il en de m\^{e}me pour $L(K)/L$. En outre, $dm(K/k)=dm(L(K)/L)$.}
\end{itemize}
\end{cor}
\pre \smartqed Il suffit en effet de remarquer que $rp(K/k)\sim K$,
$L\sim k$ et $L(K)\sim K$. 
\sk

Le r\'{e}sultat ci-dessous montre que le niveau de modularit\'{e} d'une extension ne d\'{e}passe jamais l'entier $di(rp(K/k)/k)$.
\begin{thm} {\label{t48}} Soit $K/k$ une extension $q$-finie d'exposant non born\'{e}, alors  $dm(K/k)$ $\leq di(rp(K/k)/k)$.
\end{thm}
\pre \smartqed En vertu du corollaire
\textcolor{blue}{\ref{c49}} ci-dessus, on se ram\`{e}ne au cas o\`u $K/k$ est
relativement parfaite. Soit $n=dm(K/k)$, donc $H_0(K/k)=k\subseteq
\cdots\subseteq H_n(K/k)=K$ est une suite de d\'{e}composition de rang
minimal associ\'{e}e \`{a} $K/k$, avec $H_j(K/k)/k$ est relativement
parfaite pour tout entier $j\in [1, n]$. D'o\`u d'apr\`{e}s la proposition \textcolor{blue}{\ref{pr12}}, $n\leq
\di\sum_{j=0}^{n-1}di(H_{j+1}(K/k)/H_j(K/k))=di(K/k)$, (car pour tout $j\in\{0,\cdots,n-1\}$, $H_{j+1}(K/k)\not = H_j(K/k)$, et donc
$1\leq di(H_{j+1}(K/k)/H_j(K/k))$).
Il en r\'{e}sulte que $n=dm(K/k)\leq di(rp(K/k)/k)$. \cqfd
\begin{pro} {\label{p410}}
Soient $k\subset L \subset K$ des extensions $q$-finies. Si $L/k$ est $i$-modulaire et
$K/L$ est $j$-modulaire, alors $K/k$ est $(i+j)$-modulaire. En outre,
$dm(K/k)\leq dm(L/k)+dm(K/L)$.
\end{pro}
\pre \smartqed En effet, si $k\subseteq F_1\subseteq\cdots\subseteq
F_i=L$ et $L\subseteq F'_1\subseteq\cdots\subseteq F'_j=K$ sont deux
suites de d\'{e}compositions, il en est de m\^{e}me de $k\subseteq
F_1\subseteq\cdots\subseteq F_i=L\subseteq
F'_1\subseteq\cdots\subseteq F'_j=K$. \cqfd
\sk

La $i$-modularit\'{e} est stable par le produit. Autrement dit, on a :
\begin{pro}
Soient $K_1$ et $K_2$ deux corps interm\'{e}diaires d'une extension
$q$-finie $K/k$. Si $K_1/k$ et
$K_2/k$ sont $i$-modulaires, il en est de m\^{e}me de $K_1(K_2)/k$.
En particulier, $dm(K_1(K_2)/k)\leq \sup(dm(K_1/k), dm(K_2/k))$.
\end{pro}
\pre \smartqed Imm\'{e}diat, car le produit de deux suites de
d\'{e}composition est une suite de d\'{e}composition. \cqfd
\sk

La $i$-modularit\'{e} est stable par changement du corps de base dans
le sens ascendant. Plus g\'{e}n\'{e}ralement, on a :
\begin{pro}
Soit $K/k$ une extension $i$-modulaire. Pour toute extension $L/k$,
on a $L(K)/L$ est $i$-modulaire. En outre, $dm(L(K)/L)\leq dm(K/k)$.
\end{pro}
\pre \smartqed Cela r\'{e}sulte imm\'{e}diatement de la proposition
\textcolor{blue}{\ref{p42}}. \cqfd
\sk

Comme cas particulier, on a :
\begin{cor} Soient $k\subseteq L\subseteq K$ des extensions $q$-finies.
Si $K/k$ est $i$-modulaire, il en est de m\^{e}me de $K/L$. En outre,
$dm(K/L)\leq dm(K/k)$.
\end{cor}
\pre \smartqed Imm\'{e}diat. \cqfd
\sk

Dans ce qui suit nous allons d\'{e}crire les extensions $q$-finies $K/k$ dont
le degr\'{e} de modularit\'{e} atteint la taille de $K/k$.
\subsection{Extensions $e$-ferm\'{e}es}
\begin{df} Toute extension $q$-finie $K/k$ v\'{e}rifiant $dm(K/k)=di(K/k)$
s'appelle extension $e$-ferm\'{e}e d'ordre $di(K/k)$.
\end{df}

Il est \'{e}vident que :
\sk

\begin{itemize}{\it
\item[$\bullet$] $K/K$ est $e$-ferm\'{e}e d'ordre 0.
\item[$\bullet$] $K/k$ est $e$-ferm\'{e}e d'ordre $1$ si et seulement si $K/k$ est $q$-simple.}
\end{itemize}
\begin{exe} {\label{e2}} Reprenons les notations de l'exemple \textcolor{blue}{\ref{e1}}. Soient $k_0$
un corps parfait de caract\'{e}ristique $p>0$ et $(X,Z_1,Z_2)$ une
famille alg\'{e}briquement ind\'{e}pendante sur $k_{0}$.
Posons $k=k_0(X,Z_1,Z_2)$ et $K_{n}=k(X^{p^{-2n}},\theta _{n})$, avec $\theta _{1}=Z_1^{p^{-1}}X^{p^{-2}}+Z_2^{p^{-1}}$, et pour tout entier $n\geq 2$,
\[\theta _{n}=Z_1^{p^{-1}}X^{p^{-2n}}+{(\theta _{n-1})}^{p^{-1}}=Z_1^{p^{-1}}X^{p^{-2n}}+Z_{1}^{p^{-2}}X^{p^{-2n+1}}+
\cdots+Z_1^{p^{-n}}X^{p^{-n-1}}+Z_2^{p^{-n}}.\]
Soit $K=\displaystyle\bigcup_n K_{n}$.
\end{exe}
\begin{pro} L'extension $K/k$ ci-dessus est $e$-ferm\'{e}e d'ordre 2.
\end{pro}
\pre \smartqed Imm\'{e}diat, car
$lqm(K/k)=k({X}^{p^{-\infty}})=H_1(K/k)$. \cqfd
\begin{pro} Toute extension $e$-ferm\'{e}e d'exposant non born\'{e} est relativement parfaite.
\end{pro}
\pre \smartqed On a $K/k$ est $e$-ferm\'{e}e, donc $di(K/k)=dm(K/k)$. Mais, comme $di(K/k)=di(rp(K/k)/k)+di(K/k(K^p))$ et $dm(K/k)= dm(rp(K/k)/k)\leq di(rp(K/k)/k)$, il en r\'{e}sulte que
$di(K/k(K^p))=0$, ou encore $K/k$ est relativement parfaite. \cqfd
\sk

La propri\'{e}t\'{e} "\^{e}tre $e$-ferm\'{e}e" est non seulement pr\'{e}serv\'{e}e \`{a} une extension
finie pr\`{e}s, mais elle est aussi semi-transitive. En d'autres termes, on a :
\begin{pro} {\label{p415}}
Pour toute  extension $q$-finie $K/k$, les assertions suivantes sont vraies :
\sk

\begin{itemize}{\it
\item[{\rm (i)}] Si $k\sim k_1$ et $k_1\subset K$,  $K/k$
est $e$-ferm\'{e}e d'ordre $n$ si et seulement si il en est de
m\^{e}me de $K/k_1$.
\item[{\rm (ii)}] Pour toute sous-extension relativement parfaite
$L/k$ de $K/k$, si $K/k$ est $e$-ferm\'{e}e, il en
est de m\^{e}me de $K/L$ et $L/k$.}
\end{itemize}
\end{pro}
\pre \smartqed L'assertion $(i)$ est trivialement \'{e}vidente (cf. proposition
\textcolor{blue}{\ref{p46}}).  Posons ensuite $n=dm(K/k)$. En vertu de la proposition \textcolor{blue}{\ref{p410}}, et comme $dm(L/k)\leq di(L/k)$ et $dm(K/L)\leq di(K/L)$, on aura $n\leq dm(K/L)+dm(L/k)\leq di(K/L)+di(L/k)=di(K/k)=n$,  et par suite $dm(K/L)=di(K/L)$ et $dm(L/k)=di(L/k)$. \cqfd
\sk

On en d\'{e}duit aussit\^{o}t le r\'{e}sultat suivant :
\begin{cor} Soit $K/k$ une extension  $q$-finie. Si $K/k$ est $e$-ferm\'{e}e d'ordre $n$, pour tout entier $j\in [1,n]$, $H_j(K/k)/k$ et $K/H_j(K/k)$ sont $e$-ferm\'{e}es d'ordre respectivement $j$ et $n-j$.
\end{cor}
\pre Imm\'{e}diat. \cqfd
\sk

D'apr\`{e}s le corollaire pr\'{e}c\'{e}dent, si $K/k$ est $e$-ferm\'{e}e, il
en est de m\^{e}me de $H_j(K/k)$ $/H_i(K/k)$ pour tout entiers naturels $i\leq j$ ; et
comme cons\'{e}quence on a :
\begin{thm}
Toute extension $e$-ferm\'{e}e $K/k$ d'ordre $n$ se d\'{e}compose
compl\`{e}tement en $n$ extensions $q$-simples. D'une mani\`{e}re
pr\'{e}cise, on a $H_0(K/k)=k\subseteq
H_1(K/k)\subseteq\cdots\subseteq H_{n-1}(K/k)\subseteq K=H_n(K/k)$,
avec $H_{i+1}(K/k)/H_i(K/k)$ est $q$-simple pour tout $i\in
\{0,\cdots,n-1\}$.
\end{thm}
\pre \smartqed Imm\'{e}diat, il suffit de remarquer que
$H_{j+1}(K/k)/H_j(K/k)$ est $e$-ferm\'{e}e d'ordre $1$ pour tout entier $j\in [0,n[$,
ou encore $H_{j+1}(K/k)/H_j(K/k)$ est $q$-simple. 
\sk

Voici une autre caract\'{e}risation des extensions $e$-ferm\'{e}es.  Elle permet  d'identifier  la plus petite sous-extension $m$ telle que $K/m$ est $lq$-modulaire
\begin{pro}
Soit $K/k$ une extension $e$-ferm\'{e}e d'ordre $n$, alors
$H_{n-1}(K/k)$ est la plus petite sous-extension de $K/k$ telle que
$K/H_{n-1}(K/k)$ est $lq$-modulaire.
\end{pro}
\pre \smartqed Soit $m=lqm(K/k)$, donc $m\subseteq H_ {n-1}(K/k)$,
car $K/H_{n-1}(K/k)$ est $lq$-modulai\-re. D'autre part, on a $K/m$
est $lq$-modulaire et $m/k$ est relativement parfaite, donc $K/m$
est $e$-ferm\'{e}e d'ordre $1$ en vertu de la proposition
\textcolor{blue}{\ref{p415}}. D'o\`u, par la proposition \textcolor{blue}{\ref{pr12}}, $di(K/H_{n-1}(K/k))=di(K/m)=di(K/H_{n-1}(K/k))+di(m/H_{n-1}(K/k))$,
ou encore $di(m/H_{n-1}(K/k))=0$, i.e. $m=H_{n-1}(K/k)$. \cqfd
\sk

Nous obtenons  imm\'{e}diatement le r\'{e}sultat suivant :
\begin{cor} Soit $K/k$ une extension $e$-ferm\'{e}e d'ordre $n$.
Pour tous entiers naturels $i< j\leq n$,  on a
$H_{j-1}(K/k)=lqm(H_{j}(K/k)/H_{i}(K/k))$.
\end{cor}
\pre \smartqed Imm\'{e}diat. \cqfd
\begin{pro} {\label{p417}} Soient $K/k$ une extension $q$-finie et
relativement parfaite,  et $m/k$ la plus petite
sous-extension de $K/k$ telle que $K/m$ est $lq$-modulaire. Si
$K/m$ est $q$-simple, alors $m/k$ est absorbante, c'est-\`{a}-dire
pour toute sous-extension propre $L/k$ de $K/k$,
on a $rp(L/k)\subseteq m$.
\end{pro}
\pre \smartqed Posons $di(K/k)=n$.
Comme $K/m$ est $q$-simple et $m/k$ est relativement parfaite, on a
$n=di(K/k)=di(K/m)+di(m/k)=1+di(m/k)$. Ou encore $di(m/k)=n-1$.
Consid\'{e}rons maintenant une sous-extension propre $L/k$  de $K/k$, et soit $N=lqm(K/L)$ et $M=rp(N/k)$, donc $m\subseteq M$. En
particulier, $di(m/k)=n-1\leq di(M/k)\leq di(K/k)=n$ (cf. th\'{e}or\`{e}me
\textcolor{blue}{\ref{thm6}}). Si $di(M/k)=di(K/k)$, compte tenu de la propri\'{e}t\'{e} "\^{e}tre relativement parfaite" respectivement de $K/k$ et $M/k$, on obtient
$K=M$, contradiction. D'o\`{u} $di(m/k)=di(M/k)$. Egalement, comme $m/k$ et
$M/k$ sont relativement parfaites, on aura $m=M$. Plus particuli\`{e}rement, $rp(L/k)\subseteq m$. \cqfd
\sk

Une application type de la proposition pr\'{e}c\'{e}dente est le
r\'{e}sultat suivant :
\begin{thm}
Soit $K/k$ une extension $e$-ferm\'{e}e d'ordre $n$. Alors $(H_i(K/k)/k)$,
$(1\leq i\leq n)$, sont les seules sous-extensions
relativement parfaites de $K/k$.
\end{thm}
\pre \smartqed  Soient $L/k$ une sous-extension relativement
parfaite de $K/k$ et $j$ le plus petit entier tel que $L\subseteq
H_j(K/k)$. Si $L\not= H_j(K/k)$, comme $H_j(K/k)/k$ est $e$-ferm\'{e}e
d'ordre $j$ ; d'apr\`{e}s la proposition \textcolor{blue}{\ref{p417}}
ci-dessus, $rp(L/k)=L\subseteq H_{j-1}(K/k)$, c'est une contradiction. \cqfd\sk
\begin{cor} Soit $K/k$ une extension $e$-ferm\'{e}e d'ordre $n$. Pour tout $L\in [k,K]$, il existe un entier $i\in[1,n]$ tel que $L\sim H_i(K/k)$.
\end{cor}

Le r\'{e}sultat ci-dessus peut se traduire  par les $H_i(K/k)$ sont les seuls corps interm\'{e}diaires \`{a} une une extension finie pr\`{e}s d'une extension $e$-ferm\'{e}e $K/k$.
\subsection{Existence d'une extension $e$-ferm\'{e}e d'ordre $n$}

Comme g\'{e}n\'{e}ralisation de l'exemple \textcolor{blue}{\ref{e2}}, nous allons construire une extension
$e$-ferm\'{e}e d'ordre $n>2$. Il est \`{a} signaler que cet exemple  reprend  (\cite{Che-Fli5},   Exemple \textcolor{blue}{4}) en le corrigeant.
\begin{exe} {\label{e3}}  Soient $k_0$  un corps parfait  de caract\'{e}ristique
$p>0$ et $k=k_0(X,\al_1,\cdots,$ $\al_n,\be_1\cdots,\be_n)$ le corps des fractions rationnelles aux ind\'{e}termin\'{e}es $(X,\al_1,\cdots,\al_n,$ $\be_1\cdots,\be_n)$. On pose :\\
$T_0(X)=(T_{1}^0,\cdots, T_{n}^0,\cdots)$ o\`{u}
$T_{n}^0=X^{p^{-n}}$.\\
$T_1(X)=(T_{1}^1,\cdots,T_{n}^1,\cdots)$ o\`{u}
$T_{1}^1={\al_1}^{p^{-1}}T_{2}^0+{\be_1}^{p^{-1}}$ et
$T_{n}^1={\al_1}^{p^{-1}} T_{2n}^{0}+ {(T_{n-1}^{1})}^{p^{-1}}$.
\\
Par r\'{e}currence, on note
 $T_j(X)=(T_{1}^j,\cdots, T_{n}^j,\cdots)$ o\`{u} $T_{1}^j={\al_j}^{p^{-1}} T_{2}^{j-1}+{\be_j}^{p^{-1}}$ et
$T_{n}^j={\al_j}^{p^{-1}}T_{2n}^{j-1}+ {(T_{n-1}^ {j})}^{p^{-1}}$.
\end{exe}

Pour tout entier $j\in [1,n]$, on note aussi $K_i=k(T_0(X),\cdots,T_i(X))$ , et $K_0=k$ par convention. Par construction,
on a :
\begin{eqnarray*}
T_{n}^j&=&{\al_j}^{p^{-1}}T_{2n}^{j-1}+ {(T_{n-1}^{j})}^{p^{-1}}, \\
&=& {\al_{j}}^{p^{-1}}T_{2n}^{j-1}+\dots+{\al_{j}}^{p^{-i}}{(T_{2(n-(i-1))}^{j-1})}^{p^{-i+1}}+ {(T_{n-i}^{j})}^{p^{-i}},\\
&=& {\al_{j}}^{p^{-1}}T_{2n}^{j-1}+\dots+{\al_{j}}^{p^{-n}}{(T_{2}^{j-1})}^{p^{-n+1}}+{\be_{j}}^{p^{-n}} ;
\end{eqnarray*}
donc pour tout entiers naturels $i<j$, on obtient
${(T_{n}^j)}^{p^i}=\theta_{n,j}+T_{n-i}^ {j}$, avec $\theta_{n,j}\in
K_{j-1}$. Il en r\'{e}sulte que
$K_{j-1}({(T_{n}^{j})}^{p^i})=K_{j-1}(T_{n-i}^j)$ pour tout entier $j\in [1,n[$. En particulier, $K_j/K_{j-1}$ est $q$-simple.
\begin{thm} {\label{t420}}Sous les hypoth\`{e}ses ci-dessus, $K_{r-1}/k$ est la
plus petite sous-extension de $K_r/k$ telle que $K_r/K_{r-1}$ est
modulaire pour tout entier $r\in [1,n]$.
\end{thm}

La d\'{e}monstration de ce th\'{e}or\`{e}me utilise les r\'{e}sultats auxiliaires suivants.
D'abord, pour tout entier $j\in [1,n]$, on pose $L_j=k_0(X^{p^{-\infty}}, ({\al_{i}}^{p^{-\infty}}, \be_{i}^{p^{-\infty}})_{i\in \{1,\cdots,j\}})$.
\begin{lem} {\label{lem110}} Pour tout entier $i\in [1,n]$, pour tout $s\in \N^*$, on a
$T_{s}^{i}\in L_i$. En outre, $K_i\subseteq L_i$.
\end{lem}
\pre
On va utiliser une d\'{e}monstration par r\'{e}currence.  Par construction, pour tout entier $s\in[1,n]$, $T_{s}^1={\al_{1}}^{p^{-1}}T_{2s}^{0}+\dots+ {\al_{1}}^{p^{-s}}{(T_{2}^{0})}^{p^{-s+1}}+{\be_{1}}^{p^{-s}}$, d'o\`u $T_{s}^{1}\in k_0(T_{2s}^{0},\dots, {(T_{2}^0)}^{p^{-s+1}}, {\al_{1}}^{p^{-s}}, {\be_{1}}^{p^{-s}})\subseteq k_0({X}^{p^{-\infty}}, {\al_{1}}^{p^{-\infty}},{\be_{1}}^{p^{-\infty}})=L_1$, et par suite le lemme est v\'erifi\'e pour le premier rang. Supposons que la propri\'{e}t\'{e} de r\'{e}currence s'applique jusqu'\`{a} l'ordre $i>1$. Egalement, pour tout $s\in \N^*$,  on a $T_{s}^{i+1}={\al_{i+1}}^{p^{-1}}T_{2s}^{i}+\dots+{\al_{i+1}}^{p^{-s}}{(T_{2}^i)}^{p^{-s+1}}+{\be_{i+1}}^{p^{-s}}$, et donc $T_{s}^{i+1}\in k_0(T_{2s}^{i},\dots,{(T_{2}^{i})}^{p^{-s+1}},$ $ {\al_{i+1}}^{p^{-s}}, {\be_{i+1}}^{p^{-s}})$. Or, d'apr\`{e}s l'hypoth\`{e}se de r\'{e}currence, pour tous $r\in \N^*$,
$T_{r}^{i}\in L_i$. Comme $k_0$ est parfait, i.e. ${k_0}^{p}=k_0$, alors $L_i$ est aussi parfait par construction, et donc pour tout entier naturel non nul $n$, on aura ${(T_{r}^{i})}^{p^{-n}}\in {L_i}^{p^{-n}}=L_i$. Par cons\'{e}quent, $T_{s}^{i+1}\in L_i({\al_{i+1}}^{p^{-s}}, {\be_{i+1}}^{p^{-s}})\subseteq L_{i+1}$. \cqfd
\begin{lem} Pour tout $s\in \{1,\cdots,n\}$, on a ${\al_s}^{p^{-1}}\not \in K_n$.
\end{lem}
Pour simplifier l'\'{e}criture, pour tout entier $j\in [1,n]$, on pose $N_j=k(L_j)=k(X^{p^{-\infty}}, ({\al_{i}}^{p^{-\infty}}, \be_{i}^{p^{-\infty}}$ $)_{i\in \{1,\cdots,j\}})$.
D'apr\`es le lemme ci-dessus, pour tout entier $j\in [1,n]$, on a aussit\^{o}t   $K_j\subseteq N_j$.  Notamment, pour tout entier  $j\in [1,n]$,  ${\al_j}^{p^{-1}}\in N_{j-1}({\be_{j}}^{p^{-1}},T_{1}^{j})$ (\`{a} savoir $T_{1}^j={\al_j}^{p^{-1}} T_{2}^{j-1}+{\be_j}^{p^{-1}}$).
D'autre part, comme $k=k_0(X,\al_1,\cdots,$ $\al_n,\be_1\cdots,\be_n)$ est un corps des fractions rationnelles aux ind\'{e}termin\'{e}es $(X,\al_1,\cdots,\al_n,\be_1\cdots,\be_n)$, alors pour tous entiers $1\leq s\leq j\leq n$, $di(N_s({\al_{s+1}}^{p^{-1}},$ $\cdots, {\al_j}^{p^{-1}},{\be_{s+1}}^{p^{-1}},\dots, {\be_{j}}^{p^{-1}})/N_s)=2(j-s)$. En outre $di(N_j/N_s)=2(j-s)$. \sk

\pre \smartqed
Soit maintenant  $m$ le plus petit entier tel que ${\al_s}^{p^{-1}} \in K_m$. Il clair que  $m>1$, et on distingue deux cas :
\sk

1-ier cas : si $m\leq s$, donc ${\al_s}^{p^{-1}}\in K_m \subseteq  K_s \subseteq N_{s-1}(K_s) $. Or, on a, $T_{1}^s={\al_s}^{p^{-1}}T_{2}^{s-1}+{\be_{s}}^{p^{-1}}$, on en d\'{e}duit que ${\be_{s}}^{p^{-1}} \in N_{s-1}(K_s)$. Cela conduit en vertu du th\'{e}or\`{e}me  \textcolor{blue}{\ref{thm6}} \`{a} $2=di(N_{s-1}( {\al_{s}}^{p^{-1}},{\be_{s}}^{p^{-1}})/N_{s-1})\leq
di(N_{s-1}( T_s(X))/N_{s-1})=1$,  c'est une contradiction.
\sk

2-i\`{e}me cas : $m\geq s+1$. On a ${\al_s}^{p^{-1}} \not\in K_{m-1}$ et ${\al_s}^{p^{-1}} \in K_m$, comme $K_m/K_{m-1}$ est $q$-simple et $T_{1}^{m}\not\in K_{m-1}$, (sinon on aura ${(\al_{1}^m)}^{p^{-1}}\in N_{m-1}({(\be_{1}^{m})}^{p^{-1}})$, absurde),  alors $K_{m-1} (T_{1}^{m})=K_{m-1} ({\al_s}^{p^{-1}})$. Toutefois, on a $s\leq m-1$, donc ${\al_s}^{p^{-1}}\in N_{s}\subseteq N_{m-1}$, et par suite $K_{m-1} ({\al_s}^{p^{-1}})
\subseteq N_{m-1}$, (car $K_{m-1}\subseteq N_{m-1}$). En particulier, $T_{1}^{m}\in N_{m-1}$  et donc ${(\al_{1}^{m})}^{p^{-1}}\in N_{m-1}({(\be_{1}^{m})}^{p^{-1}})$, (car $T_{1}^m={\al_m}^{p^{-1}}T_{2}^{m-1}+{\be_{m}}^{p^{-1}}$), c'est une contradiction. \cqfd
\sk

\noindent\noindent{\bf Preuve du th\'{e}or\`{e}me \textcolor{blue}{\ref{t420}}} \smartqed Soit $m=lm(K_r/k)$. Si
$K_{r-1}\not\subseteq m$, soient $s$ le plus petit entier tel que
$k(T_s(X))\not\subset m$ ($s<r$) et $j$ le premier entier tel que
$T_{j}^s\not\in m$. Comme, $T_{1}^{s+1}={\al_{s+1}}^{p^{-1}}T_{2}^{s}+
{\be_{s+1}}^{p^{-1}}$, et pour tout entier $j\geq 2$,
\begin{eqnarray*}
T_{j}^{s+1}&=&{\al_{s+1}}^{p^{-1}}T_{2j}^{s}+
{(T_{j-1}^{s+1})}^{p^{-1}},\\
&=& {\al_{s+1}}^{p^{-1}}{T_{2j}^s}+{\al_{s+1}}^{p^{-2}}
{(T_{2(j-1)}^s)}^{p^{-1}}+\cdots+{\al_{s+1}}^{p^{-j}}{(T_{2}^s)}^{p^{-j+1}}+ {\be_{s+1}}^{p^{-j}},
\end{eqnarray*}
alors ${(T_{1}^{s+1})}^p={\al_{s+1}}{(T_{2}^{s})}^p+
{\be_{s+1}}$, et pour tout entier $j\geq 2$,
${(T_{j}^{s+1})}^{p^j}={\al_{s+1}}^{p^{j-1}}{(T_{2j}^s)}^{p^j}$ $+{\al_{s+1}}^{p^{j-2}}
({{T_{2(j-1)}^s})}^{p^{j-1}}+\cdots+{\al_{s+1}}{(T_{2}^s)}^p $
$+{\be_{s+1}}$. Puisque $K_{s-1}({(T_{2j}^s)}^{p^j})=K_{s-1} (T_{j}^s)$ $\not\subseteq m$ et $K_{s-1}\subseteq m$, on en d\'{e}duit que ${(T_{2j}^s)}^{p^j}\not\in m$. D'autre part, il est trivialement \'{e}vident que  ${\al_{s+1}}^{p^{j-2}}({
T_{2(j-1)}^s)}^{p^{j-1}}+\cdots+{\al_{s+1}}{(T_{2}^s)}^p
+{\be_{s+1}}\in K_{s-1}({T_{j-1}^s})\subseteq m$. D'apr\`{e}s le
lemme \textcolor{blue}{\ref{l37}}, on a ${\al_{s+1}}^{p^{j-1}}\in
m\cap {K_r}^{p^j}$. D'o\`u $\al_{s+1}^{p^{-1}}\in K_r\subseteq K_n$, contradiction. \cqfd
\sk

D'apr\`{e}s le corollaire \textcolor{blue}{\ref{c319}}, on v\'{e}rifie
aussit\^{o}t que $K_{s-1}=rp(lm(K_s/k)/k)=rp(K_{s-1}/k)$
$=lqm(K/k)$ pour tout entier $s \in [1,n]$. Plus particuli\`{e}rement, on a :
\begin{thm} {\label{t421}} Avec les notations ci-dessus, pour tout entier $s\in [1,n]$, $K_s/k$ est une extension
$e$-ferm\'{e}e d'ordre $s$.
\end{thm}
\pre \smartqed Raisonnons par r\'{e}currence sur $s$.  Le cas $s=1$
\'{e}tant trivial, supposons que $2\leq s$, et soit $k\subseteq
F_1\subseteq\cdots\subseteq F_j=K_s$ une suite de d\'{e}composition
dont les termes sont relativement parfaits. Compte tenu de la
proposition \textcolor{blue}{\ref{p417}} et du th\'{e}or\`{e}me
\textcolor{blue}{\ref{t420}}, on a $F_{j-1}=K_{s-1}$. D'apr\`{e}s la
propri\'{e}t\'{e} de r\'{e}currence appliqu\'{e}e \`{a} $K_{s-1}/k$,
on obtient $j-1=s-1$, ou encore $j=s$. \cqfd
\sk

Dans la section qui suit nous examinons  plus particuli\`{e}rement une suite
de d\'{e}composition  dite suite $l$-$i$-modulaire. Il s'agit d'une d\'{e}composition en extensions $lq$-modulaires d\'{e}finie par cha\^{i}nage inverse.
\section{Extensions $l$-$i$-modulaires}
\begin{pro} Toute extension $q$-finie $K/k$ se d\'{e}compose sous forme
$K=m_0(K/k)\supseteq m_1(K/k)\supseteq \cdots \supseteq m_{i-1}(K/k)\supseteq m_i(K/k)=k$
o\`u
$m_{j+1}(K/k)=lqm(m_{j}(K/k)/k)$ pour tout entier $j \in [0,i[$.
\end{pro}
\pre R\'{e}sulte imm\'{e}diatement des propositions \textcolor{blue}{\ref{pr41}}, \textcolor{blue}{\ref{p318}} et du corollaire \textcolor{blue}{\ref{corrr1}}.\cqfd \sk

Il est \'{e}videment clair que $k=m_i(K/k)\subset m_{i-1}(K/k)\subseteq\cdots\subseteq m_0(K/k)=K$ est une suite de d\'{e}composition associ\'{e}e \`{a} $K/k$. En outre, cela  nous permet d'adopter la d\'{e}finition suivante :
\begin{df}
Toute extension $q$-finie $K/k$ qui se d\'{e}compose sous forme de $K=m_0(K/k)\supseteq m_1(K/k)\supseteq \cdots \supseteq m_{i-1}(K/k)\supseteq m_i(K/k)=k$ o\`u
$m_{j+1}(K/k)=lqm(m_{j}(K/k)/k)$ pour tout entier $j \in [0,i[$ sera appel\'{e}e extension $l$-$i$-modulaire. En particulier, le plus petit entier $i$ pour lequel
$K/k$ est $l$-$i$-modulaire s'appelle le degr\'{e} de modularit\'{e} inf\'{e}rieur de $K/k$, et se note $dmi(K/k)$.
\end{df}
\begin{rem} L'invariant $dmi(K/k)$  permet de mesurer le niveau  de modularit\'{e} inf\'{e}rieur de $K/k$.
\end{rem}

On montre sans peine  que :
\sk

\begin{itemize}{\it
\item $dmi(K/k)=inf\{i\in\N$ tel que $m_i(K/k)=k\}$.
\item La $l$-$i$ modularit\'{e} entra\^{\i}ne la $i$-modularit\'{e}. En particulier,
$dm(K/k)\leq dmi(K$ $/k)$.
\item Toute extension $lq$-modulaire est $l$-$i$-modulaire pour tout entier $1\leq i$.
\item $K/k$ est $l$-$j$-modulaire pour tout entier $j\geq dmi(K/k)$.
\item Si $K/k$ est $e$-ferm\'{e}e d'ordre $n$, on a $H_{n-i}(K/k)=m_i(K/k)$ pour tout entier $i\in [0, n]$.
La r\'{e}ciproque est g\'{e}n\'{e}ralement fausse  comme l'indique l'exemple ci-dessous. En outre, $K/k$ est $e$-ferm\'{e}e d'ordre $n$ si et
seulement si  $dmi(K/k)=di(K/k)=n$.}
\end{itemize}
\begin{exe} Soient $k_0$ un corps parfait de caract\'{e}ristique
$p>0$ et $(X,Y,Z_1,Z_2, $ $\al_1,\al_2)$ une famille alg\'{e}briquement ind\'{e}pendante sur
$k_{0}$. Posons $k=k_0(X,Y,Z_1,Z_2,$ $\al_1,\al_2)$ et $K_{n}=k(X^{p^{-2n}},\theta
_{n})$ o\`{u} $\theta _{1}=Z_1^{p^{-1}}X^{p^{-2}}+Z_2^{p^{-1}}$, et pour tout entier $n\geq 2$, $\theta _{n}=Z_1^{p^{-1}}X^{p^{-2n}}+{(\te_{n-1})}^{p^{-1}}=Z_1^{p^{-1}}X^{p^{-2n}}+Z_{1}^{p^{-2}}X^{p^{-2n+1}}+
\cdots+Z_1^{p^{-n}}X^{p^{-n-1}}+Z_2^{p^{-n}}$. De m\^{e}me, on note $L_n= k(Y^{p^{-2n}},\Phi
_{n})$ o\`{u} $\Phi _{1}={\al_1}^{p^{-1}}Y^{p^{-2}}+{\al_2}^{p^{-1}}$, et pour tout entier $n\geq 2$,  $\Phi _{n}={\al_1}^{p^{-1}}X^{p^{-2n}}+{(\phi_{n-1})}^{p^{-1}}={\al_1}^{p^{-1}}Y^{p^{-2n}}
+{\al_{1}}^{p^{-2}}Y^{p^{-2n+1}}+ \cdots+{\al_1}^{p^{-n}}Y^{p^{-n-1}}+{\al_2}^{p^{-n}}.$
\end{exe}

Comme $K_n\subseteq K_{n+1}$ et $L_n\subseteq L_{n+1}$ pour tout $n\in \N^*$, alors $K=\displaystyle\bigcup_n K_{n}$ et $L=\displaystyle\bigcup_n L_{n}$ sont des corps commutatifs.
Posons ensuite  $S=K(L)$, $H_1=H(S/k)$, $H_2=H(S/K)$ et $H_3=H(S/L)$.
\sk

On v\'{e}rifie imm\'{e}diatement que :
\sk

\begin{itemize}{\it
\item $ k(X^{p^{-\infty}},Y^{p^{-\infty}})\subseteq H_1\subseteq H_2\cap H_3$,
puisque $k(X^{p^{-\infty}},Y^{p^{-\infty}})/k$,  $K(H_1)/K$ et $L(H_1)/L$ sont
$lq$-modulaires et relativement parfaites en vertu de la proposition \textcolor{blue}{\ref{pr8}}.
\item $K(Y^{p^{-\infty}})\subseteq H_2$ et $L(X^{p^{-\infty}})\subseteq H_3$,
du fait que $K(Y^{p^{-\infty}})/K$ et $ L(X^{p^{-\infty}})/L$ sont $lq$-modulaires et relativement parfaites.
\item $lqm(S/k)=k(X^{p^{-\infty}}, Y^{p^{-\infty}})$, on utilise les m\^{e}mes techniques
de raisonnement que ceux de l'exemple \textcolor{blue}{\ref{e1}}. En particulier,
$k\subseteq k(X^{p^{-\infty}}, Y^{p^{-\infty}})\subseteq S$ est une suite $l$-$i$-modulaire.
\item $K(Y^{p^{-\infty}})\cap L(X^{p^{-\infty}})=k(X^{p^{-\infty}},
Y^{p^{-\infty}})$, car $$S\sim L\otimes_k K\sim
L(X^{p^{-\infty}})\otimes_{k(X^{p^{-\infty}},Y^{p^{-\infty}})} K(Y^{p^{-\infty}}).$$}
\end{itemize}
\begin{pro} Sous les notations de l'exemple ci-dessus, $k\subseteq k(X^{p^{-\infty}},
Y^{p^{-\infty}})\subseteq S$ est une suite $l$-$i$-modulaire, avec $H_1=k(X^{p^{-\infty}},
Y^{p^{-\infty}})=lqm(S/k)$. En particulier, $S/k$ n'est pas $e$-ferm\'{e}e.
\end{pro}
\pre \smartqed Puisque $k\subseteq k(X^{p^{-\infty}}, Y^{p^{-\infty}})\subseteq S$ est une suite
$l$-$i$-modulaire, il suffit
donc de montrer que $H_1=k(X^{p^{-\infty}},Y^{p^{-\infty}})$. Pour cela on va prouver
d'abord que $H_2=K(Y^{p^{-\infty}})$. Comme $S/K(Y^{p^{-\infty}})$
est $q$-simple et $ H_2/K(Y^{p^{-\infty}})$ est relativement parfaite, on a
$H_2=K(Y^{p^{-\infty}})$ ou $H_2=S$. Or, $S/K$ n'est pas $lq$-modulaire,
puisque $lqm(S/K)=K(Y^{p^{-\infty}})$ (cf. corollaire \textcolor{blue}{\ref{c319}}),
il en r\'{e}sulte que $H_2=K(Y^{p^{-\infty}})$. De la m\^{e}me fa\c{c}on, on montre
que $H_3=L(X^{p^{-\infty}})$. Par suite, $k(X^{p^{-\infty}},Y^{p^{-\infty}})\subseteq
H_1\subseteq H_2\cap H_3=K(Y^{p^{-\infty}})\cap L(X^{p^{-\infty}})=
k(X^{p^{-\infty}},Y^{p^{-\infty}})$. D'o\`{u} $H_1=k(X^{p^{-\infty}},
Y^{p^{-\infty}})=lqm(S/k)$. \cqfd
\sk

Voici une g\'{e}n\'{e}ralisation naturelle du corollaire \textcolor{blue}{\ref{c319}}.
\begin{pro} {\label{p319}}
Soient $K$ et $K'$ deux corps interm\'{e}diaires d'une extension
$q$-finie $K/k$. Alors :
\sk

\begin{itemize}{\it
\item[{\rm (i)}] Si $K\sim K'$,  pour tout entier naturel  $j$, on a $m_j(K/k)=m_j(K'/k)$.
\item[{\rm (ii)}] Pour toute sous-extension finie $k'/k$ de $K/k$, pour tout entier naturel  $j$, on a $k'(m_j(K/k))=m_j(K/k')$.}
\end{itemize}
\end{pro}
\pre \smartqed D'apr\`{e}s le corollaire \textcolor{blue}{\ref{c319}}, il suffit de remarquer que si
$K=m_0\supseteq m_1\supseteq\dots\supseteq m_{i-1}\supseteq m_i=k$
est une suite $l$-$i$-modulaire  associ\'{e}e \`{a} $K/k$, alors
$K=k'(K)=k'(m_0)\supseteq k'(m_1)\supseteq\dots\supseteq k'(m_{i-1})\supseteq k'(m_i)=k'$ et
$K'\supseteq m_1\supseteq\dots\supseteq m_{i-1}\supseteq m_i=k$ sont aussi deux suites
$l$-$i$-modulaires associ\'{e}es respectivement \`{a} $K/k'$ et $K'/k$.\cqfd
\sk

Comme cons\'{e}quence, la $l$-$i$-modularit\'{e} est respect\'{e}e
\`{a} une extension finie pr\`{e}s. Plus pr\'{e}cis\'{e}ment, on a :
\begin{pro} Soit $K/k$ une extension $q$-finie d'exposant non born\'{e}. Alors :
\sk

\begin{itemize}{\it
\item[{\rm (1)}] Si $k'\sim k$ et
$k'\subset K$, alors $K/k$ est $l$-$i$-modulaire si et seulement si il en est de
m\^{e}me de $K/k'$. En outre, $dmi(K/k)=dmi(K/k')$.
\item[{\rm (2)}] Si $K\sim K'$ et $k\subset
K'$, alors $K/k$ est $l$-$i$-modulaire si et seulement si $K'/k$ l'est aussi. Notamment, $dmi(K/k)=dmi(K'/k)$.
\item[{\rm (3)}] Si $k'\sim k$ et $K\sim K'$ avec $k'\subset
K'$, alors $K/k$ est $l$-$i$-modulaire si et seulement si il en est de
m\^{e}me de $K'/k'$. En outre, $dmi(K/k)=dmi(K'/k')$.}
\end{itemize}
\end{pro}
\pre (1) et (2) r\'{e}sultent  de la proposition ci-dessus, et  l'assertion $(3)$ d\'{e}coule
aussit\^{o}t de $(1)$ et $(2)$.\cqfd
\sk

Comme cons\'{e}quence, le r\'{e}sultat qui suit permet de ramener l'\'{e}tude de la $l$-$i$-modularit\'{e} au cas des
extensions relativement parfaites.
\begin{cor} Soit $K/k$ une extension $q$-finie d'exposant non born\'{e}. Alors :
\sk

\begin{itemize}{\it
\item[{\rm (i)}] $K/k$ est $l$-$i$-modulaire si et seulement si il en est de m\^{e}me de $rp(K/k)/k$.
En outre, $dmi(K/k)=dmi(rp(K/k)/k)$.
\item[{\rm (ii)}] Pour toute extension finie $L/k$, on a $K/k$ est $l$-$i$-modulaire si et seulement si
il en de m\^{e}me pour $L(K)/L$. En outre, $dmi(K/k)=dmi(L(K)/L)$.}
\end{itemize}
\end{cor}
\pre \smartqed Imm\'{e}diat, puisque $rp(K/k)\sim K$,  $L\sim k$ et $L(K)\sim K$.
 \cqfd
\sk

Comme dans le cas de la $i$-modularit\'{e}, le r\'{e}sultat suivant montre que le niveau  de
modularit\'{e} inf\'{e}rieur d'une extension $K/k$
est aussi  major\'{e} par l'entier $di(rp(K/k)/k)$.
\begin{thm} Pour toute  extension $q$-finie $K/k$ d'exposant non born\'{e},  on a $dmi(K/k)\leq di(rp(K/k)/k)$.
\end{thm}
\pre \smartqed  Il suffit de remarquer que $dmi(K/k)\leq dm(K/k)\leq di(rp(K/k)/k)$. \cqfd
\sk

Soient $L$ et $K$ deux corps interm\'{e}diaires d'une extension $q$-finie $\Omega/k$. On a d\'{e}j\`{a} vu (cf. proposition \textcolor{blue}{\ref{pr8}} et
proposition \textcolor{blue}{\ref{pr5}}) que :
\sk

\begin{itemize}{\it
\item[$\bullet$] Si $K/k$ et $L/k$ sont $l$-$1$-modulaires, il en est de m\^{e}me de $K(L)/k$.
\item[$\bullet$] $L(K)/L$ est $l$-$1$-modulaire si $K/k$ l'est.}
\end{itemize}
\sk

Plus g\'{e}n\'{e}ralement, on estime que la conjecture suivante est vraie.
\begin{con} Sous les notations ci-dessus, on a :
\sk

\begin{itemize}{\it
\item[{\rm (1)}] Le produit de deux extensions $l$-$i$-modulaires est $l$-$i$-modulaire. Autrement dit, si $K/k$ et
$L/k$ sont $l$-$i$-modulaires, il en est de m\^{e}me de $K(L)/k$.
\item[{\rm (2)}] La $l$-$i$-modularit\'{e} est respect\'{e}e si on change le
corps de base dans le sens ascendant. Par ailleurs,
$L(K)/L$ est $l$-$i$-modulaire si $K/k$ l'est.}
\end{itemize}
\end{con}

Si la conjecture ci-dessus est vraie, on en d\'{e}duit imm\'{e}diatement que :
\sk

\begin{itemize}{\it
\item[{\rm (1)}] $dmi(K(L)/k)\leq \sup(dmi(K/k), dmi(L/k))$.
\item[{\rm (2)}] $dmi(K(L)/L)\leq dmi(K/k)$.
\item[{\rm (3)}] Pour toute sous-extension $F/k$ de $K/k$, si $K/k$ est $l$-$i$-modulaire, il en est de m\^{e}me de $K/F$.}
\end{itemize}
\subsection{La $l$-$i$-modularit\'{e} et la stabilit\'{e}}

Soit $K/k$ une extension $q$-finie.
L'intersection portant sur $k$ ne respecte pas la
$l$-$i$-modularit\'{e} comme le montre l'exemple suivant :
\begin{exe} {\label{e5}}  On reprend  (\cite{Che-Fli5}, Exemple \textcolor{blue}{8}). Soient $k_0$
un corps parfait de caract\'{e}ristique $p>0$ et
$k=k_0(X,\al_1,\al_2,$ $\be_1,\be_2)$ le corps des fractions rationnelles aux ind\'{e}termin\'{e}es $(X,\al_1,\al_2,$ $\be_1,\be_2)$. Posons : \\
$T_0(X)=(T_{1}^0,\cdots, T_{n}^0,\cdots)$ o\`{u} $T_{n}^0=X^{p^{-n}}$.\\
$T_1(X)=( T_{1}^1,\cdots, T_{n}^1,\cdots)$ o\`{u}
$T_{1}^1={\al_1}^{p^{-1}}
T_{2}^0+{\be_1}^{p^{-1}}$ et $T_{n}^1={\al_1}^{p^{-1}}T_{2n}^{0}+{(T_{n-1}^1)}^{p^{-1}}.$\\
$T_2(X)=( T_{1}^2,\cdots, T_{n}^2,\cdots)$ o\`{u}
$T_{1}^2={\al_2}^{p^{-1}}T_{2} ^1+{\be_2}^{p^{-1}}$ et
$T_{n}^2={\al_2}^{p^{-1}}T_{2n}^{1}+{(T_{n-1}^2)}^{p^{-1}}.$
\end{exe}

Soient $K=k(T_0(X),T_1(X),T_2(X))$ et
$K_1=K({\al_1}^{p^{-\infty}})$. Puisque,
\begin{eqnarray*}
T_{n}^1&=&{\al_{1}}^{p^{-1}}T_{2n}^{0}+{(T_{n-1}^1)}^{p^{-1}},\\
&=& {\al_{1}}^{p^{-1}}T_{2n}^{0}+
{\al_{1}}^{p^{-2}}{(T_{2(n-1)}^{0})}^{p^{-1}}+\cdots+{\al_{1}}^{p^{-n}}{(T_{2}^{0})}^{p^{-n+1}}+
{\be_1}^{p^{-n}},
\end{eqnarray*}
et ${\al_1}^{p^{-n}}\in K_1$, alors ${\be_1}^{p^{-n}}\in K_1$ ; et par suite
$$K_1=k(X^{p^{-\infty}},{\al_{1}}^{p^{-\infty}},{\be_{1}}^{p^{-\infty}})(T_2(X)).$$
\begin{thm} {\label{t51}} Sous les hypoth\`{e}ses  ci-dessus,  les extensions $K_1/k({\al_1}^{p^{-\infty}},{\be_1}^{p^{-\infty}})$
et $K_1/k(T_0(X))$ sont $l$-$2$-modulaires, mais
$K_1/k({\al_1}^{p^{-\infty}},{\be_1}^{p^{-\infty}})\cap k(T_0(X))$
ne l'est pas.
\end{thm}

La d\'{e}monstration utilise plusieurs r\'{e}sultats pr\'{e}liminaires.
\begin{lem} {\label{l52}}
Sous les hypoth\`{e}ses cit\'{e}es en haut, on a :
\begin{itemize}{\it
\item[{\rm (i)}] $k( {(T_{i}^1)}^{p^{i}})= k( {(T_{i}^{2})}^{p^{3i}})= k( {(T_{i}^0)})$.
\item[{\rm (ii)}] ${\al_{2}}^{p^{-1}}\not\in K_1$ et ${\be_{2}}^{p^{-1}}\not\in K_1$.}
\end{itemize}
\end{lem}
\pre \smartqed Par construction, on a :
\begin{eqnarray*}
T_{i}^1 &=& {\al_{1}}^{p^{-1}}T_{2i}^0+(T_{i-1}^1)^{p^{-1}}, \\
        &=& {\al_{1}}^{p^{-1}}{T_{2i}^0}+{\al_{1}}^{p^{-2}}( {T_{2(i-1)}^0)}
^{p^{-1}}+\cdots+{\al_{1}}^{p^{-i}}{(T_{2}^0)}^{p^{-i+1}}+{\be_{1}}^{p^{-i}} ;
\end{eqnarray*}
on en d\'{e}duit que
${(T_{i}^1)}^{p^i}={\al_{1}}^{p^{i-1}}{(T_{2i}^0)}^{p^i}+{\al_{1}}^{p^{i-2}}{(T_{2(i-1)}^0)}
^{p^{i-1}}+\cdots+{\al_{1}}{(T_{2}^0)} ^{p}+\be_{1}$, et par suite  $k({(T_{2i}^0)}^{p^i},
{(T_{2(i-1)}^0)}^{p^{i-1}},\cdots,{(T_{2}^0)} ^{p})=k({(T_{i}^1)}^{p^i},
{(T_{2(i-1)}^0)}^{p^{i-1}},\cdots,{(T_{2}^0)}^{p})$. Comme $k({(T_{2(i-1)}^0)}^{p^{i-1}},\cdots,{(T_{2}^0)} ^{p})=k(T_{i-1}^0)=k(X^{p^{-i+1}})\varsubsetneq k(X^{p^{-i}})=k(T_{i}^{0})$, alors
$k({(T_{i}^{1})}^{p^i})=k(T_{i}^{0})$. Egalement,
${(T_{i}^2)}^{p^{2i}}={\al_{1}}^{p^{2i-1}}{(T_{2i}^1)}^{p^{2i}}$ $+{\al_{1}}^{p^{2i-2}}{(T_{2(i-1)}^1)}
^{p^{2i-1}}$ $+\cdots+{\al_{1}^{p^i}}{(T_{2}^1)}
^{p^{i+1}}+{\be_{1}}^{p^i}$. Comme $k({(T_{2(i-j)}^1)}^{p^{2i-j}})\subseteq
k(({T_{2(i-j)}^1)}^{p^{2(i-j)}})=k({T_{2(i-j)}^0})$ pour tout entier $j \in [1, i-1]$, on aura $
k({(T_{i}^{2})}^{p^{2i}})=k({(T_{2i}^{1})}^{p^{2i}})=k(T_{2i}^0)$ ;
et par suite, $
k({(T_{i}^{2})}^{p^{3i}})=k({(T_{2i}^0)}^{p^i})=k(T_{i}^0)$. Si
${\al_2}^{p^{-1}}\in K_1$ ou ${\be_2}^{p^{-1}}\in K_1$, puisque
$T_{1}^{2}={\al_2}^{p^{-1}}T_{2}^{1}+{\be_2}^{p^{-1}}$, on aura
$k(X^{p^{-1}},{\al_1}^{p^{-1}},$ ${\be_1}^{p^{-1}},{\al_2}^{p^{-1}},{\be_2}^{p^{-1}})\subseteq
K_1$, et par suite $5=di(k(X^{p^{-1}},{\al_1}^{p^{-1}},$
${\be_1}^{p^{-1}},{\al_2}^{p^{-1}},
 {\be_2}^{p^{-1}})/k)\leq di(K_1/k)=4$, contradiction. \cqfd
\begin{pro}
Sous les m\^{e}mes hypoth\`{e}ses de l'exemple ci-dessus, on a
$k(T_0(X),$ $T_1(X))$ est la plus petite sous-extension de $K_1/k$
telle que $K_1/k(T_0(X),$ $T_1(X))$ est modulaire (i.e. $lm(K_1/k)=k(T_0(X),$ $T_1(X))$).
\end{pro}
\pre \smartqed Il est clair que $K_1/k(T_0(X),T_1(X))$ est
modulaire. Soit $m=lm(K_1/k)$, donc $m\subseteq k(T_0(X),T_1(X))$.
Si $k(T_0(X))\not \subseteq m$, soit $j$ le plus petit entier tel
que $T_{j}^0\not\in m$. Par construction,  ${(T_{1}^2)}^{p^{3}}={\al_2}^{p^{2}}{(T_{2}^{1})}^{p^{3}}+{\be_2}^{p^{2}}$ et
${(T_{j}^2)}^{p^{3j}}={\al_2}^{p^{3j-1}}{(T_{2j}^{1})}^{p^{3j}}+{(T_{j-1}^2)}^{p^{3j-1}}$ si $j\not =1$, de plus en vertu du
lemme \textcolor{blue}{\ref{l52}},
$k({(T_{2j}^1)}^{p^{3j}})=k(T_{j}^0)\not\subseteq m $ et
$k({(T_{j-1}^{2})}^{p^{3j-1}})\subseteq
k({(T_{j-1}^{2})}^{p^{3(j-1)}})=k(T_{j-1}^0)\subseteq m$. Comme $K_1/m$ est modulaire ;
d'apr\`{e}s le lemme \textcolor{blue}{\ref{l37}}, dans les deux cas
${\al_2}^{p^{-1}}\in K_1$, contradiction. D'o\`u
$k(T_0(X))\subseteq m$. De m\^{e}me, si $k(T_1(X))\not\subseteq m$,
soit $j$ le plus petit entier tel que $T_{j}^1\not\in m$. Egalement, ${(T_{1}^2)}^{p}={\al_2}{(T_{2}^{1})}^{p}+{\be_2}$ et
${(T_{j}^2)}^{p^{j}}={\al_2}^{p^{j-1}}{(T_{2j}^{1})}^{p^{j}}+{(T_{j-1}^2)}^{p^{j-1}}$ si $j\not=1$,
avec  $k(T_0(X))({(T_{2j}^1)}^{p^{j}})=k(T_0(X))({(T_{j}^{1})}) \not \subseteq m$
et $k(T_0(X))({(T_{j-1}^{2})}^{p^{j-1}})\subseteq
k(T_0(X))({(T_{j-1}^{1})})\subseteq m$. D'apr\`{e}s le lemme
\textcolor{blue}{\ref{l37}}, dans les deux cas ${\al_2}^{p^{-1}}\in K_1$, contradiction. Il en r\'{e}sulte que $m=k(T_0(X),T_1(X))$.\cqfd
\begin{pro} Sous les notations de l'exemple ci-dessus, on a  $k(T_0(X),{\al_{1}}^{p^{-\infty}},
$ ${\be_{1}}^{p^{-\infty}})/k({\al_{1}}^{p^{-\infty}},{\be_{1}}^{p^{-\infty}})$ (respectivement,
$k(T_0(X),T_1(X))/k(T_0(X)$) est la plus petite sous-extension de
$K_1/k({\al_{1}}^{p^{-\infty}},{\be_{1}}^{p^{-\infty}})$ (respectivement,  de
$K_1/k(T_0(X))$) telle que $K_1/k(T_0(X),
{\al_{1}}^{p^{-\infty}},{\be_{1}}^{p^{-\infty}})$ (respectivement,
$K_1/k(T_0(X), T_1(X))$) est modulaire.
\end{pro}
\pre \smartqed Cela r\'{e}sulte imm\'{e}diatement du corollaire
\textcolor{blue}{\ref{c319}} et du fait que  $K_1/k(T_0(X),$ $ {\al_{1}}^{p^{-\infty}},{\be_{1}}^{p^{-\infty}})$ et
$K_1/k(T_0(X), T_1(X))$ sont modulaires.
 \cqfd
\sk

D'apr\`{e}s la proposition \textcolor{blue}{\ref{p318}}, on a aussit\^{o}t :
\sk

\begin{itemize}{\it
\item[$\bullet$] $lqm(K_1/k({\al_{1}}^{p^{-\infty}},{\be_{1}}^{p^{-\infty}}))=
rp(lm(K_1/k({\al_{1}}^{p^{-\infty}},{\be_{1}}^{p^{-\infty}}))/
k({\al_{1}}^{p^{-\infty}},{\be_{1}}^{p^{-\infty}})) $ $=rp(k(T_0(X),{\al_{1}}^{p^{-\infty}}, {\be_{1}}^{p^{-\infty}})$
$/k({\al_{1}}^{p^{-\infty}},{\be_{1}}^ {p^{-\infty}}))=
k(T_0(X),{\al_{1}}^{p^{-\infty}}, {\be_{1}}^{p^{-\infty}})$ et
$lqm(K_1/k(T_0(X)))=rp(lm(K_1/k(T_0(X)))/k(T_0(X)))=k(T_0(X),T_1(X))$.
\item[$\bullet$] $k(T_0(X))=lqm(k(T_1(X),T_0(X))/k)$, d\'{e}monstration identique \`{a} celle du
th\'{e}or\`{e}me \textcolor{blue}{\ref{t420}}.}
\end{itemize}
\sk

\noindent\noindent {\bf Preuve du th\'{e}or\`{e}me \textcolor{blue}{\ref{t51}}.} \smartqed
Par construction, on a
$k({\al_1}^{p^{-\infty}},{\be_1}^{p^{-\infty}})\subseteq k(T_0(X),$
${\al_1}^{p^{-\infty}},{\be_1}^{p^{-\infty}})\subseteq K_1$
et $k(T_0(X))\subseteq k(T_0(X))(T_1(X))\subseteq K_1$ sont deux
suites de d\'{e}composition associ\'{e}es respectivement \`{a}
$K_1/k({\al_1}^{p^{-\infty}}, {\be_1}^{p^{-\infty}})$ et $K_1/k($
$T_0(X))$, avec $k(T_0(X),{\al_1}^{p^{-\infty}},
{\be_1}^{p^{-\infty}})=lqm(K_1/ k({\al_1}^{p^{-\infty}},
{\be_1}^{p^{-\infty}}))$ et $k(T_0(X),T_1(X))=lqm$ $(K_1/k($
$T_0(X)))$ ; donc
$K_1/k({\al_1}^{p^{-\infty}},{\be_1}^{p^{-\infty}})$ et
$K_1/k(T_0(X))$ sont $l$-$2$-modulaires. D'autre part, comme $k({\al_1}^{p^{-\infty}},{\be_1}^{p^{-\infty}},T_0(X))\sim
k({\al_1}^{p^{-\infty}},{\be_1}^{p^{-\infty}})\otimes_k k(T_0(X))
$, on a $k=
k({\al_1}^{p^{-\infty}},$ ${\be_1}^{p^{-\infty}})\cap k(T_0(X))$. De plus,
$k\subseteq k(T_0(X))\subseteq k(T_0(X),T_1(X))$ $\subseteq
K_1$ est une suite de d\'{e}composition associ\'{e}e \`{a} $K_1/k$, avec
$k(T_0(X),T_1(X))=lqm(K_1/k(T_0(X)))$ et $k(T_0(X))$
$=lqm(k(T_0(X),T_1(X))/k)$, on en d\'{e}duit que $K_1/k$ est
$l$-$3$-modulaire. \cqfd
\begin{rem} L'exemple pr\'{e}c\'{e}dent est aussi bon pour montrer que
$dm(K_1/k)\not =dmi(K_1/k)$. Plus pr\'{e}cis\'{e}ment, on a $2=
dm(K_1/k)<dmi(K_1/k)=3$. En particulier, la $i$-modularit\'{e} n'entra\^{\i}ne pas
la $l$-$i$-modularit\'{e}. Donc la $i$-modularit\'{e} est moins fine que
la $l$-$i$-modularit\'{e}.
\end{rem}


\begin{thebibliography}{999}
\bibitem[\textcolor{blue}{1}]{Che-Fli5} {M. Chellali and E. Fliouet}, {\textit{ Extensions $i$-Modulaires}}, International Journal of Algebra, Vol. {\bf 6}, (2012), no. 10, 457--492
\bibitem[\textcolor{blue}{2}]{Che-Fli4}M. Chellali et E. Fliouet, {\textit Sur la tour des cl\^{o}tures modulaires,} An. St. Univ. Ovidius Constanta Vol. {\bf 14(1)}, (2006), 45-66
\bibitem[\textcolor{blue}{3}]{Che-Fli} M. Chellali et E. Fliouet, {\textit{Sur les extensions purement
ins\'{e}parable,}} Arch. Math.  Vol {\bf 81}, (2003), 369-382
\bibitem[\textcolor{blue}{4}]{Che-Fli1}M. Chellali et E. Fliouet, {{\textit Extension presque modulaire, }}
Ann. Sci. Math Qu\'{e}bec Vol {\bf 28} no. 1-2, (2004), 65-75
\bibitem[\textcolor{blue}{5}]{Che-Fli2}M. Chellali et E. Fliouet, {\textit Extensions purement
ins\'{e}parables d'exposant non born\'{e},}  Archivum Mathematicum
{\bf 40}, (2004), 129-159.
\bibitem[\textcolor{blue}{6}]{Che-Fli7} M. Chellali et E. Fliouet, {\textit{Extensions semi-simples,}} An. Univ. Vest Timis., Ser. Mat.-Inform. 44, No. 2, 51-92 (2006). Vol {\bf 44}, no 2, (2006), 51--92
\bibitem[\textcolor{blue}{7}]{Che-Fli3}M. Chellali et E. Fliouet, {\textit{ Th\'{e}or\`{e}me
de la cl\^{o}ture $lq$-modulaire et applications,}} Colloq. Math. {\bf  122},
(2011), 275-287
\bibitem[\textcolor{blue}{8}]{Dev1}J.K.  Deveney, {\textit An intermediate theory for a purely inseparable Galois  theory,} Trans. Amer. Math. Soc.   {\bf 198}, (1975), 287-295
\bibitem[\textcolor{blue}{9}]{Dev2}J.K. Deveny, {\textit{$w_0$-generated field extensions,}} Arch. Math. {\bf  47}, (1986), 410-412
\bibitem[\textcolor{blue}{10}]{Che1} {E. Fliouet}, {\textit{Absolutely $lq$-finite  extensions}}, https://arxiv.org/pdf/1701.05430.pdf
\bibitem[\textcolor{blue}{11}]{Kim}L.A. Kime, {\emph Purely inseparable modular extensions of
unbounded exponent,}  Trans. Amer. Math. Soc  {\bf  176},  (1973) ,
335-349
\bibitem[\textcolor{blue}{12}]{Pic}G. Pickert, {\textit Inseparable
K\"{o}rperweiterungen,} Math. Z.   {\bf  52}, (1949),  81-135
\bibitem[\textcolor{blue}{13}]{Swe}M.E.  Sweedler,
{\textit Structure  of inseparable extensions},  Ann. Math.   {\bf  87}
(2),  (1968), 401-410
\bibitem[\textcolor{blue}{14}]{Wat}W.C. Waterhouse, {\textit{ The structure of inseparable field
extensions,}}  Trans. Am. Math. Soc.  {\bf  211},  (1975), 39-56
\bibitem[\textcolor{blue}{15}]{N.B} {N. Bourbaki}, {\textit{El\'{e}ments de Math\'{e}matique Th\'{e}orie des ensembles}}, Springer-Verlag Berlin Heidelberg 2006
\bibitem[\textcolor{blue}{16}]{N.B2} N. Bourbaki, {\textit{Alg\`{e}bre, Chapitre 1 \`{a} 3,}} Springer-Verlag Berlin Heidelberg 2007
\bibitem[\textcolor{blue}{17}]{Mor-Vin}J.N.  Mordeson and B.Vinograde,  {\textit Structure of arbitrary purely inseparable extension fields}, Springer-Verlag, Berlin, LNM {\bf  173},  (1970)
\end{thebibliography}
\end{document}